\documentclass[reqno,10pt,centertags,draft]{amsart}
\usepackage{amsmath,amsthm,amscd,amssymb,latexsym,upref}
\newcommand{\bbN}{{\mathbb{N}}}
\newcommand{\bbR}{{\mathbb{R}}}

\newcommand{\bbZ}{{\mathbb{Z}}}
\newcommand{\bbC}{{\mathbb{C}}}

\newcommand{\cF}{{\mathcal F}}

\newcommand{\cH}{{\mathcal H}}

\newcommand{\cM}{{\mathcal M}}

\newcommand{\cU}{{\mathcal U}}

\newcommand{\dott}{\,\cdot\,}
\newcommand{\no}{\notag}
\newcommand{\lb}{\label}
\newcommand{\f}{\frac}

\newcommand{\ol}{\overline}
\newcommand{\ti}{\tilde}
\newcommand{\wti}{\widetilde}

\newcommand{\Oh}{O}
\newcommand{\oh}{o}

\newcommand{\Arc}{\text{\rm{Arc}}}
\newcommand{\dom}{\text{\rm{dom}}}
\newcommand{\ess}{\text{\rm{ess}}}

\newcommand{\supp}{\text{\rm{supp}}}

\newcommand{\bi}{\bibitem}
\newcommand{\hatt}{\widehat}
\newcommand{\beq}{\begin{equation}}
\newcommand{\eeq}{\end{equation}}
\newcommand{\ba}{\begin{align}}
\newcommand{\ea}{\end{align}}

\newcommand{\abs}[1]{\lvert#1\rvert}

\renewcommand{\Re}{\text{\rm Re}}
\renewcommand{\Im}{\text{\rm Im}}
\renewcommand{\ln}{\text{\rm ln}}
\newcommand{\essran}{\text{\rm ess.ran}}

\newcommand{\norm}[1]{\left\Vert#1\right\Vert}
\newcommand{\Om}{\Omega}
\newcommand{\si}{\sigma}

\newcommand{\la}{\lambda}
\newcommand{\al}{\alpha}

\newcommand{\de}{\delta}
\newcommand{\te}{\theta}
\newcommand{\ze}{\zeta}

\newcommand{\C}{\mathbb{C}}
\newcommand{\Cz}{\C\backslash\{0\}}

\newcommand{\D}{\mathbb{D}}
\newcommand{\dD}{{\partial\hspace*{.2mm}\mathbb{D}}}
\newcommand{\Z}{\mathbb{Z}}

\newcommand{\UU}{\mathbb{U}}
\newcommand{\Cl}{\mathbb{C}_{\ell}}
\newcommand{\Cr}{\mathbb{C}_{r}}

\newcommand{\ph}[1]{\phantom{#1}}

\newcommand{\spn}{{\text{\rm span}}}

\newcommand{\ltz}{{\ell^2(\Z)}}
\newcommand{\Lt}[1]{{L^2(\dD;d\mu_{#1}(\cdot,k_0))}}
\newcommand{\st}{\,|\,}

\allowdisplaybreaks \numberwithin{equation}{section}
\newtheorem{theorem}{Theorem}[section]

\newtheorem{lemma}[theorem]{Lemma}
\newtheorem{corollary}[theorem]{Corollary}
\newtheorem{hypothesis}[theorem]{Hypothesis}
\theoremstyle{definition}
\newtheorem{definition}[theorem]{Definition}
\newtheorem{remark}[theorem]{Remark}

\begin{document}

\title[Weyl--Titchmarsh Theory for CMV Operators]{Weyl--Titchmarsh
Theory for CMV Operators \\ Associated with Orthogonal Polynomials \\ on
the Unit Circle}
\author[F.\ Gesztesy, and M.\ Zinchenko]{Fritz
Gesztesy and Maxim Zinchenko}
\address{Department of Mathematics,
University of Missouri, Columbia, MO 65211, USA}
\email{fritz@math.missouri.edu}
\urladdr{http://www.math.missouri.edu/personnel/faculty/gesztesyf.html}
\address{Department of Mathematics,
University of Missouri, Columbia, MO 65211, USA}
\email{maxim@math.missouri.edu}
\dedicatory{Dedicated with great pleasure to Barry Simon on the
occasion of his 60th birthday.}

\thanks{Based upon work supported by the National Science
Foundation under Grant No.\ DMS-0405526.}
%\date{November 15, 2004}
\thanks{\it J.\ Approx.\ Theory {\bf 139}, 172--213 (2006).}
%%%%%%%%%%%%%%%%%%%%%%%%%%%%%%%%%%%%%%%%%%%%%%%%%%%%%%%%
\begin{abstract}
We provide a detailed treatment of Weyl--Titchmarsh theory for
half-lattice and full-lattice CMV operators and discuss their systems
of orthonormal Laurent polynomials on the unit circle, spectral
functions, variants of Weyl--Titchmarsh functions, and Green's functions.
In particular, we discuss the corresponding spectral representations of
half-lattice and full-lattice CMV operators.
\end{abstract}
%%%%%%%%%%%%%%%%%%%%%%%%%%%%%%%%%%%%%%%%%%%%%%%%%%%%%%%%

\maketitle

%%%%%%%%%%%%%%%%%%%%%%%%%%%%%%%%%%%%%%%%%%%%%%%%%%%%%%%%%
\section{Introduction}\lb{s1}
%%%%%%%%%%%%%%%%%%%%%%%%%%%%%%%%%%%%%%%%%%%%%%%%%%%%%%%%%

The aim of this paper is to develop Weyl--Titchmarsh theory for
a special class of unitary doubly infinite five-diagonal
matrices. The corresponding unitary semi-infinite five-diagonal matrices
were recently introduced by Cantero, Moral, and
Vel\'azquez (CMV) \cite{CMV03} in 2003. In \cite[Sects.\ 4.5,
10.5]{Si04}, Simon introduced the corresponding notion of unitary doubly
infinite five-diagonal matrices and coined the term ``extended'' CMV
matrices. To simplify notations we will often just speak of CMV operators
whether or not they are half-lattice or full-lattice operators indexed by
$\bbN$ or $\bbZ$, respectively.

CMV operators on $\bbZ$ are intimately related to a completely
integrable version of the defocusing nonlinear Schr\"odinger equation
(continuous in time but discrete in space), a special case of the
Ablowitz--Ladik system. Relevant references in this context are, for
instance, \cite{AL75}, \cite{APT04}, \cite{GGH05}, \cite{GH05},
\cite{MEKL95}--\cite{Ne05a}, and the literature cited therein.
A recent application to a Borg-type theorem (an inverse spectral
result), which motivated us to write this paper, appeared in
\cite{GZ05}. For more details we refer to Theorem \ref{t1.1} at the end
of this introduction.

We denote by $\D$ the open unit disk in $\bbC$ and let $\alpha$ be a
sequence of complex numbers in $\D$, $\alpha=\{\al_k\}_{k \in \Z}
\subset \D$. The unitary CMV operator $U$ on $\ell^2(\bbZ)$ then can
be written as a special five-diagonal doubly infinite matrix in the
standard basis of $\ell^2(\bbZ)$ according to
\cite[Sects.\ 4.5, 10.5]{Si04}) as
\begin{align}
U = \begin{pmatrix} \ddots &&\hspace*{-8mm}\ddots
&\hspace*{-10mm}\ddots &\hspace*{-12mm}\ddots
&\hspace*{-14mm}\ddots &&&
\raisebox{-3mm}[0mm][0mm]{\hspace*{-6mm}{\Huge $0$}}
\\
&0& -\al_{0}\rho_{-1} & -\ol{\al_{-1}}\al_{0} &
-\al_{1}\rho_{0} & \rho_{0}\rho_{1}
\\
&& \rho_{-1}\rho_{0} &\ol{\al_{-1}}\rho_{0} &
-\ol{\al_{0}}\al_{1} & \ol{\al_{0}}\rho_{1} & 0
\\
&&&0& -\al_{2}\rho_{1} & -\ol{\al_{1}}\al_{2} &
-\al_{3}\rho_{2} & \rho_{2}\rho_{3}
\\
&&\raisebox{-4mm}[0mm][0mm]{\hspace*{-6mm}{\Huge $0$}} &&
\rho_{1}\rho_{2} & \ol{\al_{1}}\rho_{2} & -\ol{\al_{2}}\al_{3}
& \ol{\al_{2}}\rho_{3}&0
\\
&&&&&\hspace*{-14mm}\ddots &\hspace*{-14mm}\ddots
&\hspace*{-14mm}\ddots &\hspace*{-8mm}\ddots &\ddots
\end{pmatrix}. \lb{1.1}
\end{align}
Here the sequence of positive real numbers
$\{\rho_k\}_{k\in\bbZ}$ is defined by
\begin{equation}
\rho_k = \sqrt{1-\abs{\al_k}^2}, \quad k\in\bbZ, \label{1.2}
\end{equation}
and terms of the form $-\ol{\alpha_k}\alpha_{k+1}$,
$k\in\Z$, represent the $k,k$-diagonal entries in the infinite
matrix \eqref{1.1}. For the corresponding half-lattice CMV
operators $U^{(s)}_{+,k_0}$, $s\in[0,2\pi)$ in
$\ell^2([k_0,\infty)\cap\bbZ)$ we refer to \eqref{2.30}.

The relevance of this unitary operator $U$ on
$\ell^2(\bbZ)$, more precisely, the relevance of the
corresponding half-lattice CMV operator $U_{+,0}$ in
$\ell^2(\bbN_0)$ (cf.\ \eqref{2.31}) is derived from its
intimate relationship with the trigonometric moment problem
and hence with finite measures on the unit circle $\dD$.
(Here $\bbN_0=\bbN\cup\{0\}$.)  Let
$\{\alpha_k\}_{k\in\bbN}\subset\D$  and define the transfer
matrix
\begin{equation}
S(\zeta,k)=\begin{pmatrix} \zeta & \alpha_k \\ \overline{\alpha_k}
\zeta &1
\end{pmatrix}, \quad \zeta\in\dD,\; k\in\bbN, \lb{1.3}
\end{equation}
with spectral parameter $\zeta\in\dD$. Consider the system of difference
equations
\begin{equation}
\begin{pmatrix}\varphi_+(\zeta,k) \\ \varphi^*_+(\zeta,k)\end{pmatrix}
= S(\zeta,k)\begin{pmatrix} \varphi_+(\zeta,k-1) \\
\varphi^*_+(\zeta,k-1) \end{pmatrix}, \quad \zeta\in\dD,\;
k\in\bbN  \lb{1.4}
\end{equation}
with initial condition
\begin{equation}
\begin{pmatrix} \varphi_+(\zeta,0)\\ \varphi^*_+(\zeta,0) \end{pmatrix}
=\begin{pmatrix} 1 \\ 1 \end{pmatrix}, \quad \zeta\in\dD. \lb{1.4a}
\end{equation}
Then $\varphi_+ (\dott,k)$ are monic
polynomials of degree $k$ and
\begin{equation}
\varphi^*_+ (\zeta,k)= \zeta^k
\ol{\varphi_+(1/\ol{\zeta},k)}, \quad \zeta\in\dD,\;
k\in\bbN_0, \lb{1.4b}
\end{equation}
the reversed ${}^*$-polynomial of $\varphi_+(\cdot,k)$, is at
most of degree $k$. These polynomials were first introduced by Szeg\H o in
the 1920's in his work on the asymptotic distribution of eigenvalues of
sections of Toeplitz forms \cite{Sz20}, \cite{Sz21}  (see
also \cite[Chs.\ 1--4]{GS84}, \cite[Ch.\ XI]{Sz78}). Szeg\H o's
point of departure was the trigonometric moment problem and hence the
theory of orthogonal polynomials on the unit circle: Given a probability
measure $d\sigma_+$  supported on an infinite set on the unit circle, find
monic polynomials of degree $k$ in $\zeta=e^{i\theta}$, $\theta\in
[0,2\pi]$, such that
\begin{equation}
\int_{0}^{2\pi} d\sigma_+(e^{i\theta}) \,
\overline{\varphi_+ (e^{i\theta},k)} \varphi_+ (e^{i\theta},k')
=\gamma_{k}^{-2} \delta_{k,k'}, \quad k,k'\in\bbN_0,
\lb{1.4c}
\end{equation}
where (cf.\ \eqref{1.2})
\begin{equation}
\gamma_k^2=\begin{cases} 1, & k=0, \\
\prod_{j=1}^k \rho_j^{-2}, & k\in\bbN.
\end{cases} \lb{1.4d}
\end{equation}
One then also infers
\begin{equation}
\int_{0}^{2\pi} d\sigma_+(e^{i\theta}) \,
\overline{\varphi^*_+ (e^{i\theta},k)} \varphi^*_+
(e^{i\theta},k') = \gamma_{k''}^{-2}, \quad k''=\max\{k,k'\},
\; k,k' \in\bbN_0 \lb{1.4e}
\end{equation}
and obtains that $\varphi_+(\cdot,k)$ is orthogonal to
$\{\zeta^j\}_{j=0,\dots,k-1}$ in $L^2(\dD;d\sigma_+)$ and
$\varphi^*_+(\cdot,k)$ is orthogonal to
$\{\zeta^j\}_{j=1,\dots,k}$ in $L^2(\dD;d\sigma_+)$. Additional
comments in this context will be provided in Remark \ref{r2.8}.
For a detailed account of the relationship of
$U_{+,0}$ with orthogonal polynomials on the unit circle we refer to the
monumental two-volume treatise by Simon \cite{Si04} (see also
\cite{Si04b} and \cite{Si05} for a description of some of the principal results
in \cite{Si04}) and the exhaustive bibliography therein. For classical
results on orthogonal polynomials on the unit circle we refer, for
instance, to \cite{Ak65}, \cite{Ge46}--\cite{Ge61}, \cite{GS84},
\cite{Kr45}, \cite{Sz20}--\cite{Ve36}. More recent
references relevant to the spectral theoretic content of this paper are
\cite{GJ96}--\cite{GT94}, \cite{GZ05}, \cite{GN01}, \cite{Lu04}, \cite{PY04},
\cite{Si04a}.

We note that $S(\zeta,k)$ in \eqref{1.3} is not the
transfer matrix that leads to the half-lattice CMV operator
$U_{+,0}$ in $\ell^2(\bbN_0)$ (cf.\ \eqref{2.30}). After a
suitable change of basis introduced by Cantero, Moral, and
Vel\'azquez \cite{CMV03}, the transfer matrix $S(\zeta,k)$
turns into $T(\zeta,k)$ as defined in \eqref{A.31}.

In Section \ref{s2} we provide an extensive treatment of
Weyl--Titchmarsh theory for half-lattice CMV operators
$U_{+,k_0}$ on $\ell([k_0,\infty)\cap\bbZ)$ and discuss
various systems of orthonormal Laurent polynomials on the
unit circle, the half-lattice spectral function of
$U_{+,k_0}$,  variants of half-lattice Weyl--Titchmarsh
functions, and the Green's function of $U_{+,k_0}$. In
particular, we discuss the spectral representation of
$U_{+,k_0}$. While many of these results can be found in
Simon's two-volume treatise \cite{Si04}, we survey some of this
material here from an operator theoretic point of view, starting
directly from the CMV operator. Section \ref{s3} then contains our new
results on Weyl--Titchmarsh theory for full-lattice CMV
operators $U$ on $\ell^2(\bbZ)$. Again we discuss systems
of orthonormal Laurent polynomials on the unit circle, the
$2\times 2$ matrix-valued spectral and Weyl--Titchmarsh
functions of $U$, its Green's matrix, and the spectral
representation of $U$. Finally, Appendix \ref{A} summarizes
basic facts on Caratheodory and Schur functions relevant to
this paper.

We conclude this introduction with citing a Borg-type (inverse spectral)
result from our paper \cite{GZ05}, which motivated us to write the present
paper.

First we introduce our notation for closed arcs on the unit circle $\dD$,
\begin{equation}
\Arc\big(\big[e^{i\theta_1},e^{i\theta_2}\big]\big)
=\big\{e^{i\theta}\in\dD\,|\,
\theta_1\leq\theta\leq \theta_2\big\}, \quad \theta_1 \in
[0,2\pi), \; \theta_1\leq \theta_2\leq \theta_1+2\pi \lb{5.1}
\end{equation}
and similarly for open arcs on $\dD$.

%%%%%%%%%%%%%%%%%%%%%%%%%%%%%%%%%%%%%%%%%%%%%%%%%%%%%%%%%
\begin{theorem}  \lb{t1.1}
Let $\alpha=\{\alpha_k\}_{k\in\bbZ}\subset\D$ be a
reflectionless sequence of Verblunsky coefficients. Let $U$
be the associated unitary CMV operator \eqref{1.1} $($cf.\
also \eqref{A.19}--\eqref{A.22}$)$ on $\ell^2(\bbZ)$ and
suppose that the spectrum of $U$ consists of a connected arc
of $\dD$,
\begin{equation}
\sigma(U)=\Arc\big(\big[e^{i\theta_0},e^{i\theta_1}\big]\big)
    \lb{5.12}
\end{equation}
with $\theta_0 \in [0,2\pi]$,
$\theta_0<\theta_1\leq\theta_0+2\pi$, and hence
$e^{i(\theta_0+\theta_1)/2}\in\Arc\big(\big(e^{i\theta_0},
e^{i\theta_1}\big)\big)$. Then
$\alpha=\{\alpha_k\}_{k\in\bbZ}$ is of the form,
\begin{equation}
\alpha_k=\alpha_0 g^k, \quad k\in\bbZ, \lb{5.13}
\end{equation}
where
\begin{equation}
g=-\exp(i(\theta_0+\theta_1)/2) \, \text{ and } \,
|\alpha_0|=\cos((\theta_1-\theta_0)/4). \lb{5.14}
\end{equation}
\end{theorem}
%%%%%%%%%%%%%%%%%%%%%%%%%%%%%%%%%%%%%%%%%%%%%%%%%%%%%%%%%

Here the sequence $\alpha=\{\alpha_k\}_{k\in\bbZ}\subset\D$ is called
{\it reflectionless} if
\begin{equation}
\text{for all $k\in\bbZ$, } \, M_+(\zeta,k)=-\ol{M_-(\zeta,k)} \,
\text{ for $\mu_0$-a.e.\ $\zeta\in \sigma_{\ess}(U)$,} \lb{1.14}
\end{equation}
where $M_\pm(\cdot,k)$, $k\in\Z$, denote the half-lattice Weyl--Titchmarsh
functions  of $U$ in \eqref{A.68} (cf.\ \cite{GZ05} for further details).
The case of reflectionless Verblunsky coefficients includes the
periodic case and certain quasi-periodic and almost periodic cases.

%%%%%%%%%%%%%%%%%%%%%%%%%%%%%%%%%%%%%%%%%%%%%%%%%%%%%%%%%
\section{Weyl--Titchmarsh Theory for CMV Operators on Half-Lattices}\lb{s2}
%%%%%%%%%%%%%%%%%%%%%%%%%%%%%%%%%%%%%%%%%%%%%%%%%%%%%%%%%

In this section we describe the Weyl--Titchmarsh theory for
CMV operators on half-lattices.

In the following, let $\ltz$ be the usual Hilbert space of
all square summable complex-valued sequences with scalar
product $(\cdot,\cdot)$ linear in the second argument. The
{\it standard basis} in $\ell^2(\bbZ)$ is denoted by
\begin{equation}
\{\delta_k\}_{k\in\bbZ}, \quad
\delta_k=(\dots,0,\dots,0,\underbrace{1}_{k},0,\dots,0,\dots)^\top,
\; k\in\bbZ.
\end{equation}
$\ell^\infty_0(\bbZ)$ denotes the set of sequences of
compact support (i.e.,
$f=\{f(k)\}_{k\in\bbZ}\in\ell^\infty_0(\bbZ)$ if there
exist $M(f), N(f) \in\bbZ$ such that $f(k)=0$ for $k<M(f)$
and $k>N(f)$). We use the analogous notation for compactly
supported sequences on half-lattices
$[k_0,\pm\infty)\cap\bbZ$, $k_0\in\bbZ$, and then write
$\ell^\infty_0([k_0,\pm\infty)\cap\bbZ)$, etc. For
$J\subseteq\bbR$ an interval, we will identify
$\ell^2(J\cap\bbZ)\oplus\ell^2(J\cap\bbZ)$ and
$\ell^2(J\cap\bbZ)\otimes\bbC^2$ and then use the
simplified notation $\ell^2(J\cap\bbZ)^2$. For simplicity,
the identity operator on $\ell^2(J\cap\bbZ)$ is abbreviated
by $I$ without separately indicating its dependence on $J$.

Moreover, we denote by $\D=\{ z\in\C \st \abs{z} < 1 \}$
the open unit disk in the complex plane $\C$, by $\dD=\{
\ze\in\C \st \abs{\ze} = 1 \}$ its counterclockwise oriented boundary, and we
freely use the notation employed in Appendix \ref{A}. By a {\it Laurent
polynomial}  we denote a finite linear combination of terms $z^k$,
$k\in\bbZ$, with complex-valued coefficients.

Throughout this paper we make the following basic assumption:

%%%%%%%%%%%%%%%%%%%%%%%%%%%%%%%%%%%%%%%%%%%%%%%%%%%%%%%%%
\begin{hypothesis} \lb{h2.1}
Let $\alpha$ be a sequence of complex numbers such that
\begin{equation} \label{A.15}
\alpha=\{\al_k\}_{k \in \Z} \subset \D.
\end{equation}
\end{hypothesis}
%%%%%%%%%%%%%%%%%%%%%%%%%%%%%%%%%%%%%%%%%%%%%%%%%%%%%%%%%

Given a sequence $\alpha$ satisfying \eqref{A.15}, we
define the sequence of positive real numbers
$\{\rho_k\}_{k\in\bbZ}$ and two sequences of complex
numbers with positive real parts $\{a_k\}_{k\in\bbZ}$ and
$\{b_k\}_{k\in\bbZ}$ by
\begin{align} \label{A.16}
\rho_k &= \sqrt{1-\abs{\al_k}^2}, \quad k\in\bbZ,
\\ \label{A.17}
a_k &= 1+\al_k, \quad k\in\bbZ,
\\ \label{A.18}
b_k &= 1-\al_k, \quad k \in \Z.
\end{align}
Following Simon \cite{Si04}, we call $\alpha_k$ the
Verblunsky coefficients in honor of Verblunsky's pioneering
work in the theory of orthogonal polynomials on the unit
circle \cite{Ve35}, \cite{Ve36}.

Next, we also introduce a sequence of $2\times 2$ unitary
matrices $\te_k$ by
\begin{equation} \label{A.19}
\te_k = \begin{pmatrix} -\al_k & \rho_k \\ \rho_k &
\ol{\al_k}
\end{pmatrix},
\quad k \in \Z,
\end{equation}
and two unitary operators $V$ and $W$ on $\ltz$ by their
matrix representations in the standard basis of
$\ell^2(\bbZ)$ as follows,
\begin{align} \label{A.20}
V &= \begin{pmatrix} \ddots & & &
\raisebox{-3mm}[0mm][0mm]{\hspace*{-5mm}\Huge $0$}  \\ &
\te_{2k-2} & & \\ & & \te_{2k} & & \\ &
\raisebox{0mm}[0mm][0mm]{\hspace*{-10mm}\Huge $0$} & &
\ddots
\end{pmatrix}, \quad
W &= \begin{pmatrix} \ddots & & &
\raisebox{-3mm}[0mm][0mm]{\hspace*{-5mm}\Huge $0$}
\\ & \te_{2k-1} &  &  \\ &  & \te_{2k+1} &  & \\ &
\raisebox{0mm}[0mm][0mm]{\hspace*{-10mm}\Huge $0$} & &
\ddots
\end{pmatrix},
\end{align}
where
\begin{align}
\begin{pmatrix}
V_{2k-1,2k-1} & V_{2k-1,2k} \\ V_{2k,2k-1}   & V_{2k,2k}
\end{pmatrix} =  \te_{2k},
\quad
\begin{pmatrix}
W_{2k,2k} & W_{2k,2k+1} \\ W_{2k+1,2k}  & W_{2k+1,2k+1}
\end{pmatrix} =  \te_{2k+1},
\quad k\in\Z.
\end{align}
Moreover, we introduce the unitary operator $U$ on
$\ltz$ by
\begin{equation} \label{A.22}
U = VW,
\end{equation}
or in matrix form, in the standard basis of $\ell^2(\bbZ)$,
by
\begin{align}
U = \begin{pmatrix} \ddots &&\hspace*{-8mm}\ddots
&\hspace*{-10mm}\ddots &\hspace*{-12mm}\ddots
&\hspace*{-14mm}\ddots &&&
\raisebox{-3mm}[0mm][0mm]{\hspace*{-6mm}{\Huge $0$}}
\\
&0& -\al_{0}\rho_{-1} & -\ol{\al_{-1}}\al_{0} &
-\al_{1}\rho_{0} & \rho_{0}\rho_{1}
\\
&& \rho_{-1}\rho_{0} &\ol{\al_{-1}}\rho_{0} &
-\ol{\al_{0}}\al_{1} & \ol{\al_{0}}\rho_{1} & 0
\\
&&&0& -\al_{2}\rho_{1} & -\ol{\al_{1}}\al_{2} &
-\al_{3}\rho_{2} & \rho_{2}\rho_{3}
\\
&&\raisebox{-4mm}[0mm][0mm]{\hspace*{-6mm}{\Huge $0$}} &&
\rho_{1}\rho_{2} & \ol{\al_{1}}\rho_{2} & -\ol{\al_{2}}\al_{3}
& \ol{\al_{2}}\rho_{3}&0
\\
&&&&&\hspace*{-14mm}\ddots &\hspace*{-14mm}\ddots
&\hspace*{-14mm}\ddots &\hspace*{-8mm}\ddots &\ddots
\end{pmatrix}. \lb{A.23}
\end{align}
Here terms of the form $-\ol{\alpha_k}\alpha_{k+1}$,
$k\in\Z$, represent the diagonal $k,k$-entries in the infinite
matrix \eqref{A.23}. We will call the operator $U$ on
$\ell^2(\bbZ)$ the CMV operator since
\eqref{A.19}--\eqref{A.23} in the context of the
semi-infinite (i.e., half-lattice) case were first obtained
by Cantero, Moral, and Vel\'azquez in \cite{CMV03}.

Finally, let $\UU$ denote the unitary operator on $\ltz^2$
defined by
\begin{equation} \label{A.24}
\UU = \begin{pmatrix} U & 0 \\ 0 & U^\top
\end{pmatrix}
=
\begin{pmatrix}
VW & 0 \\ 0 & WV
\end{pmatrix}
=
\begin{pmatrix}
0 & V \\ W & 0
\end{pmatrix}^2.
\end{equation}
One observes remnants of a certain ``supersymmetric'' structure in
$\left(\begin{smallmatrix} 0 & V \\ W & 0 \end{smallmatrix}\right)$
which is also reflected in the following result.

%%%%%%%%%%%%%%%%%%%%%%%%%%%%%%%%%%%%%%%%%%%%%%%%%%%%%%%%%%%%%%
\begin{lemma} \label{lA.1}
Let $z\in\bbC\backslash\{0\}$ and $\{u(z,k)\}_{k\in\bbZ},
\{v(z,k)\}_{k\in\bbZ}$ be sequences of complex functions.
Then the following items $(i)$--$(vi)$ are equivalent:
\begin{align}
& (i) \quad U u(z,\cdot) = z u(z,\cdot), \quad (W
u)(z,\cdot)=z v(z,\cdot). \label{A.25} \\ & (ii) \quad U^\top
v(z,\cdot) = z v(z,\cdot), \quad (V v)(z,\cdot)= u(z,\cdot).
\label{A.25a} \\ & (iii) \quad (W u)(z,\cdot) = z v(z,\cdot),
\quad (V v)(z,\cdot)= u(z,\cdot). \label{A.25b} \\ & (iv)
\quad \UU \binom{u(z,\cdot)}{v(z,\cdot)} = z
\binom{u(z,\cdot)}{v(z,\cdot)}, \quad (W u)(z,\cdot)=z
v(z,\cdot). \label{A.26} \\ & (v) \quad \UU
\binom{u(z,\cdot)}{v(z,\cdot)} = z
\binom{u(z,\cdot)}{v(z,\cdot)}. \quad (V
v)(z,\cdot)=u(z,\cdot). \label{A.27} \\ & (vi) \quad
\binom{u(z,k)}{v(z,k)} = T(z,k) \binom{u(z,k-1)}{v(z,k-1)},
\quad k\in\Z,  \label{A.28}
\end{align}
where the transfer matrices $T(z,k)$, $z\in\Cz$, $k\in\Z$,
are given by
\begin{equation}
T(z,k) = \begin{cases} \frac{1}{\rho_{k}} \begin{pmatrix}
\al_{k} & z \\ 1/z & \ol{\al_{k}} \end{pmatrix},  &
\text{$k$ odd,}  \\ \frac{1}{\rho_{k}} \begin{pmatrix}
\ol{\al_{k}} & 1 \\ 1 & \al_{k} \end{pmatrix}, & \text{$k$
even.} \end{cases} \label{A.31}
\end{equation}
\end{lemma}
%%%%%%%%%%%%%%%%%%%%%%%%%%%%%%%%%%%%%%%%%%%%%%%%%%%%%%%%%%%%%%
\begin{proof}
The equivalence of \eqref{A.25} and \eqref{A.25b} follows
from \eqref{A.22} after one defines $v(z,\cdot) =
\frac{1}{z} (W u)(z,\cdot)$. Since $\te_k^\top=\te_k$, one
has $V^\top=V$, $W^\top=W$ and hence, $U^\top = (VW)^\top =
WV$. Thus, defining $u(z,\cdot) = (V v)(z,\cdot)$, one gets
the equivalence of \eqref{A.25a} and \eqref{A.25b}. The
equivalence of \eqref{A.25b}, \eqref{A.26}, and
\eqref{A.27} follows immediately from \eqref{A.24}.

Next, we will prove that \eqref{A.25b} is equivalent to
\eqref{A.28}. Assuming $k$ to be odd one obtains the equivalence of
the following items $(i)$--$(v)$:
\begin{align}
& (i) \quad  \binom{u(z,k)}{v(z,k)} = T(z,k) \binom{u(z,k-1)}{v(z,k-1)}.
\\ & (ii) \quad
\rho_k \binom{u(z,k)}{v(z,k)} =\begin{pmatrix} \al_k & z \\
1/z & \ol{\al_k}
\end{pmatrix}
\binom{u(z,k-1)}{v(z,k-1)}. \\
& (iii) \quad
\begin{cases}
z v(z,k-1) = - \al_k u(z,k-1) + \rho_k u(z,k), \\ z\rho_k
v(z,k) = u(z,k-1) + \ol{\al_k} z v(z,k-1).
\end{cases} \\
& (iv) \quad
\begin{cases}
z v(z,k-1) = - \al_k u(z,k-1) + \rho_k u(z,k), \\ z v(z,k)
= \rho_k u(z,k-1) + \ol{\al_k} u(z,k).
\end{cases} \\
& (v) \quad
z \binom{v(z,k-1)}{v(z,k)} = \te_k
\binom{u(z,k-1)}{u(z,k)}.
\end{align}
If $k$ is even, one similarly proves that the following
items $(vi)$--$(viii)$ are equivalent:
\begin{align}
& (vi) \quad \binom{u(z,k)}{v(z,k)} = T(z,k)
\binom{u(z,k-1)}{v(z,k-1)}.
\\
& (vii) \quad
\rho_k \binom{v(z,k)}{u(z,k)} =
\begin{pmatrix}
\al_k & 1 \\ 1 & \ol{\al_k}
\end{pmatrix}
\binom{v(z,k-1)}{u(z,k-1)}. \\
& (viii) \quad
\binom{u(z,k-1)}{u(z,k)} = \te_k \binom{v(z,k-1)}{v(z,k)}.
\end{align}
Thus, taking into account \eqref{A.20}, one concludes that
\begin{align}
& \begin{cases} W u(z,\cdot) = z v(z,\cdot), \\
V v(z,\cdot) = u(z,\cdot) \end{cases}
\intertext{is equivalent to}
& \binom{u(z,k)}{v(z,k)} = T(z,k)
\binom{u(z,k-1)}{v(z,k-1)}, \quad k\in\Z.
\end{align}
\end{proof}
%%%%%%%%%%%%%%%%%%%%%%%%%%%%%%%%%%%%%%%%%%%%%%%%%%%%%%%%%%%%%%

We note that in studying solutions of $Uu(z,\cdot)=zu(z,\cdot)$ as in
Lemma \ref{lA.1}\,$(i)$, the purpose of the additional relation
$(Wu)(z,\cdot)=zv(z,\cdot)$ in \eqref{A.25} is to introduce a new variable
$v$ that improves our understanding of the structure of such solutions $u$.
An analogous comment applies to solutions of $U^\top
v(z,\cdot)=zv(z,\cdot)$ and the relation $(V v)(z,\cdot)=u(z,\cdot)$ in
Lemma \ref{lA.1}\,$(ii)$.

If one sets $\al_{k_0} = e^{is}$, $s\in [0,2\pi)$, for some
reference point $k_0\in\Z$, then the operator $U$ splits
into a direct sum of two half-lattice operators
$U_{-,k_0-1}^{(s)}$ and $U_{+,k_0}^{(s)}$ acting on
$\ell^2((-\infty,k_0-1]\cap\Z)$ and on
$\ell^2([k_0,\infty)\cap\Z)$, respectively. Explicitly, one
obtains
\begin{align}
\begin{split}
& U=U_{-,k_0-1}^{(s)} \oplus U_{+,k_0}^{(s)} \, \text{ in }
\, \ell^2((-\infty,k_0-1]\cap\Z) \oplus
\ell^2([k_0,\infty)\cap\Z) \\ & \text{if } \, \al_{k_0} =
e^{is}, \; s\in [0,2\pi).  \lb{2.30}
\end{split}
\end{align}
(Strictly speaking, setting $\al_{k_0} = e^{is}$, $s\in [0,2\pi)$, for some
reference point $k_0\in\Z$ contradicts our basic Hypothesis \ref{h2.1}.
However, as long as the exception to Hypothesis \ref{h2.1} refers to only
one or two sites (cf.\ also \eqref{2.177}), we will safely ignore this
inconsistency in favor of the notational simplicity it provides by avoiding
the introduction of a properly modified hypothesis on
$\{\alpha_k\}_{k\in\bbZ}$.) Similarly, one obtains $W_{-,k_0-1}^{(s)}$,
$V_{-,k_0-1}^{(s)}$ and $W_{+,k_0}^{(s)}$,
$V_{+,k_0}^{(s)}$ such that
\begin{equation}
U_{\pm,k_0}^{(s)} = V_{\pm,k_0}^{(s)} W_{\pm,k_0}^{(s)}.
\end{equation}
For simplicity we will abbreviate
\begin{equation}
U_{\pm,k_0} =
U_{\pm,k_0}^{(s=0)}=V_{\pm,k_0}^{(s=0)}W_{\pm,k_0}^{(s=0)}
=V_{\pm,k_0} W_{\pm,k_0}. \lb{2.31}
\end{equation}
In addition, we introduce on
$\ell^2([k_0,\pm\infty)\cap\Z)^2$ the half-lattice
operators $\UU_{\pm,k_0}^{(s)}$ by
\begin{equation}
\UU_{\pm,k_0}^{(s)} = \begin{pmatrix} U_{\pm,k_0}^{(s)} & 0
\\ 0 & (U_{\pm,k_0}^{(s)})^\top
\end{pmatrix}
=\begin{pmatrix} V_{\pm,k_0}^{(s)} W_{\pm,k_0}^{(s)} & 0 \\
0 & W_{\pm,k_0}^{(s)} V_{\pm,k_0}^{(s)}
\end{pmatrix}.
\end{equation}
By $\UU_{\pm,k_0}$ we denote the half-lattice operators
defined
for $s=0$,
\begin{equation}
\UU_{\pm,k_0} = \UU_{\pm,k_0}^{(s=0)} = \begin{pmatrix}
U_{\pm,k_0} & 0 \\ 0 & (U_{\pm,k_0})^\top
\end{pmatrix}
=\begin{pmatrix}
V_{\pm,k_0}W_{\pm,k_0} & 0 \\ 0 & W_{\pm,k_0}V_{\pm,k_0}
\end{pmatrix}.  \lb{2.33}
\end{equation}

%%%%%%%%%%%%%%%%%%%%%%%%%%%%%%%%%%%%%%%%%%%%%%%%%%%%%%%%%%%%
\begin{lemma} \label{lA.2}
Let $z\in\bbC\backslash\{0\}$, $k_0\in\bbZ$, and
 $\{\hatt p_+(z,k,k_0)\}_{k\geq k_0}$, $\{\hatt
r_+(z,k,k_0)\}_{k\geq k_0}$ be sequences of complex functions.
Then, the following items $(i)$--$(vi)$ are equivalent:
\begin{align}
& (i) \quad U_{+,k_0} \hatt p_+(z,\cdot,k_0) = z \hatt
p_+(z,\cdot,k_0), \quad W_{+,k_0}\hatt p_+(z,\cdot,k_0) = z
\hatt r_+(z,\cdot,k_0). \label{A.36} \\ & (ii) \quad
(U_{+,k_0})^\top \hatt r_+(z,\cdot,k_0) = z \hatt
r_+(z,\cdot,k_0), \quad V_{+,k_0}\hatt r_+(z,\cdot,k_0) =
\hatt p_+(z,\cdot,k_0). \label{A.37} \\ & (iii) \quad
W_{+,k_0} \hatt p_+(z,\cdot,k_0) = z \hatt r_+(z,\cdot,k_0),
\quad V_{+,k_0} \hatt r_+(z,\cdot,k_0) = \hatt
p_+(z,\cdot,k_0). \label{A.38} \\ & (iv) \quad \UU_{+,k_0}
\binom{\hatt p_+(z,\cdot,k_0)}{\hatt r_+(z,\cdot,k_0)} = z
\binom{\hatt p_+(z,\cdot,k_0)}{\hatt r_+(z,\cdot,k_0)}, \quad
W_{+,k_0}\hatt p_+(z,\cdot,k_0) = z \hatt r_+(z,\cdot,k_0).
\label{A.39} \\ & (v) \quad \UU_{+,k_0} \binom{\hatt
p_+(z,\cdot,k_0)}{\hatt r_+(z,\cdot,k_0)} = z \binom{\hatt
p_+(z,\cdot,k_0)}{\hatt r_+(z,\cdot,k_0)}, \quad
V_{+,k_0}\hatt r_+(z,\cdot,k_0) = \hatt p_+(z,\cdot,k_0).
\label{A.40} \\ & (vi) \quad \binom{\hatt p_+(z,k,k_0)}{\hatt
r_+(z,k,k_0)} = T(z,k) \binom{\hatt p_+(z,k-1,k_0)}{\hatt
r_+(z,k-1,k_0)}, \quad k
>      k_0, \label{A.41}  \\
& \ph{(vi)} \quad \;\, \text{ assuming } \,
\hatt p_+(z,k_0,k_0) = \begin{cases}
z\hatt r_+(z,k_0,k_0), & \text{$k_0$ odd}, \\ \hatt
r_+(z,k_0,k_0), & \text{$k_0$ even}. \end{cases} \lb{A.42}
\end{align}
Next, consider sequences $\{\hatt p_-(z,k,k_0)\}_{k\leq
k_0}$, $\{\hatt r_-(z,k,k_0)\}_{k\leq k_0}$. Then, the
following items $(vii)$--$(xii)$ are equivalent:
\begin{align}
& (vii) \quad U_{-,k_0} \hatt p_-(z,\cdot,k_0) = z \hatt
p_-(z,\cdot,k_0), \quad W_{-,k_0}\hatt p_-(z,\cdot,k_0) = z \hatt
r_-(z,\cdot,k_0). \\
& (viii) \quad (U_{-,k_0})^\top \hatt
r_-(z,\cdot,k_0) = z \hatt r_-(z,\cdot,k_0), \quad
V_{-,k_0}\hatt r_-(z,\cdot,k_0) = \hatt p_-(z,\cdot,k_0). \\
& (ix)
\quad W_{-,k_0} \hatt p_-(z,\cdot,k_0) = z \hatt
r_-(z,\cdot,k_0), \quad V_{-,k_0} \hatt r_-(z,\cdot,k_0) =
\hatt p_-(z,\cdot,k_0). \\
& (x) \quad \UU_{-,k_0}
\binom{\hatt p_-(z,\cdot,k_0)}{\hatt r_-(z,\cdot,k_0)} = z
\binom{\hatt p_-(z,\cdot,k_0)}{\hatt r_-(z,\cdot,k_0)},
\quad W_{-,k_0}\hatt p_-(z,\cdot,k_0) = z \hatt
r_-(z,\cdot,k_0). \\
& (xi) \quad \UU_{-,k_0} \binom{\hatt
p_-(z,\cdot,k_0)}{\hatt r_-(z,\cdot,k_0)} = z \binom{\hatt
p_-(z,\cdot,k_0)}{\hatt r_-(z,\cdot,k_0)}, \quad
V_{-,k_0}\hatt r_-(z,\cdot,k_0) = \hatt p_-(z,\cdot,k_0). \\
& (xii) \quad \binom{\hatt p_-(z,k-1),k_0}{\hatt
r_-(z,k-1,k_0)}  = T(z,k)^{-1} \binom{\hatt
p_-(z,k,k_0)}{\hatt r_-(z,k,k_0)}, \quad k \leq k_0,
\label{A.46}  \\
& \ph{(xii)} \quad \;\, \text{ assuming }
\hatt p_-(z,k_0,k_0) =\begin{cases} -\hatt r_-(z,k_0,k_0),
& \text{$k_0$ odd,}
\\ -z\hatt r_-(z,k_0,k_0), & \text{$k_0$ even.} \end{cases}
\lb{A.47}
\end{align}
\end{lemma}
%%%%%%%%%%%%%%%%%%%%%%%%%%%%%%%%%%%%%%%%%%%%%%%%%%%%%%%%%%%%
\begin{proof}
Repeating the first part of the proof of Lemma \ref{lA.1}
one obtains the equivalence of \eqref{A.36}, \eqref{A.37},
\eqref{A.38}, \eqref{A.39}, and \eqref{A.40}. Moreover,
repeating the second part of the proof of Lemma \ref{lA.1}
one obtains that
\begin{align}
(W_{+,k_0} \hatt p_+(z,\cdot,k_0))(k) &= z \hatt
r_+(z,k,k_0),
\\
(V_{+,k_0} \hatt r_+(z,\cdot,k_0))(k) &= \hatt
p_+(z,k,k_0), \quad k > k_0
\end{align}
is equivalent to
\begin{align}
\binom{\hatt p_+(z,k,k_0)}{\hatt r_+(z,k,k_0)} &=  T(z,k)
\binom{\hatt p_+(z,k-1,k_0)}{\hatt r_+(z,k-1,k_0)}, \quad k
>      k_0.
\end{align}
If $k_0$ is odd, then the operators $V_{+,k_0}$ and
$W_{+,k_0}$ have the following structure,
\begin{align}
V_{+,k_0} =
\begin{pmatrix}
\te_{k_0+1} & &
\raisebox{-1mm}[0mm][0mm]{\hspace*{-2mm}{\huge $0$}}
\\ & \te_{k_0+3} & \\
\raisebox{-1mm}[0mm][0mm]{\hspace*{0mm}{\huge $0$}} &&
\ddots
\end{pmatrix},
\quad
W_{+,k_0} =
\begin{pmatrix}
1 & & \raisebox{-1mm}[0mm][0mm]{\hspace*{-2mm}{\huge $0$}}
\\ & \te_{k_0+2} & \\ &
\raisebox{-1mm}[0mm][0mm]{\hspace*{-14mm}{\huge $0$}}
&\ddots
\end{pmatrix},
\end{align}
and hence,
\begin{equation}
W_{+,k_0} \hatt p_+(z,\cdot,k_0))(k_0) = z \hatt r_+(z,k_0,k_0)
\end{equation}
is equivalent to
\begin{equation}
\hatt p_+(z,k_0,k_0) = z\hatt r_+(z,k_0,k_0).
\end{equation}
Thus, one infers that \eqref{A.38} is equivalent to
\eqref{A.41}, \eqref{A.42} for $k_0$ odd. If $k_0$ is even,
then the operators $V_{+,k_0}$ and $W_{+,k_0}$ have the
following structure,
\begin{align}
V_{+,k_0} =
\begin{pmatrix}
1 & & \raisebox{-1mm}[0mm][0mm]{\hspace*{-2mm}{\huge $0$}}
\\ & \te_{k_0+2} & \\ &
\raisebox{-1mm}[0mm][0mm]{\hspace*{-14mm}{\huge $0$}}
&\ddots
\end{pmatrix},
\quad
W_{+,k_0} =
\begin{pmatrix}
\te_{k_0+1} & &
\raisebox{-1mm}[0mm][0mm]{\hspace*{-2mm}{\huge $0$}} \\ &
\te_{k_0+3} & \\
\raisebox{-1mm}[0mm][0mm]{\hspace*{0mm}{\huge $0$}} & &
\ddots
\end{pmatrix},
\end{align}
and hence,
\begin{equation}
(V_{+,k_0} \hatt r_+(z,\cdot,k_0))(k_0) = \hatt
p_+(z,k_0,k_0)
\end{equation}
is equivalent to
\begin{equation}
\hatt p_+(z,k_0,k_0) = \hatt r_+(z,k_0,k_0).
\end{equation}
Thus, one infers that \eqref{A.38} is equivalent to
\eqref{A.41}, \eqref{A.42} for $k_0$ even.

The results for $\hatt p_-(z,\cdot,k_0)$ and $\hatt
r_-(z,\cdot,k_0)$ are proved analogously.
\end{proof}
%%%%%%%%%%%%%%%%%%%%%%%%%%%%%%%%%%%%%%%%%%%%%%%%%%%%%%%%%%%%

Analogous comments to those made right after the proof of Lemma \ref{lA.1}
apply in the present context of Lemma \ref{lA.2}.

%%%%%%%%%%%%%%%%%%%%%%%%%%%%%%%%%%%%%%%%%%%%%%%%%%%%%%%%%%%%
\begin{definition} \label{dA.2}
We denote by $\Big(\begin{smallmatrix}p_+(z,k,k_0) \\
r_+(z,k,k_0)\end{smallmatrix}\Big)_{k\geq k_0}$ and
$\Big(\begin{smallmatrix} q_+(z,k,k_0) \\ s_+(z,k,k_0)
\end{smallmatrix}\Big)_{k\geq k_0}$,
$z\in\bbC\backslash\{0\}$, two linearly independent
solutions
of \eqref{A.41} with the following initial conditions:
\begin{align}
\binom{p_+(z,k_0,k_0)}{r_+(z,k_0,k_0)} = \begin{cases} \binom{z}{1}, &
\text{$k_0$ odd,} \\[1mm]
\binom{1}{1}, & \text{$k_0$ even,} \end{cases} \quad
\binom{q_+(z,k_0,k_0)}{s_+(z,k_0,k_0)} = \begin{cases} \binom{z}{-1}, &
\text{$k_0$ odd,} \\[1mm]
\binom{-1}{1}, & \text{$k_0$ even.} \end{cases} \label{A.48}
\end{align}
Similarly, we denote by
$\Big(\begin{smallmatrix}p_-(z,k,k_0)\\r_-(z,k,k_0)
\end{smallmatrix}\Big)_{k\leq k_0}$ and
$\Big(\begin{smallmatrix}q_-(z,k,k_0)\\
s_-(z,k,k_0)\end{smallmatrix}\Big)_{k\leq k_0}$,
$z\in\bbC\backslash\{0\}$, two linearly independent
solutions of \eqref{A.46} with the following initial
conditions:
\begin{align}
\binom{p_-(z,k_0,k_0)}{r_-(z,k_0,k_0)} = \begin{cases} \binom{1}{-1}, &
\text{$k_0$ odd,} \\[1mm]
\binom{-z}{1}, & \text{$k_0$ even,} \end{cases} \quad
\binom{q_-(z,k_0,k_0)}{s_-(z,k_0,k_0)} = \begin{cases} \binom{1}{1}, &
\text{$k_0$ odd, }\\[1mm]
\binom{z}{1}, & \text{$k_0$ even.} \end{cases} \label{A.50}
\end{align}
Using \eqref{A.28} one extends
$\Big(\begin{smallmatrix}p_+(z,k,k_0)\\r_+(z,k,k_0)
\end{smallmatrix}\Big)_{k\geq k_0}$,
$\Big(\begin{smallmatrix}q_+(z,k,k_0)\\s_+(z,k,k_0)
\end{smallmatrix}\Big)_{k\geq k_0}$, $z\in\bbC\backslash\{0\}$,
to $k < k_0$. In the same manner, one extends
$\Big(\begin{smallmatrix}p_-(z,k,k_0)\\r_-(z,k,k_0)
\end{smallmatrix}\Big)_{k\leq k_0}$ and
$\Big(\begin{smallmatrix}q_-(z,k,k_0)\\s_-(z,k,k_0)
\end{smallmatrix}\Big)_{k\leq k_0}$, $z\in\bbC\backslash\{0\}$,
to $k> k_0$. These extensions will be denoted by
$\Big(\begin{smallmatrix}p_\pm(z,k,k_0)\\r_\pm
(z,k,k_0)\end{smallmatrix}\Big)_{k\in\Z}$ and
$\Big(\begin{smallmatrix}q_\pm (z,k,k_0)\\s_\pm
(z,k,k_0)\end{smallmatrix}\Big)_{k\in\Z}$.\ Moreover, it
follows from \eqref{A.28} that $p_\pm(z,k,k_0)$,
$q_\pm(z,k,k_0)$, $r_\pm(z,k,k_0)$, and $s_\pm(z,k,k_0)$,
$k,k_0\in\Z$, are Laurent polynomials in $z$.
\end{definition}
%%%%%%%%%%%%%%%%%%%%%%%%%%%%%%%%%%%%%%%%%%%%%%%%%%%%%%%%%%%%%%

\smallskip
In particular, one computes
\smallskip
%%%%%%%%%%%%%%%%%%%%%%%%%%%%%%%%%%%%%%%%%%%%%%%%%%%%%%%%%%%%%%
\begin{align*}
\begin{array}{|c|c|c|c|}
\hline k & k_0-1 & k_0 \text{ odd} & k_0+1
\\ \hline \ph{\Bigg|}
\displaystyle \binom{p_+(z,k,k_0)}{r_+(z,k,k_0)} &
\displaystyle \frac{1}{\rho_{k_0}}
\binom{z(1-\ol{\al_{k_0}})}{1-\al_{k_0}} & \displaystyle
\binom{z}{1} & \displaystyle \frac{1}{\rho_{k_0+1}}
\binom{1+\ol{\al_{k_0+1}}z}{z+\al_{k_0+1}}
\\ \hline \ph{\Bigg|}
\displaystyle \binom{q_+(z,k,k_0)}{s_+(z,k,k_0)} &
\displaystyle
\frac{1}{\rho_{k_0}}\binom{z(-1-\ol{\al_{k_0}})}{1+\al_{k_0}}
& \displaystyle \binom{z}{-1} & \displaystyle
\frac{1}{\rho_{k_0+1}}
\binom{-1+\ol{\al_{k_0+1}}z}{z-\al_{k_0+1}}
\\ \hline \ph{\Bigg|}
\displaystyle \binom{p_-(z,k,k_0)}{r_-(z,k,k_0)} &
\displaystyle
\frac{1}{\rho_{k_0}}\binom{-z-\ol{\al_{k_0}}}{1/z+\al_{k_0}}
& \displaystyle \binom{1}{-1} & \displaystyle
\frac{1}{\rho_{k_0+1}}
\binom{-1+\ol{\al_{k_0+1}}}{1-\al_{k_0+1}}
\\ \hline \ph{\Bigg|}
\displaystyle \binom{q_-(z,k,k_0)}{s_-(z,k,k_0)} &
\displaystyle
\frac{1}{\rho_{k_0}}\binom{z-\ol{\al_{k_0}}}{1/z-\al_{k_0}}
& \displaystyle \binom{1}{1} & \displaystyle
\frac{1}{\rho_{k_0+1}}
\binom{1+\ol{\al_{k_0+1}}}{1+\al_{k_0+1}}
\\ \hline \hline
k & k_0-1 & k_0 \text{ even} & k_0+1
\\ \hline \ph{\Bigg|}
\displaystyle \binom{p_+(z,k,k_0)}{r_+(z,k,k_0)} &
\displaystyle
\frac{1}{\rho_{k_0}}\binom{1-\al_{k_0}}{1-\ol{\al_{k_0}}} &
\displaystyle \binom{1}{1} & \displaystyle
\frac{1}{\rho_{k_0+1}}
\binom{z+\al_{k_0+1}}{1/z+\ol{\al_{k_0+1}}}
\\ \hline \ph{\Bigg|}
\displaystyle \binom{q_+(z,k,k_0)}{s_+(z,k,k_0)} &
\displaystyle
\frac{1}{\rho_{k_0}}\binom{1+\al_{k_0}}{-1-\ol{\al_{k_0}}}
& \displaystyle \binom{-1}{1} & \displaystyle
\frac{1}{\rho_{k_0+1}}
\binom{z-\al_{k_0+1}}{-1/z+\ol{\al_{k_0+1}}}
\\ \hline \ph{\Bigg|}
\displaystyle \binom{p_-(z,k,k_0)}{r_-(z,k,k_0)} &
\displaystyle
\frac{1}{\rho_{k_0}}\binom{1+\al_{k_0}z}{-z-\ol{\al_{k_0}}}
& \displaystyle \binom{-z}{1} & \displaystyle
\frac{1}{\rho_{k_0+1}}
\binom{z(1-\al_{k_0+1})}{-1+\ol{\al_{k_0+1}}}
\\ \hline \ph{\Bigg|}
\displaystyle \binom{q_-(z,k,k_0)}{s_-(z,k,k_0)} &
\displaystyle
\frac{1}{\rho_{k_0}}\binom{1-\al_{k_0}z}{z-\ol{\al_{k_0}}}
& \displaystyle \binom{z}{1} & \displaystyle
\frac{1}{\rho_{k_0+1}}
\binom{z(1+\al_{k_0+1})}{1+\ol{\al_{k_0+1}}}
\\ \hline
\end{array}
\end{align*}
%%%%%%%%%%%%%%%%%%%%%%%%%%%%%%%%%%%%%%%%%%%%%%%%%%%%%%%%%%%%%%

\medskip

%%%%%%%%%%%%%%%%%%%%%%%%%%%%%%%%%%%%%%%%%%%%%%%%%%%%%%%%%%%%%%
\begin{remark} \lb{rA.2a}
We note that Lemmas \ref{lA.1} and \ref{lA.2} are crucial for many of the
proofs to follow. For instance, we note that the equivalence of items $(i)$
and $(vi)$ in Lemma \ref{lA.1} proves that for each
$z\in\bbC\backslash\{0\}$, the solutions $\{u(z,k)\}_{k\in\bbZ}$ of $U
u(z,\cdot)=zu(z,\cdot)$ form a two-dimensional space, which implies that such
solutions are linear combinations of $\{p_{\pm}(z,k,k_0)\}_{k\in\bbZ}$ and
$\{q_{\pm}(z,k,k_0)\}_{k\in\bbZ}$ (with $z$-dependent coefficients). This
equivalence also proves that any solution of $U u(z,\cdot)=zu(z,\cdot)$ is
determined by its values at a site $k_0$ of $u$ and the auxiliary variable
$v$. Moreover, taking into account item $(vi)$ of Lemma \ref{lA.1}, this
also implies that such a solution is determined by its values at two
consecutive sites $k_0-1$ and $k_0$. Similar comments apply to the
solutions of $U^\top v(z,\cdot)=zv(z,\cdot)$. In the context of Lemma
\ref{lA.2}, we remark that its importance lies in the fact that it shows
that in the case of half-lattice CMV operators, the analogous equations have
a one-dimensional space of solutions for each $z\in\bbC\backslash\{0\}$, due
to the restriction on $k_0$ that appears in items $(vi)$ and $(xii)$ of Lemma
\ref{lA.2}. As a consequence, the corresponding solutions are determined by
their value at a single site $k_0$.
\end{remark}
%%%%%%%%%%%%%%%%%%%%%%%%%%%%%%%%%%%%%%%%%%%%%%%%%%%%%%%%%%%%%%

Next, we introduce the following modified Laurent
polynomials $\wti p_\pm(z,k,k_0)$ and $\wti
q_\pm(z,k,k_0)$, $z\in\Cz$, $k,k_0\in\Z$, as follows,
\begin{align}
\wti p_+(z,k,k_0) &= \begin{cases} p_+(z,k,k_0)/z, &
\text{$k_0$ odd,}
\\ p_+(z,k,k_0), &  \text{$k_0$ even,} \end{cases} \lb{2.99} \\
\wti q_+(z,k,k_0) &= \begin{cases} q_+(z,k,k_0)/z, &
\text{$k_0$ odd,} \\ q_+(z,k,k_0), & \text{$k_0$ even,}
\end{cases} \lb{2.100} \\ \wti p_-(z,k,k_0) &=
\begin{cases} p_-(z,k,k_0), & \text{$k_0$ odd,} \\
p_-(z,k,k_0)/z, & \text{$k_0$ even,} \end{cases} \lb{2.101}
\\ \wti q_-(z,k,k_0) &= \begin{cases} q_-(z,k,k_0), &
\text{$k_0$ odd,} \\ q_-(z,k,k_0)/z, & \text{$k_0$ even.}
\end{cases} \lb{2.102}
\end{align}

%%%%%%%%%%%%%%%%%%%%%%%%%%%%%%%%%%%%%%%%%%%%%%%%%%%%%%%%%%%%%%%
\begin{remark} \label{rA.3}
By Lemma \ref{lA.2},
$\Big(\begin{smallmatrix}p_\pm(z,k,k_0)\\
r_\pm(z,k,k_0)\end{smallmatrix}\Big)_{k\gtreqless k_0}$,
$z\in\bbC\backslash\{0\}$, $k_0\in\Z$, are generalized
eigenvectors of the operators $\UU_{\pm,k_0}$. Moreover,
by Lemma \ref{lA.1}, $\Big(\begin{smallmatrix}p_\pm
(z,k,k_0)\\r_\pm (z,k,k_0)\end{smallmatrix}\Big)_{k\in\Z}$
and $\Big(\begin{smallmatrix}q_\pm (z,k,k_0)\\s_\pm
(z,k,k_0)\end{smallmatrix}\Big)_{k\in\Z}$,
$z\in\bbC\backslash\{0\}$, $k_0\in\Z$, are generalized
eigenvectors of $\UU$.
\end{remark}
%%%%%%%%%%%%%%%%%%%%%%%%%%%%%%%%%%%%%%%%%%%%%%%%%%%%%%%%%%%%%%%

%%%%%%%%%%%%%%%%%%%%%%%%%%%%%%%%%%%%%%%%%%%%%%%%%%%%%%%%%%%%%%%
\begin{lemma} \label{lA.3a}
The Laurent polynomials $\wti p_\pm(z,k,k_0)$,
$r_\pm(z,k,k_0)$, $\wti q_\pm(z,k,k_0)$, and
$s_\pm(z,k,k_0)$ satisfy the following relations for all
$z\in\bbC\backslash\{0\}$ and $k,k_0\in\bbZ$,
\begin{align}
r_+(z,k,k_0) &= \ol{\wti p_+(1/\ol{z},k,k_0)}, \lb{2.59}
\\
s_+(z,k,k_0) &= -\ol{\wti q_+(1/\ol{z},k,k_0)}, \lb{2.60}
\\
r_-(z,k,k_0) &= -\ol{\wti p_-(1/\ol{z},k,k_0)},
\\
s_-(z,k,k_0) &= \ol{\wti q_-(1/\ol{z},k,k_0)}.
\end{align}
\end{lemma}
%%%%%%%%%%%%%%%%%%%%%%%%%%%%%%%%%%%%%%%%%%%%%%%%%%%%%%%%%%%%%%%
\begin{proof}
Let $\{u(z,k)\}_{k\in\Z}$, $\{v(z,k)\}_{k\in\Z}$ be two
sequences of complex functions, then the following items
$(i)$--$(iii)$ are seen to be equivalent:
\begin{align}
& (i) \quad W u(z,\cdot) = z v(z,\cdot), \quad V v(z,\cdot)
= u(z,\cdot).  \lb{2.70a} \\
& (ii) \quad \frac{1}{z} u(z,\cdot) = W^{*}v(z,\cdot),
\quad v(z,\cdot) = V^{*}u(z,\cdot).  \lb{2.71a} \\
& (iii) \quad \frac{1}{\ol{z}} \ol{u(z,\cdot)} = W
\ol{v(z,\cdot)}, \quad \ol{v(z,\cdot)} = V\ol{u(z,\cdot)}, \lb{2.72a}
\end{align}
where equations \eqref{2.70a}--\eqref{2.72a} are meant in the
algebraic sense and hence $V$, $V^*$, $W$, and $W^*$ are
considered as difference expressions rather than difference
operators. Thus, the assertion of the Lemma follows from
Lemma \ref{lA.2}, Definition \ref{dA.2}, and equalities
\eqref{2.99}--\eqref{2.102}.
\end{proof}
%%%%%%%%%%%%%%%%%%%%%%%%%%%%%%%%%%%%%%%%%%%%%%%%%%%%%%%%%%%%%%%

%%%%%%%%%%%%%%%%%%%%%%%%%%%%%%%%%%%%%%%%%%%%%%%%%%%%%%%%%%%%%%%
\begin{lemma} \label{lA.4}
Let $k_0\in\bbZ$. Then the sets of Laurent polynomials
$\{p_+(\cdot,k,k_0)\}_{k\geq k_0}$ $($resp.,
$\{p_-(\cdot,k,k_0)\}_{k\leq k_0}$$)$ and
$\{r_+(\cdot,k,k_0)\}_{k\geq k_0}$ $($resp.,
$\{r_-(\cdot,k,k_0)\}_{k\leq k_0}$$)$ form orthonormal
bases in $\Lt{+}$ $($resp., $\Lt{-}$$)$, where
\begin{equation}
d\mu_\pm(\zeta,k_0) =
d(\de_{k_0},E_{U_{\pm,k_0}}(\zeta)
\de_{k_0})_{\ell^2([k_0,\pm\infty)\cap\Z)},
\quad \zeta\in\dD, \label{A.51}
\end{equation}
and $dE_{U_{\pm,k_0}}(\cdot)$ denote the operator-valued
spectral measures of the operators $U_{\pm,k_0}$,
\begin{equation}
U_{\pm,k_0}=\oint_{\dD} dE_{U_{\pm,k_0}}(\zeta)\,\zeta.
\end{equation}
\end{lemma}
%%%%%%%%%%%%%%%%%%%%%%%%%%%%%%%%%%%%%%%%%%%%%%%%%%%%%%%%%%%%%%%
\begin{proof}
It follows from the definition of the transfer matrix
$T(z,k)$ in \eqref{A.31} and the recursion relations
\eqref{A.41} and \eqref{A.46} that
\begin{align}
\begin{split}
&\ol{\spn{\{p_\pm(\cdot,k,k_0)\}_{k\gtreqless k_0}}} =
\ol{\spn{\{r_\pm(\cdot,k,k_0)\}_{k\gtreqless k_0}}}  \\ & \quad
= \ol{\spn{\{\zeta^k\}_{k\in\Z}}} = L^2(\dD;d\mu), \lb{2.72}
\end{split}
\end{align}
where $d\mu$ is any finite (nonnegative) Borel measure on
$\dD$. Thus, one concludes that the systems of Laurent
polynomials $\{p_\pm(\cdot,k,k_0)\}_{k\gtreqless k_0}$ and
$\{r_\pm(\cdot,k,k_0)\}_{k\gtreqless k_0}$ are complete in
$\Lt{\pm}$.

Next, consider the following equations
\begin{align}
(U_{+,k_0})^\top \de_k &= \sum_{j=k-2}^{k+2} (
U_{+,k_0})^\top(j,k)\de_j = \sum_{j=k-2}^{k+2} (
U_{+,k_0})(k,j)\de_j, \lb{2.73}
\\
(U_{+,k_0}) \de_k &= \sum_{j=k-2}^{k+2} ( U_{+,k_0})(j,k)\de_j =
\sum_{j=k-2}^{k+2} ( U_{+,k_0})^\top(k,j)\de_j, \lb{2.74}
\end{align}
and
\begin{align}
z\hatt p_+(z,k,k_0) &= (U_{+,k_0} \hatt p_+(z,\cdot,k_0))(k) =
\sum_{j=k-2}^{k+2} (U_{+,k_0})(k,j) \hatt p_+(z,j,k_0),
\lb{2.75}
\\
z \hatt r_+(z,k,k_0) &= ((U_{+,k_0})^\top \hatt
r_+(z,\cdot,k_0))(k) = \sum_{j=k-2}^{k+2} (U_{+,k_0})^\top(k,j)
\hatt r_+(z,j,k_0). \lb{2.76}
\end{align}
By Lemma \ref{lA.2} the latter ones have unique solutions $\wti
p_+(z,k,k_0)$ and $r_+(z,k,k_0)$ satisfying $\wti p_+(z,k_0,k_0)
= r_+(z,k_0,k_0) = 1$. Moreover, due to the algebraic nature of
the proof of Lemma \ref{lA.2}, \eqref{2.75} and \eqref{2.76}
remain valid if $z\in\bbC\backslash\{0\}$ is replaced by a
unitary operator on a Hilbert space and the left- and right-hand
sides are applied to the vector $\de_{k_0}$. Thus, $\{\wti
p_+((U_{+,k_0})^\top,k,k_0)\delta_{k_0}\}_{k\geq k_0}$ and
$\{r_+(U_{+,k_0},k,k_0)\delta_{k_0}\}_{k\geq k_0}$ are the
unique solutions of
\begin{align}
(U_{+,k_0})^\top \hatt p_+((U_{+,k_0})^\top,k,k_0)\delta_{k_0}
&= \sum_{j=k-2}^{k+2} (U_{+,k_0})(k,j) \hatt
p_+((U_{+,k_0})^\top,j,k_0) \delta_{k_0}, \lb{2.77}
\\
U_{+,k_0} \hatt r_+(U_{+,k_0},k,k_0)\delta_{k_0} &=
\sum_{j=k-2}^{k+2} (U_{+,k_0})^\top(k,j) \hatt
r_+(U_{+,k_0},j,k_0) \delta_{k_0} \lb{2.78}
\end{align}
%\begin{equation}
%(U_{+,k_0})^\top p((U_{+,k_0})^\top,\cdot,k_0)
%=U_{+,k_0} p((U_{+,k_0})^\top,\cdot,k_0)
%\end{equation}
%and
%\begin{equation}
%U_{+,k_0} r(U_{+,k_0},\cdot,k_0)=(U_{+,k_0})^\top r(U_{+,k_0},\cdot,k_0)
%\end{equation}
with value $\delta_{k_0}$ at $k=k_0$, respectively. In
particular, comparing \eqref{2.73}, \eqref{2.74} with
\eqref{2.77}, \eqref{2.78}, one concludes that for $k\geq k_0$,
\begin{align}
\de_k &= \wti{p}_+((U_{+,k_0})^\top,k,k_0) \de_{k_0},
\\
\de_k &= r_+(U_{+,k_0},k,k_0) \de_{k_0}.
\end{align}
Using the spectral representation for the operators
$U_{+,k_0}$ and $( U_{+,k_0})^\top$ one obtains (all scalar
products $(\cdot,\cdot)$ in the remainder of this proof are
with respect to the Hilbert space
$\ell^2([k_0,\pm\infty)\cap\Z)$ and for simplicity we omit
the corresponding subscript in $(\cdot,\cdot)$),
\begin{align}
(\de_k,\de_\ell) &= \oint_\dD
d(\de_{k_0},E_{(U_{+,k_0})^\top}(\ze)\de_{k_0})\,
\ol{p_+(\ze,k,k_0)} p_+(\ze,\ell,k_0), \\ (\de_k,\de_\ell)
&= \oint_\dD d(\de_{k_0},E_{ U_{+,k_0}}(\ze)
\de_{k_0})\,\ol{r_+(\ze,k,k_0)} r_+(\ze,\ell,k_0), \quad
k,\ell\in\Z.
\end{align}
Finally,  one notes that
\begin{align}
d\mu_{+}(\ze,k_0) =
d(\de_{k_0},E_{U_{+,k_0}}(\ze)\de_{k_0}) =
d(\de_{k_0},E_{(U_{+,k_0})^\top}(\ze)\de_{k_0})
\end{align}
since
\begin{align}
\begin{split}
& \oint_\dD d\mu_{+}(\ze,k_0)\, \ze^k = \Big(\de_{k_0},
U_{+,k_0}^k\de_{k_0}\Big) = \left(\de_{k_0},\left(
U_{+,k_0}^k\right)^\top\de_{k_0}\right) \\
& \quad = \Big(\de_{k_0},\left(
U_{+,k_0}^\top\right)^k\de_{k_0}\Big) = \oint_\dD
d(\de_{k_0},E_{(U_{+,k_0})^\top}(\ze)\de_{k_0}) \, \ze^k,
\quad k\in\Z.
\end{split}
\end{align}
Thus, the Laurent polynomials $\{p_+(\cdot,k,k_0)\}_{k\geq
k_0}$ and $\{r_+(\cdot,k,k_0)\}_{k\geq k_0}$ are orthonormal in
$\Lt{+}$.

The results for
$\{p_-(\cdot,k,k_0)\}_{k\leq
k_0}$ and $\{r_-(\cdot,k,k_0)\}_{k\leq k_0}$ are proved similarly.
\end{proof}
%%%%%%%%%%%%%%%%%%%%%%%%%%%%%%%%%%%%%%%%%%%%%%%%%%%%%%%%%%%%%%%

We note that the measures $d\mu_\pm(\cdot,k_0)$, $k_0\in\bbZ$, are
not only nonnegative but also supported on an infinite set.

%%%%%%%%%%%%%%%%%%%%%%%%%%%%%%%%%%%%%%%%%%%%%%%%%%%%%%%%%%%%%%%
\begin{remark} \label{r2.8}
In connection with our introductory remarks in
\eqref{1.3}--\eqref{1.4e} we note that
$d\sigma_+=d\mu_+(\cdot,0)$ and
\begin{align}
\begin{split}
p_+(\zeta,k,0)&=\begin{cases} \gamma_k \zeta^{-(k-1)/2}
\varphi_+(\zeta,k), & \text{$k$ odd,} \\[1mm]
\gamma_k \zeta^{-k/2} \varphi^*_+(\zeta,k), & \text{$k$ even,}
\end{cases}  \\
r_+(\zeta,k,0)&=\begin{cases} \gamma_{k} \zeta^{-(k+1)/2}
\varphi^*_+(\zeta,k), & \text{$k$ odd,} \\[1mm]
\gamma_k \zeta^{-k/2} \varphi_+(\zeta,k), & \text{$k$ even;}
\end{cases} \quad \zeta\in\dD.
\end{split}
\end{align}
\end{remark}
%%%%%%%%%%%%%%%%%%%%%%%%%%%%%%%%%%%%%%%%%%%%%%%%%%%%%%%%%%%%%%%

Let $\phi \in C(\dD)$ and define the operator of multiplication by
$\phi$, $M_{\pm,k_0}(\phi)$, in $L^2(\dD;d\mu_\pm(\cdot,k_0))$ by
\begin{equation}
(M_{\pm,k_0}(\phi)f)(\zeta)=\phi(\zeta)f(\zeta), \quad
f\in  L^2(\dD;d\mu_\pm(\cdot,k_0)).   \lb{2.86}
\end{equation}
In the special case $\phi= id$ (where $id(\zeta)=\zeta$, $\zeta\in\dD$),
the corresponding multiplication operator is denoted by $M_{\pm,k_0}(id)$.
The spectrum of
$M_{\pm,k_0}(\phi)$ is given by
\begin{equation}
\sigma(M_{\pm,k_0}(\phi))=\essran_{d\mu_{\pm}(\cdot,k_0)}(\phi),
\end{equation}
where the essential range of $\phi$ with respect to a
measure $d\mu$ on $\dD$ is defined by
\begin{equation}
\essran_{d\mu}(\phi)=\{z\in\bbC\,|\, \text{for all
$\varepsilon>0$,} \, \mu(\{\zeta\in\dD \,|\,
|\phi(\zeta)-z|<\varepsilon\})>0\}.   \lb{2.87}
\end{equation}

%%%%%%%%%%%%%%%%%%%%%%%%%%%%%%%%%%%%%%%%%%%%%%%%%%%%%%%%%%%%%%%
\begin{corollary} \label{cA.4}
Let $k_0\in\bbZ$ and $\phi\in C(\dD)$. Then the operators
$\phi(U_{\pm,k_0})$ and $\phi(U_{\pm,k_0}^\top)$ are
unitarily equivalent to the operators $M_{\pm,k_0}(\phi)$
of multiplication by $\phi$ on $\Lt{\pm}$. In particular,
\begin{align}
& \si(\phi(U_{\pm,k_0})) =\si(\phi(U_{\pm,k_0}^\top))
= \essran_{d\mu_{\pm}(\cdot,k_0)}(\phi),  \\
& \si(U_{\pm,k_0}) =\si(U_{\pm,k_0}^\top) =\supp \, (d\mu_\pm(\cdot,k_0))
\end{align}
and the spectrum of $U_{\pm,k_0}$ is simple.
\end{corollary}
%%%%%%%%%%%%%%%%%%%%%%%%%%%%%%%%%%%%%%%%%%%%%%%%%%%%%%%%%%%%%%%
\begin{proof}
Consider the following linear maps $\dot \cU_\pm$ from
$\ell^\infty_0([k_0,\pm\infty)\cap\bbZ)$ into the set of
Laurent polynomials on $\dD$ defined by
\begin{equation}
(\dot \cU_\pm f)(\ze) = \sum_{k=k_0}^{\pm\infty}
r_\pm(\ze,k,k_0) f(k), \quad f\in
\ell^\infty_0([k_0,\pm\infty)\cap\bbZ).
\end{equation}
A simple calculation for $F(\ze) = (\dot \cU_\pm f)(\ze)$,
$f\in \ell^\infty_0([k_0,\pm\infty)\cap\bbZ)$, shows that
\begin{align}
\sum_{k=k_0}^{\pm\infty} \abs{f(k)}^2 = \oint_\dD
d\mu_\pm(\ze,k_0) \, \abs{F(\ze)}^2.
\end{align}
Since $\ell^\infty_0([k_0,\pm\infty)\cap\bbZ)$ is dense in
$\ell^2([k_0,\pm\infty)\cap\bbZ)$, $\dot \cU_\pm$ extend to
bounded linear operators $\cU_\pm\colon
\ell^2([k_0,\pm\infty)\cap\bbZ) \to \Lt{\pm}$. Since by
\eqref{2.72}, the sets of Laurent  polynomials are dense in
$\Lt{\pm}$, the maps $\cU_\pm$ are onto and one infers
\begin{equation}
(\cU_\pm^{-1}F)(k) = \oint_\dD
d\mu_{\pm}(\ze,k_0)\,\ol{r_\pm(\ze,k,k_0)} F(\ze), \quad
F\in \Lt{\pm}.
\end{equation}
In particular, $\cU_\pm$ are unitary. Moreover, we claim
that $\cU_\pm$ map the operators $\phi(U_{\pm,k_0})$ on
$\ell^2([k_0,\pm\infty)\cap\bbZ)$ to the operators
$M_{\pm,k_0}(\phi)$ of multiplication by $\phi$ on
$\Lt{\pm}$,
\begin{align}
\cU_\pm \phi(U_{\pm,k_0}) \cU_\pm^{-1} = M_{\pm,k_0}(\phi).
\end{align}
Indeed,
\begin{align}
& (\cU_\pm \phi(U_{\pm,k_0}) \cU_\pm^{-1} F(\cdot))(\ze) =
(\cU_\pm \phi(U_{\pm,k_0}) f(\cdot))(\ze) \no \\ & \quad =
\sum_{k=k_0}^{\pm\infty}
(\phi(U_{\pm,k_0})f(\cdot))(k)r_\pm(\ze,k,k_0) =
\sum_{k=k_0}^{\pm\infty} (\phi(U_{\pm,k_0}^\top)
r_\pm(\ze,\cdot,k_0))(k)f(k) \no \\ &\quad =
\sum_{k=k_0}^{\pm\infty} \phi(\ze) r_\pm(\ze,k,k_0) f(k) =
\phi(\ze) F(\ze)  \no \\ & \quad
=(M_{\pm,k_0}(\phi)F)(\ze), \quad F \in \Lt{\pm}.
\end{align}

The result for $\phi(U_{\pm,k_0}^\top)$ is proved analogously.
\end{proof}
%%%%%%%%%%%%%%%%%%%%%%%%%%%%%%%%%%%%%%%%%%%%%%%%%%%%%%%%%%%%%%%

%%%%%%%%%%%%%%%%%%%%%%%%%%%%%%%%%%%%%%%%%%%%%%%%%%%%%%%%%%%%%%%
\begin{corollary} \label{cA.5}
Let $k_0\in\bbZ$. \\ The
Laurent polynomials $\{p_+(\cdot,k,k_0)\}_{k\geq k_0}$ can
be constructed by Gram--Schmidt orthogonalizing
\begin{equation}
\begin{cases} \ze,\,1,\,\ze^2,\,\ze^{-1},\,\ze^3,\,\ze^{-2},\dots, &
\text{$k_0$ odd,} \\
1,\,\ze,\,\ze^{-1},\,\ze^2,\,\ze^{-2},\ze^3,\,\dots, &
\text{$k_0$ even}
\end{cases}
\end{equation}
in $\Lt{+}$. \\ The Laurent polynomials
$\{r_+(\cdot,k,k_0)\}_{k\geq k_0}$ can be constructed by
Gram--Schmidt orthogonalizing
\begin{equation}
\begin{cases} 1,\,\ze,\,\ze^{-1},\,\ze^2,\,\ze^{-2},\ze^3,\,\dots, &
\text{$k_0$ odd,} \\
1,\,\ze^{-1},\,\ze,\,\ze^{-2},\,\ze^2,\ze^{-3},\,\dots, &
\text{$k_0$ even}
\end{cases}
\end{equation}
in $\Lt{+}$. \\ The Laurent polynomials
$\{p_-(\cdot,k,k_0)\}_{k\leq k_0}$ can be constructed by
Gram--Schmidt orthogonalizing
\begin{equation}
\begin{cases} 1,\,-\ze,\,\ze^{-1},\,-\ze^2,\,\ze^{-2},-\ze^3,\,\dots, &
\text{$k_0$ odd,} \\
-\ze,\,1,\,-\ze^2,\,\ze^{-1},\,-\ze^3,\,\ze^{-2},\dots, &
\text{$k_0$ even}
\end{cases}
\end{equation}
in $\Lt{-}$. \\ The Laurent polynomials
$\{r_-(\cdot,k,k_0)\}_{k\leq k_0}$ can be constructed by
Gram--Schmidt orthogonalizing
\begin{equation}
\begin{cases} -1,\,\ze^{-1},\,-\ze,\,\ze^{-2},\,-\ze^2,\ze^{-3},\,\dots, &
\text{$k_0$ odd,} \\
1,\,-\ze,\,\ze^{-1},\,-\ze^2,\,\ze^{-2},-\ze^3,\,\dots, &
\text{$k_0$ even}
\end{cases}
\end{equation}
in $\Lt{-}$.
\end{corollary}
%%%%%%%%%%%%%%%%%%%%%%%%%%%%%%%%%%%%%%%%%%%%%%%%%%%%%%%%%%%%%%%
\begin{proof}
The statements follow from Definition \ref{dA.2} and Lemma \ref{lA.4}.
\end{proof}
%%%%%%%%%%%%%%%%%%%%%%%%%%%%%%%%%%%%%%%%%%%%%%%%%%%%%%%%%%%%%%%

The following result clarifies which measures arise as spectral measures
of half-lattice CMV operators and it yields the reconstruction of
Verblunsky coefficients from the spectral measures and the corresponding
orthogonal polynomials.

%%%%%%%%%%%%%%%%%%%%%%%%%%%%%%%%%%%%%%%%%%%%%%%%%%%%%%%%%%%%%%%%%%
\begin{theorem} Let $k_0\in\Z$ and
$d\mu_\pm(\cdot,k_0)$ be nonnegative finite measures on
$\dD$ which are supported on infinite sets and normalized
by
\begin{align}
\oint_{\dD} d\mu_\pm(\ze,k_0) = 1.
\end{align}
Then $d\mu_\pm(\cdot,k_0)$ are necessarily the spectral
measures for some half-lattice CMV operators $U_{\pm,k_0}$ with
coefficients $\{\al_k\}_{k \geq k_0+1}$, respectively
$\{\al_k\}_{k \leq k_0}$, defined as follows,
\begin{equation}
\al_k = -
\begin{cases}
\big(p_+(\cdot,k-1,k_0),M_{\pm,k_0}(id)
r_+(\cdot,k-1,k_0)\big)_{\Lt{+}}, & k \text{ odd,} \\
\big(r_+(\cdot,k-1,k_0),p_+(\cdot,k-1,k_0)\big)_{\Lt{+}}, &
k \text{ even}
\end{cases}   \lb{A.202}
\end{equation}
for all $k \geq k_0+1$ and
\begin{equation}
\al_k = - \begin{cases}
\big(p_-(\cdot,k-1,k_0),M_{\pm,k_0}(id)
r_-(\cdot,k-1,k_0)\big)_{\Lt{-}}, & k \text{ odd,} \\
\big(r_-(\cdot,k-1,k_0),p_-(\cdot,k-1,k_0)\big)_{\Lt{-}}, &
k \text{ even} \end{cases} \label{A.203}
\end{equation}
for all $k \leq k_0$. Here the Laurent polynomials
$\{p_+(\cdot,k,k_0),r_+(\cdot,k,k_0)\}_{k\geq k_0}$ and
$\{p_-(\cdot,k,k_0),r_-(\cdot,k,k_0)\}_{k\leq k_0}$ denote
the orthonormal polynomials constructed in Corollary
\ref{cA.5}.
\end{theorem}
%%%%%%%%%%%%%%%%%%%%%%%%%%%%%%%%%%%%%%%%%%%%%%%%%%%%%%%%%%%%%%%%%
\begin{proof}
Using Corollary \ref{cA.5} one constructs the orthonormal
Laurent polynomials
$\{p_+(\ze,k,k_0),r_+(\ze,k,k_0)\}_{k\geq k_0}$, $\zeta\in\dD$. Because of
their orthogonality properties one concludes
\begin{equation}
r_+(\ze,k,k_0) = \begin{cases} \ze \ol{p_+(\ze,k,k_0)}, &
\text{$k_0$ odd,} \\
              \ol{p_+(\ze,k,k_0)}, & \text{$k_0$ even,} \end{cases}
\quad \ze\in\dD,\; k \geq k_0. \label{A.205}
\end{equation}
Next we will establish the recursion relation \eqref{A.41}.
Consider the following Laurent polynomial $p(\zeta)$, $\zeta\in\dD$, for
some fixed $k > k_0$,
\begin{equation}
p(\ze) = \begin{cases} \rho_k p_+(\ze,k,k_0) - \ze
r_+(\ze,k-1,k_0), & \text{$k$ odd}, \\ \rho_k
p_+(\ze,k,k_0) - r_+(\ze,k-1,k_0), & \text{$k$ even},
\end{cases} \quad \zeta\in\dD,
\end{equation}
where $\rho_k \in (0,\infty)$ is chosen such that the
leading term of $p_+(\cdot,k,k_0)$ cancels the leading term
of $r_+(\cdot,k-1,k_0)$. Using Corollary \ref{cA.5} one
checks that the Laurent polynomial $p(\cdot)$ is proportional
to $p_+(\cdot,k-1,k_0)$. Hence, one arrives at the following
recursion relation,
\begin{equation}
\rho_k p_+(\ze,k,k_0) = \begin{cases}\al_k p_+(\ze,k-1,k_0)
+ \ze r_+(\ze,k-1,k_0), & \text{$k$ odd}, \\ \ol{\al_k}
p_+(\ze,k-1,k_0) + r_+(\ze,k-1,k_0), & \text{$k$ even},
\end{cases} \quad \zeta\in\dD,
\end{equation}
where $\al_k \in \C$ is the proportionality constant.
Taking the scalar product of both sides with
$p_+(\ze,k-1,k_0)$ yields the expressions for $\alpha_k$,
$k\geq k_0 +1$, in \eqref{A.202}. Moreover, applying \eqref{A.205} one obtains
\begin{equation}
\rho_k r_+(\ze,k,k_0) = \begin{cases} \ol{\al_k}
r_+(\ze,k-1,k_0) + \frac{1}{\ze} p_+(\ze,k-1,k_0), &
\text{$k$ odd}, \\ \al_k r_+(\ze,k-1,k_0) +
p_+(\ze,k-1,k_0), & \text{$k$ even},
\end{cases} \quad \zeta\in\dD,
\end{equation}
and hence \eqref{A.41}. Since $\rho_k>0$, $k\in\Z$, it
remains to show that $\rho_k^2 = 1 - \abs{\al_k}^2$ and
hence that $\abs{\al_k}<1$. This follows from the
orthonormality of Laurent polynomials
$\{p_+(\cdot,k,k_0)\}_{k\geq k_0}$ in $\Lt{+}$,
\begin{align}
\abs{\al_k}^2 &= \norm{\al_k p_+(\cdot,k-1,k_0)}_{\Lt{+}}^2  \no \\
& = \norm{\rho_k p_+(\cdot,k,k_0) - id(\cdot)
r_+(\cdot,k-1,k_0)}_{\Lt{+}}^2 \no \\
&= \rho_k^2 + 1 - 2\Re\Big(\big(\rho_k p_+(\cdot,k,k_0), id(\cdot)
r_+(\cdot,k-1,k_0)\big)_{\Lt{+}}\Big)  \no
\\
&= \rho_k^2 + 1 \no \\
& \quad - 2\Re\Big(\big(\rho_k p_+(\cdot,k,k_0),[\rho_k
p_+(\cdot,k,k_0) - \al_k p_+(\cdot,k-1,k_0)] \big)_{\Lt{+}}\Big)  \no
\\
&= 1 - \rho_k^2, \, \text{ $k$ odd}.
\end{align}
Similarly one treats the case $k$ even.
Finally, using Lemma \ref{lA.2} one concludes that
$\Big(\begin{smallmatrix}p_+(z,k,k_0)\\r_+(z,k,k_0)
\end{smallmatrix}\Big)_{k\geq k_0}$,
$z\in\bbC\backslash\{0\}$, $k_0\in\Z$, is a generalized
eigenvector of the operator $\UU_{+,k_0}$ defined in
\eqref{2.33} associated with the coefficients
$\al_k,\rho_k$ introduced above. Thus, the measure
$d\mu_+(\cdot,k_0)$ is the spectral measure of the operator
$U_{+,k_0}$ in \eqref{2.31}.

Similarly one proves the result for $d\mu_-(\cdot,k_0)$ and
\eqref{A.203} for $k\leq k_0$.
\end{proof}
%%%%%%%%%%%%%%%%%%%%%%%%%%%%%%%%%%%%%%%%%%%%%%%%%%%%%%%%%%%%%%%

%%%%%%%%%%%%%%%%%%%%%%%%%%%%%%%%%%%%%%%%%%%%%%%%%%%%%%%%%%%%%%%
\begin{lemma} \label{lA.6}
Let $z\in\C\backslash(\dD\cup\{0\})$ and $k_0\in\bbZ$. Then
the sets of two-dimensional Laurent polynomials
$\Big(\begin{smallmatrix}\wti
p_\pm(z,k,k_0)\\r_\pm(z,k,k_0)
\end{smallmatrix}\Big)_{k\gtreqless
k_0}$ and $\Big(\begin{smallmatrix}\wti
q_\pm(z,k,k_0)\\s_\pm(z,k,k_0)
\end{smallmatrix}\Big)_{k\gtreqless
k_0}$ are related by,
\begin{align}
& \binom{\wti q_\pm(z,k,k_0)}{s_\pm(z,k,k_0)} = \pm
\oint_\dD d\mu_{\pm}(\zeta,k_0)\, \frac{\zeta+z}{\zeta-z}
\left( \binom{\wti p_\pm(\zeta,k,k_0)}{r_\pm(\zeta,k,k_0)}-
\binom{\wti p_\pm(z,k,k_0)}{r_\pm(z,k,k_0)} \right), \no \\
& \hspace*{9.5cm} k \gtrless k_0. \lb{2.102a}
\end{align}
\end{lemma}
%%%%%%%%%%%%%%%%%%%%%%%%%%%%%%%%%%%%%%%%%%%%%%%%%%%%%%%%%%%%%%%
\begin{proof}
First, we prove \eqref{2.102a} for $k_0$ even, which by
\eqref{2.99}--\eqref{2.102} is equivalent to
\begin{align}
&\binom{q_+(z,k,k_0)}{s_+(z,k,k_0)} = \oint_\dD
\frac{\ze+z}{\ze-z} \bigg(
\binom{p_+(\ze,k,k_0)}{r_+(\ze,k,k_0)} -
\binom{p_+(z,k,k_0)}{r_+(z,k,k_0)} \bigg)
d\mu_{+}(\ze,k_0),  \no \\[1mm]
& \hspace*{5.2cm} z\in\C\backslash(\dD\cup\{0\}),
\; k > k_0, \text{ $k_0$ even}.  \lb{2.103}
\end{align}
Let $k_0\in\Z$ be even. It suffices to show that the
right-hand side of \eqref{2.103}, temporarily denoted by
the symbol $RHS(z,k,k_0)$, satisfies
\begin{align}
&  T(z,k+1)^{-1} RHS(z,k+1,k_0) = RHS(z,k,k_0),
\quad k >  k_0, \lb{2.104} \\
&  T(z,k_0+1)^{-1}RHS(z,k_0+1,k_0) =
\binom{q_+(z,k_0,k_0)}{s_+(z,k_0,k_0)} = \binom{-1}{1}. \lb{2.105}
\end{align}
One verifies these statements using the following equality,
\begin{align}
& T(z,k+1)^{-1} RHS(z,k+1,k_0) = RHS(z,k,k_0) \no
\\ \no
& \quad + \oint_\dD \frac{\ze+z}{\ze-z} \left(
T(z,k+1)^{-1} - T(\ze,k+1)^{-1}\right)
\binom{p_+(\ze,k+1,k_0)}{r_+(\ze,k+1,k_0)}d\mu_{+}(\ze,k_0),
\\
& \hspace*{10cm} k\in\Z. \label{A.55}
\end{align}
For $k>k_0$, the last term on the right-hand side of
\eqref{A.55} is equal to zero since for $k$ odd, $T(z,k+1)$
does not depend on $z$, and for $k$ even, by Corollary
\ref{cA.5}, $p_+(\ze,k+1,k_0)$ and $r_+(\ze,k+1,k_0)$ are
orthogonal in $\Lt{+}$ to $\spn\{1,\ze\}$ and
$\spn\{1,\ze^{-1}\}$, respectively. Indeed,
\begin{align}
&\oint_\dD \frac{\ze+z}{\ze-z} \Big( T(z,k+1)^{-1}-
T(\ze,k+1)^{-1}\Big)
\binom{p_+(\ze,k+1,k_0)}{r_+(\ze,k+1,k_0)}d\mu_{+}(\ze,k_0)
\no
\\
& \quad = \oint_\dD
\frac{\ze+z}{\ze-z}\;\frac{1}{\rho_{k+1}}
\begin{pmatrix}0 & z-\ze \\ (1/z) - (1/\ze) & 0\end{pmatrix}
\binom{p_+(\ze,k+1,k_0)}{r_+(\ze,k+1,k_0)}d\mu_{+}(\ze,k_0)
\no
\\
& \quad = \frac{1}{\rho_{k+1}}\oint_\dD
\begin{pmatrix}0 & -(\ze+z) \\ (1/\ze) + (1/z) & 0\end{pmatrix}
\binom{p_+(\ze,k+1,k_0)}{r_+(\ze,k+1,k_0)}d\mu_{+}(\ze,k_0)
\no
\\ & \quad = \frac{1}{\rho_{k+1}}\oint_\dD \binom
{-\ol{((1/\ze)+\ol{z})}r_+(\ze,k,k_0)}
{\ol{(\ze+(1/\ol{z}))}p_+(\ze,k,k_0)} d\mu_{+}(\ze,k_0) =
\binom{0}{0}.
\end{align}
This proves \eqref{2.104}.

For $k=k_0$ one obtains $RHS(z,k_0,k_0)=0$ since
$p_+(\ze,k_0,k_0)=r_+(\ze,k_0,k_0)=1$. By Corollary
\ref{cA.5}, $p_+(\ze,k_0+1,k_0)$ and $r_+(\ze,k_0+1,k_0)$
are orthogonal to constants in $\Lt{+}$ and by the
recursion relation \eqref{A.28},
\begin{align}
p_+(\ze,k_0+1,k_0) = (\ze+\al_{k_0+1})/\rho_{k_0+1}, \quad
r_+(\ze,k_0+1,k_0) = ((1/\ze)+\ol{\al_{k_0+1}})/\rho_{k_0+1}.
\end{align}
Thus,
\begin{align}
&\oint_\dD \frac{\ze+z}{\ze-z} \Big( T(z,k_0+1)^{-1}-
T(\ze,k_0+1)^{-1}\Big)
\binom{p_+(\ze,k_0+1,k_0)}{r_+(\ze,k_0+1,k_0)}d\mu_{+}(\ze,k_0)
\no
\\
& \quad = \oint_\dD \frac{1}{\rho_{k_0+1}} \binom
{-\ol{((1/\ze)+\ol{z})},r_+(\ze,k_0+1,k_0)}
{\ol{(\ze+(1/\ol{z}))},p_+(\ze,k_0+1,k_0)}
d\mu_{+}(\ze,k_0) \no
\\ & \quad = \binom {-\norm{r_+(\ze,k_0+1,k_0)}^2_{\Lt{+}}} {
\norm{p_+(\ze,k_0+1,k_0)}^2_{\Lt{+}}} = \binom{-1}{1}.
\end{align}
This proves \eqref{2.105}.

Next, we prove that
\begin{align}
&\binom{s_+(z,k,k_0)}{\wti q_+(z,k,k_0)} = \oint_\dD
\frac{\ze+z}{\ze-z} \bigg( \binom{r_+(\ze,k,k_0)}{\wti
p_+(\ze,k,k_0)} - \binom{r_+(z,k,k_0)}{\wti p_+(z,k,k_0)}
\bigg) d\mu_{+}(\ze,k_0),   \no \\[1mm]
& \hspace*{5.4cm} z \in \C\backslash(\dD\cup\{0\}),
\; k > k_0, \text{ $k_0$ odd}.  \lb{2.109}
\end{align}
Let $k_0\in\Z$ be odd. We note that
\begin{align}
\binom{u(z,k)}{v(z,k)} &=  T(z,k) \binom{u(z,k-1)}{v(z,k-1)} \no
\\
\intertext{is equivalent to }
\binom{v(z,k)}{\wti u(z,k)} &= \wti T(z,k)
\binom{v(z,k-1)}{\wti u(z,k-1)},
\end{align}
where
\begin{equation}
\wti u(z,k) = u(z,k)/z, \quad \wti T(z,k) =
\begin{pmatrix}0 & 1 \\ 1/z & 0\end{pmatrix} T(z,k)
\begin{pmatrix}0 & z \\ 1 & 0\end{pmatrix}.
\end{equation}
Thus, it suffices to show that the right-hand side of
\eqref{2.109}, temporarily denoted by $\wti{RHS}(z,k,k_0)$,
satisfies
\begin{align}
& \wti T(z,k+1)^{-1} \wti{RHS} (z,k+1,k_0) =
\wti{RHS}(z,k,k_0), \quad k > k_0, \lb{2.113}
\\
& \wti T(z,k_0+1)^{-1} \wti{RHS}(z,k_0+1,k_0) =
\binom{s_+(z,k_0,k_0)}{\wti q_+(z,k_0,k_0)} =
\binom{-1}{1}. \lb{2.114}
\end{align}
At this point one can follow the first part of the proof
replacing $T$ by $\widetilde T$, $\Big(\begin{smallmatrix}
p_+ \\[1mm] r_+\end{smallmatrix}\Big)$ by
$\Big(\begin{smallmatrix} r_+ \\ \wti
p_+\end{smallmatrix}\Big)$, $\Big(\begin{smallmatrix} q_+
\\[1mm] s_+\end{smallmatrix}\Big)$ by $\Big(\begin{smallmatrix}
s_+ \\ \wti q_+\end{smallmatrix}\Big)$, etc.

The result for the remaining polynomials $\wti
p_-(z,k,k_0)$, $r_-(z,k,k_0)$, $\wti q_-(z,k,k_0)$, and
$s_-(z,k,k_0)$ follows similarly.
\end{proof}
%%%%%%%%%%%%%%%%%%%%%%%%%%%%%%%%%%%%%%%%%%%%%%%%%%%%%%%%%%%%%%%

%%%%%%%%%%%%%%%%%%%%%%%%%%%%%%%%%%%%%%%%%%%%%%%%%%%%%%%%%%%%%%%
\begin{corollary} \label{cA.7}
Let $k_0\in\bbZ$. Then the sets of two-dimensional Laurent
polynomials
$\Big(\begin{smallmatrix}p_\pm(z,k,k_0)\\
r_\pm(z,k,k_0)\end{smallmatrix}\Big)_{k\gtreqless
k_0}$ and $\Big(\begin{smallmatrix}q_\pm(z,k,k_0)\\
s_\pm(z,k,k_0)\end{smallmatrix}\Big)_{k\gtreqless
k_0}$ satisfy the relation
\begin{align} \label{A.63}
\binom{q_\pm(z,\cdot,k_0)}{s_\pm(z,\cdot,k_0)} +
m_\pm(z,k_0)
\binom{p_\pm(z,\cdot,k_0)}{r_\pm(z,\cdot,k_0)} \in
\ell^2([k_0,\pm\infty)\cap\Z)^2, \quad
z\in\bbC\backslash(\dD\cup\{0\}),
\end{align}
for some coefficients $m_\pm(z,k_0)$ given by
\begin{align}
m_\pm(z,k_0) &= \pm
(\delta_{k_0},(U_{\pm,k_0}+zI)(U_{\pm,k_0}-zI)^{-1}
\delta_{k_0})_{\ell^2([k_0,\pm\infty)\cap\Z)}
\lb{A.66} \\
& =\pm \oint_\dD d\mu_{\pm}(\zeta,k_0)\,
\frac{\zeta+z}{\zeta-z}, \quad z\in\bbC\backslash\dD \lb{A.67}
\intertext{with}
m_\pm(0,k_0)&=\pm\oint_{\dD} d\mu_\pm(\zeta,k_0)=\pm 1.
\end{align}
\end{corollary}
%%%%%%%%%%%%%%%%%%%%%%%%%%%%%%%%%%%%%%%%%%%%%%%%%%%%%%%%%%%%%%%
\begin{proof}
Consider the operator
\begin{align}
& C_{\pm,k_0}(z) = \begin{cases}
\begin{pmatrix}
I & 0 \\ 0 & \pm I
\end{pmatrix}
((\UU_{\pm,k_0})^\top+zI)((\UU_{\pm,k_0})^\top-zI)^{-1}, &
k_0 \text{ odd}, \\
\begin{pmatrix}
\pm I & 0 \\ 0 & I
\end{pmatrix}
((\UU_{\pm,k_0})^\top+zI)((\UU_{\pm,k_0})^\top-zI)^{-1}, &
k_0 \text{ even}, \end{cases}
\\
& \hspace*{8.4cm}  z\in\bbC\backslash\dD, \no
\end{align}
on $\ltz^2$. Since $C_{\pm,k_0}(z)$ is bounded for
$z\in\bbC\backslash\dD$ one has
\begin{align}
\left\{\left(
\binom{\de_{k_0}}{\de_{k_0}},C_{\pm,k_0}(z)\binom{\de_k}{\de_k}
\right)\right\}_{k\in\bbZ}
=
\left\{\left(
C_{\pm,k_0}(z)^*\binom{\de_{k_0}}{\de_{k_0}},\binom{\de_k}{\de_k}
\right)\right\}_{k\in\bbZ} \in\ltz^2.
\end{align}
Using the spectral representation for the operator
$C_{\pm,k_0}(z)$, Lemma \ref{lA.6}, and
\eqref{2.99}--\eqref{2.102} one obtains
\begin{align}
& \bigg( \binom{\de_{k_0}}{\de_{k_0}},
C_{\pm,k_0}(z)\binom{\de_k}{\de_k} \bigg)
=
\oint_\dD d\mu_{\pm}(\ze,k_0)\, \frac{\ze+z}{\ze-z}
\binom{\wti p_\pm(\ze,k,k_0)}{r_\pm(\ze,k,k_0)} \no \\ &
\quad = \pm\left[ \binom{\wti q_\pm(z,k,k_0)}{s_\pm(z,k,k_0)}
+ m_\pm(z,k_0) \binom{\wti
p_\pm(z,k,k_0)}{r_\pm(z,k,k_0)}\right], \quad k \gtrless k_0,
\end{align}
where $m_\pm(z,k_0) = \pm \int_\dD
d\mu_{\pm}(\ze,k_0)\,\frac{\ze+z}{\ze-z}$.
\end{proof}
%%%%%%%%%%%%%%%%%%%%%%%%%%%%%%%%%%%%%%%%%%%%%%%%%%%%%%%%%%%%%%%

%%%%%%%%%%%%%%%%%%%%%%%%%%%%%%%%%%%%%%%%%%%%%%%%%%%%%%%%%%%%%%%
\begin{lemma} \label{lA.8}
Let $k_0\in\bbZ$. Then relation \eqref{A.63} uniquely
determines the functions $m_\pm(\cdot,k_0)$ on
$\bbC\backslash\dD$.
\end{lemma}
%%%%%%%%%%%%%%%%%%%%%%%%%%%%%%%%%%%%%%%%%%%%%%%%%%%%%%%%%%%%%%%
\begin{proof}
We will prove the lemma by contradiction. Assume that there
are two functions $m_+(z,k_0)$ and $\wti m_+(z,k_0)$
satisfying \eqref{A.63} such that $m_+(z_0,k_0) \neq \wti
m_+(z_0,k_0)$ for some $z_0\in\C\backslash\dD$. Then there
are $\la_1,\la_2 \in\C$ such that the following vector
\begin{align}
\binom{w_1(z_0,\cdot,k_0)}{w_2(z_0,\cdot,k_0)} &= (\la_1
m_+(z_0,k_0) + \la_2 \wti m_+(z_0,k_0))
\binom{p_+(z_0,\cdot,k_0)}{r_+(z_0,\cdot,k_0)}
\\ & \quad +
(\la_1+\la_2)
\binom{q_+(z_0,\cdot,k_0)}{s_+(z_0,\cdot,k_0)}
\in \ell^2([k_0,\infty)\cap\Z)^2
\end{align}
is nonzero and satisfies
\begin{equation}
w_1(z_0,k_0,k_0) = \begin{cases} z_0 w_2(z_0,k_0,k_0), &
\text{$k_0$ odd,}
\\
w_2(z_0,k_0,k_0), & \text{$k_0$ even.}\end{cases}
\end{equation}
By Lemma \ref{lA.2},
$\Big(\begin{smallmatrix}w_1(z_0,k,k_0)\\
w_2(z_0,k,k_0)\end{smallmatrix}\Big)_{k\geq k_0}$ is an
eigenvector of the operator $\UU_{+,k_0}$ and
$z_0\in\C\backslash\dD$ is the corresponding eigenvalue
which is impossible since $\UU_{+,k_0}$ is unitary.

Similarly, one proves the result for $m_-(z,k_0)$.
\end{proof}
%%%%%%%%%%%%%%%%%%%%%%%%%%%%%%%%%%%%%%%%%%%%%%%%%%%%%%%%%%%%%%%

%%%%%%%%%%%%%%%%%%%%%%%%%%%%%%%%%%%%%%%%%%%%%%%%%%%%%%%%%%%%%%%
\begin{corollary} \label{cA.9}
There are solutions $\Big(\begin{smallmatrix}\psi_\pm(z,k)\\
\chi_\pm(z,k)\end{smallmatrix}\Big)_{k\in\Z}$ of
\eqref{A.28}, unique up to constant multiples, so that for
some $($and hence for all\,$)$ $k_1\in\Z$,
\begin{align}
\binom{\psi_\pm(z,\cdot)}{\chi_\pm(z,\cdot)} \in
\ell^2([k_1,\pm\infty)\cap\Z)^2, \quad
z\in\C\backslash(\dD\cup\{0\}).
\end{align}
\end{corollary}
%%%%%%%%%%%%%%%%%%%%%%%%%%%%%%%%%%%%%%%%%%%%%%%%%%%%%%%%%%%%%%%
\begin{proof}
Since any solution of \eqref{A.28} can be expressed as a
linear combination of the polynomials
$\Big(\begin{smallmatrix}p_\pm(z,k,k_0)\\
r_\pm(z,k,k_0)\end{smallmatrix}\Big)_{k\in\Z}$ and
$\Big(\begin{smallmatrix}q_\pm(z,k,k_0)\\
s_\pm(z,k,k_0)\end{smallmatrix}\Big)_{k\in\Z}$, existence and
uniqueness of the solutions
$\Big(\begin{smallmatrix}\psi_\pm(z,\cdot)\\
\chi_\pm(z,\cdot)\end{smallmatrix}\Big)_{k\in\Z}$ follow from
Corollary \ref{cA.7} and Lemma \ref{lA.8}, respectively.
\end{proof}
%%%%%%%%%%%%%%%%%%%%%%%%%%%%%%%%%%%%%%%%%%%%%%%%%%%%%%%%%%%%%%%

%%%%%%%%%%%%%%%%%%%%%%%%%%%%%%%%%%%%%%%%%%%%%%%%%%%%%%%%%%%%%%%
\begin{lemma} \label{lA.9}
Let $z\in\bbC\backslash\{0\}$ and $k_0\in\bbZ$. Then the
two-dimensional Laurent polynomials
$\Big(\begin{smallmatrix}p_+(z,k,k_0)\\
r_+(z,k,k_0)\end{smallmatrix}\Big)_{k\in\bbZ}$,
$\Big(\begin{smallmatrix}q_+(z,k,k_0)\\
s_+(z,k,k_0)\end{smallmatrix}\Big)_{k\in\bbZ}$,
$\Big(\begin{smallmatrix}p_-(z,k,k_0-1)\\
r_-(z,k,k_0-1)\end{smallmatrix}\Big)_{k\in\bbZ}$,
$\Big(\begin{smallmatrix}q_-(z,k,k_0-1)\\
s_-(z,k,k_0-1)\end{smallmatrix}\Big)_{k\in\bbZ}$ are
connected by the following relations:
\begin{align}
\binom{p_-(z,k,k_0-1)}{r_-(z,k,k_0-1)} &= \frac{i\Im(
b_{k_0})}{\rho_{k_0}} \binom{p_+(z,k,k_0)}{r_+(z,k,k_0)} +
\frac{\Re( b_{k_0})}{\rho_{k_0}}
\binom{q_+(z,k,k_0)}{s_+(z,k,k_0)}, \lb{A.67b}
\\
\binom{q_-(z,k,k_0-1)}{s_-(z,k,k_0-1)} &= \frac{\Re
(a_{k_0})}{\rho_{k_0}} \binom{p_+(z,k,k_0)}{r_+(z,k,k_0)} +
\frac{i\Im (a_{k_0})}{\rho_{k_0}}
\binom{q_+(z,k,k_0)}{s_+(z,k,k_0)}, \quad k\in\Z.
\label{A.67a}
\end{align}
\end{lemma}
%%%%%%%%%%%%%%%%%%%%%%%%%%%%%%%%%%%%%%%%%%%%%%%%%%%%%%%%%%%%%%%
\begin{proof}
It follows from Definition \ref{dA.2} that the left- and
right-hand sides of \eqref{A.67b} and \eqref{A.67a} satisfy
the same recursion relation \eqref{A.28}. Hence, it
suffices to check \eqref{A.67b} and \eqref{A.67a} at one
point, say, the point $k=k_0$. Using \eqref{A.17},
\eqref{A.18}, \eqref{A.28}, and \eqref{A.50}, one finds the
following expressions for the left-hand sides of
\eqref{A.67b} and \eqref{A.67a},
\begin{align}
& \binom{p_-(z,k_0,k_0-1)}{r_-(z,k_0,k_0-1)} =
\frac{1}{\rho_{k_0}} \binom{zb_{k_0}}{-\ol{b_{k_0}}}, \quad
\binom{q_-(z,k_0,k_0-1)}{s_-(z,k_0,k_0-1)} =
\frac{1}{\rho_{k_0}} \binom{za_{k_0}}{\ol{a_{k_0}}}, \\ &
\hspace*{9.9cm} k_0 \text{ odd} \no
\end{align}
and
\begin{align}
& \binom{p_-(z,k_0,k_0-1)}{r_-(z,k_0,k_0-1)} =
\frac{1}{\rho_{k_0}} \binom{-\ol{b_{k_0}}}{b_{k_0}}, \quad
\binom{q_-(z,k_0,k_0-1)}{s_-(z,k_0,k_0-1)} =
\frac{1}{\rho_{k_0}} \binom{\ol{a_{k_0}}}{a_{k_0}}, \\ &
\hspace*{9.65cm} k_0 \text{ even}. \no
\end{align}
The same result also follows for the right-hand side of
\eqref{A.67b}, \eqref{A.67a} using \eqref{A.17},
\eqref{A.18}, and the initial conditions \eqref{A.48}.
\end{proof}
%%%%%%%%%%%%%%%%%%%%%%%%%%%%%%%%%%%%%%%%%%%%%%%%%%%%%%%%%%%%%%%

%%%%%%%%%%%%%%%%%%%%%%%%%%%%%%%%%%%%%%%%%%%%%%%%%%%%%%%%%%%%%%%
\begin{theorem} \lb{tA.10}
Let $k_0\in\bbZ$. Then there exist
unique functions $M_\pm(\cdot,k_0)$ such that
\begin{align}
&\binom{u_\pm(z,\cdot,k_0)}{v_\pm(z,\cdot,k_0)} =
\binom{q_+(z,\cdot,k_0)}{s_+(z,\cdot,k_0)} + M_\pm(z,k_0)
\binom{p_+(z,\cdot,k_0)}{r_+(z,\cdot,k_0)} \in
\ell^2([k_0,\pm\infty)\cap\Z)^2, \no \\
& \hspace*{8cm} z\in\bbC\backslash(\dD\cup\{0\}). \label{A.68}
\end{align}
\end{theorem}
%%%%%%%%%%%%%%%%%%%%%%%%%%%%%%%%%%%%%%%%%%%%%%%%%%%%%%%%%%%%%%%
\begin{proof}
Assertion \eqref{A.68} follows from
\eqref{2.99}--\eqref{2.102}, Corollaries \ref{cA.7} and
\ref{cA.9}, and Lemmas \ref{lA.8} and \ref{lA.9}.
\end{proof}
%%%%%%%%%%%%%%%%%%%%%%%%%%%%%%%%%%%%%%%%%%%%%%%%%%%%%%%%%%%%%%%

We will call $u_\pm(z,\cdot,k_0)$ (resp.,
$v_\pm(z,\cdot,k_0)$) {\it Weyl--Titchmarsh solutions}
of $U$ (resp., $U^\top$). By Corollary \ref{cA.9},
$u_\pm(z,\cdot,k_0)$ and $v_\pm(z,\cdot,k_0)$ are constant
multiples of $\psi_\pm(z,\cdot,k_0)$ and
$\chi_\pm(z,\cdot,k_0)$. Similarly, we will call
$m_\pm(z,k_0)$ as well as $M_\pm(z,k_0)$ the {\it half-lattice
Weyl--Titchmarsh $m$-functions} associated with
$U_{\pm,k_0}$. (See also \cite{Si04a} for a comparison of
various alternative notions of Weyl--Titchmarsh $m$-functions for
$U_{+,k_0}$.)

It follows from Corollary \ref{cA.7} and \ref{cA.9} and
Lemma \ref{lA.9} that
\begin{align}
M_+(z,k_0) &= m_+(z,k_0), \quad z\in\bbC\backslash\dD,
\lb{A.69}
\\ M_+(0,k_0) &=1, \lb{A.70} \\ M_-(z,k_0) &=
\frac{\Re(a_{k_0}) +
i\Im(b_{k_0})m_-(z,k_0-1)}{i\Im(a_{k_0})
+ \Re(b_{k_0})m_-(z,k_0-1)}, \quad z\in\bbC\backslash\dD,
\lb{A.71} \\ M_-(0,k_0)
&=\f{\alpha_{k_0}+1}{\alpha_{k_0}-1}.
\lb{A.72}
\end{align}
In particular, one infers that $M_\pm$ are analytic at
$z=0$.

Since \eqref{A.68} singles out $p_+(z,\cdot,k_0)$,
$q_+(z,\cdot,k_0)$, $r_+(z,\cdot,k_0)$, and
$s_+(z,\cdot,k_0)$, we now add the following observation.

%%%%%%%%%%%%%%%%%%%%%%%%%%%%%%%%%%%%%%%%%%%%%%%%%%%%%%%%%%%%%%%
\begin{remark} \label{rA.13}
One can also define functions $\hatt M_\pm(\cdot,k_0)$ such
that the following relation holds
\begin{align}
& \binom{\hatt u_\pm (z,\cdot,k_0)}{\hatt v_\pm
(z,\cdot,k_0)}= \binom{q_-(z,\cdot,k_0)}{s_-(z,\cdot,k_0)}
+ \hatt M_\pm(z,k_0)
\binom{p_-(z,\cdot,k_0)}{r_-(z,\cdot,k_0)} \in
\ell^2([k_0,\pm\infty)\cap\bbZ)^2,  \no \\ & \hspace*{8cm}
z\in\bbC\backslash(\dD\cup\{0\}). \lb{A.73}
\end{align}
Applying Corollary \ref{cA.9}, $\hatt
u_\pm(z,\cdot,k_0)$ and
$\hatt v_\pm(z,\cdot,k_0)$ are also constant multiples of
$\psi_\pm(z,\cdot,k_0)$ and $\chi_\pm(z,\cdot,k_0)$ (hence
they are constant multiples of $u_\pm (z,\cdot,k_0)$ and
$v_\pm (z,\cdot,k_0)$). It follows from Corollaries \ref{cA.7}
and \ref{cA.9} and Lemmas \ref{lA.8} and \ref{lA.9}, that
$\hatt M_\pm(\cdot,k_0)$ are uniquely defined and
satisfy the
relations
\begin{align}
\hatt M_+(z,k_0-1) &= \frac{\Re(a_{k_0}) -
i\Im(a_{k_0})m_+(z,k_0)} {-i\Im(b_{k_0}) +
\Re(b_{k_0})m_+(z,k_0)}, \quad z\in\bbC\backslash\dD,
\lb{A.74} \\ \hatt M_-(z,k_0) &= m_-(z,k_0), \quad
z\in\bbC\backslash\dD. \lb{A.75}
\end{align}
Moreover, one derives from \eqref{A.71} and
\eqref{A.75} that
\begin{align} \label{A.76}
M_\pm(z,k_0) &= \frac{\Re(a_{k_0}) + i\Im(b_{k_0})\hatt
M_\pm(z,k_0-1)} {i\Im(a_{k_0}) + \Re(b_{k_0})\hatt
M_\pm(z,k_0-1)}, \quad z\in\bbC\backslash\dD.
\end{align}
In this paper we will only use
$\Big(\begin{smallmatrix}u_\pm(z,\cdot,k_0)
\\ v_\pm(z,\cdot,k_0)\end{smallmatrix}\Big)$ and
$M_\pm(z,k_0)$.
\end{remark}
%%%%%%%%%%%%%%%%%%%%%%%%%%%%%%%%%%%%%%%%%%%%%%%%%%%%%%%%%%%%%%%

%%%%%%%%%%%%%%%%%%%%%%%%%%%%%%%%%%%%%%%%%%%%%%%%%%%%%%%%%%%%%%%
\begin{lemma} \label{lA.16}
Let $k\in\bbZ$. Then the functions $M_+(\cdot,k)|_{\D}$
$($resp., $M_-(\cdot,k)|_{\D}$$)$ are Caratheodory
$($resp.,
anti-Caratheodory\,$)$ functions. Moreover, $M_\pm$ satisfy
the following Riccati-type equation
\begin{align}
&(z\ol{b_k}-b_k)M_\pm(z,k-1)M_\pm(z,k)+(z\ol{b_k}+b_k)M_\pm(z,k)
-(z\ol{a_k}+a_k)M_\pm(z,k-1) \no \\
& \quad =z\ol{a_k}-a_k, \quad z\in\bbC\backslash\dD. \lb{2.148}
\end{align}
\end{lemma}
%%%%%%%%%%%%%%%%%%%%%%%%%%%%%%%%%%%%%%%%%%%%%%%%%%%%%%%%%%%%%%%
\begin{proof}
It follows from \eqref{A.67} and Theorem \ref{tA.2} that
$m_\pm(z,k_0)$ are Caratheodory and anti-Caratheodory
functions, respectively. From \eqref{A.69} one concludes
that $M_+(z,k_0)$ is also a Caratheodory function. Using
\eqref{A.71} one verifies that $M_-(z,k_0)$ is analytic in
$\D$ since $\Re(m_-(z,k_0)) < 0$ and that
\begin{align}
\Re(M_-(z,k_0)) &= \Re\left( \frac{\Re(a_{k_0}) +
i\Im(b_{k_0})m_-(z,k_0-1)} {i\Im(a_{k_0}) +
\Re(b_{k_0})m_-(z,k_0-1)}\right)  \no
\\
&= \frac{\Re(a_{k_0})\Re(b_{k_0})+\Im(a_{k_0})\Im(b_{k_0})}
{\abs{i\Im(a_{k_0}) + \Re(b_{k_0})m_-(z,k_0-1)}^2}
\Re(m_-(z,k_0-1))  \no
\\
&= \frac{\rho_{k_0}^2 \Re(m_-(z,k_0-1))}
{\abs{i\Im(a_{k_0})
+ \Re(b_{k_0})m_-(z,k_0-1)}^2} < 0.
\end{align}
Hence, $M_-(z,k_0)$ is an anti-Caratheodory function.

Next, consider the $2\times 2$ matrix
\begin{align}
D(z,k_0) = \big(d_{\ell,\ell'}(z,k_0)\big)_{\ell,\ell'=1,2}
&= \frac{1}{2\rho_{k_0}}
\begin{cases}
\begin{pmatrix}
\ol{a_{k_0}} + a_{k_0}/z & \ol{a_{k_0}} - a_{k_0}/z \\
\ol{b_{k_0}} - b_{k_0}/z & \ol{b_{k_0}} + b_{k_0}/z
\end{pmatrix}, & k_0 \text{ odd,}
\\
\begin{pmatrix}
z\ol{a_{k_0}} + a_{k_0} & z\ol{a_{k_0}} - a_{k_0} \\
z\ol{b_{k_0}} - b_{k_0} & z\ol{b_{k_0}} + b_{k_0}
\end{pmatrix}, & k_0 \text{ even,}
\end{cases} \no
\\
& \hspace*{33mm} z\in\Cz, \; k_0\in\Z.
\end{align}
It follows from \eqref{A.17}, \eqref{A.18}, and Definition
\ref{dA.2} that $D(z,k_0)$ satisfies
\begin{align}
\begin{pmatrix}
p_+(z,\cdot,k_0-1) & q_+(z,\cdot,k_0-1)
\\
r_+(z,\cdot,k_0-1) & s_+(z,\cdot,k_0-1)
\end{pmatrix}
=
\begin{pmatrix}
p_+(z,\cdot,k_0) & q_+(z,\cdot,k_0)
\\
r_+(z,\cdot,k_0) & s_+(z,\cdot,k_0)
\end{pmatrix}
D(z,k_0).
\end{align}
Thus, using Theorem \ref{tA.10} one finds
\begin{align}
M_\pm(z,k_0) =
\frac{d_{1,2}(z,k_0)+d_{1,1}(z,k_0)M_\pm(z,k_0-1)}
{d_{2,2}(z,k_0)+d_{2,1}(z,k_0)M_\pm(z,k_0-1)}.
\end{align}
\end{proof}
%%%%%%%%%%%%%%%%%%%%%%%%%%%%%%%%%%%%%%%%%%%%%%%%%%%%%%%%%%%%%%%

In addition, we introduce the functions
$\Phi_\pm(\cdot,k)$, $k\in\bbZ$, by
\begin{align}
\Phi_\pm(z,k) = \frac{M_\pm(z,k)-1}{M_\pm(z,k)+1}, \quad
z\in\C\backslash\dD. \lb{A.78}
\end{align}
One then verifies,
\begin{equation}
M_\pm(z,k) = \frac{1+\Phi_\pm(z,k)}{1-\Phi_\pm(z,k)}, \quad
z\in\C\backslash\dD. \lb{A.79}
\end{equation}
Moreover, we extend these functions to the unit circle
$\dD$ by taking the radial limits which exist and are finite
for $\mu_0$-almost every $\ze\in\dD$,
\begin{align}
M_\pm(\ze,k) &= \lim_{r \uparrow 1} M_\pm(r\ze,k),
\\
\Phi_\pm(\ze,k) &= \lim_{r \uparrow 1} \Phi_\pm(r\ze,k),
\quad k\in\Z.
\end{align}

%%%%%%%%%%%%%%%%%%%%%%%%%%%%%%%%%%%%%%%%%%%%%%%%%%%%%%%%
\begin{lemma} \label{l2.19}
Let $z\in\C\backslash(\dD\cup\{0\})$, $k_0, k\in\Z$. Then
the functions $\Phi_\pm(\cdot,k)$ satisfy
\begin{equation}
\Phi_\pm(z,k) = \begin{cases}
z\frac{v_\pm(z,k,k_0)}{u_\pm(z,k,k_0)}, &\text{$k$ odd,}
\\
\frac{u_\pm(z,k,k_0)}{v_\pm(z,k,k_0)}, & \text{$k$ even,}
\end{cases}
\end{equation}
where $u_\pm(\cdot,k,k_0)$ and $v_\pm(\cdot,k,k_0)$ are the
polynomials defined in \eqref{A.68}.
\end{lemma}
%%%%%%%%%%%%%%%%%%%%%%%%%%%%%%%%%%%%%%%%%%%%%%%%%%%%%%%%
\begin{proof}
Using Corollary \ref{cA.9} it suffices to assume $k=k_0$. Then the
statement follows immediately from \eqref{A.48} and
\eqref{A.78}.
\end{proof}
%%%%%%%%%%%%%%%%%%%%%%%%%%%%%%%%%%%%%%%%%%%%%%%%%%%%%%%%

%%%%%%%%%%%%%%%%%%%%%%%%%%%%%%%%%%%%%%%%%%%%%%%%%%%%%%%%%%%%%%%
\begin{lemma}
Let $k\in\bbZ$. Then the functions $\Phi_+(\cdot,k)|_{\D}$
$($resp., $\Phi_-(\cdot,k)|_{\D}$$)$ are Schur $($resp.,
anti-Schur\,$)$ functions. Moreover, $\Phi_\pm$ satisfy the
following Riccati-type equation
\begin{equation}
\alpha_k
\Phi_\pm(z,k-1)\Phi_\pm(z,k)-\Phi_\pm(z,k-1)+z\Phi_\pm(z,k)
=\ol{\alpha_k}z, \quad z\in\bbC\backslash\dD,\; k\in\Z.
\lb{A.80a}
\end{equation}
\end{lemma}
%%%%%%%%%%%%%%%%%%%%%%%%%%%%%%%%%%%%%%%%%%%%%%%%%%%%%%%%%%%%%%%
\begin{proof}
It follows from Lemma \ref{lA.16} and \eqref{A.78} that the
functions $\Phi_+(\cdot,k)|_{\D}$ $($resp.,
$\Phi_-(\cdot,k)|_{\D}$$)$ are Schur $($resp.,
anti-Schur\,$)$ functions.

Let $k$ be odd. Then applying Lemma \ref{l2.19}
and the recursion relation \eqref{A.28} one obtains
\begin{align}
\Phi_\pm(z,k) &= \frac{zv_\pm(z,k,k_0)}{u_\pm(z,k,k_0)} =
\frac{u_\pm(z,k-1,k_0)+z\ol{\al_k}v_\pm(z,k-1,k_0)}
{\al_{k}u_\pm(z,k-1,k_0)+zv_\pm(z,k-1,k_0)}  \no
\\ &=
\frac{\Phi_\pm(z,k-1)+z\ol{\al_k}}{\al_k\Phi_\pm(z,k-1)+z}.
\end{align}
For $k$ even, one similarly obtains
\begin{align}
\Phi_\pm(z,k) &= \frac{u_\pm(z,k,k_0)}{v_\pm(z,k,k_0)} =
\frac{\ol{\al_k}u_\pm(z,k-1,k_0)+v_\pm(z,k-1,k_0)}
{u_\pm(z,k-1,k_0)+\al_{k}v_\pm(z,k-1,k_0)}  \no
\\ &=
\frac{z\ol{\al_k}+\Phi_\pm(z,k-1)}{z+\al_k\Phi_\pm(z,k-1)}.
\end{align}
\end{proof}
%%%%%%%%%%%%%%%%%%%%%%%%%%%%%%%%%%%%%%%%%%%%%%%%%%%%%%%%%%%%%%%

%%%%%%%%%%%%%%%%%%%%%%%%%%%%%%%%%%%%%%%%%%%%%%%%%%%%%%%%%%%%%%%
\begin{remark} \lb{r2.18}
$(i)$ In the special case $\alpha=\{\alpha_k\}_{k\in\Z}=0$, one obtains
\begin{equation}
M_\pm(z,k) = \pm 1, \quad \Phi_+(z,k)=0, \quad 1/\Phi_-(z,k)=0,
\quad  z\in\C, \; k\in\Z.
\end{equation}
Thus, strictly speaking, one should always consider
$1/\Phi_-$ rather than $\Phi_-$ and hence refer to the
Riccati-type equation of $1/\Phi_-$,
\begin{equation}
\ol{\alpha_k}z\f{1}{\Phi_-(z,k-1)}\f{1}{\Phi_-(z,k)}
+\f{1}{\Phi_-(z,k)} -z\f{1}{\Phi_-(z,k-1)}=\alpha_k, \quad
z\in\C\backslash\dD, \; k\in\bbZ,  \lb{2.163}
\end{equation}
rather than that of $\Phi_-$, etc. For simplicity of
notation, we will avoid this distinction between
$\Phi_-$ and $1/\Phi_-$ and usually just invoke $\Phi_-$
whenever confusions are unlikely. \\ $(ii)$ We note that
$M_\pm(z,k)$ and $\Phi_\pm(z,k)$, $z\in\dD$, $k\in\bbZ$,
have nontangential limits to $\dD$ $\mu_0$-a.e.\ In
particular, the Riccati-type equations
\eqref{2.148}, \eqref{A.80a}, and
\eqref{2.163} extend to $\dD$ $\mu_0$-a.e.
\end{remark}
%%%%%%%%%%%%%%%%%%%%%%%%%%%%%%%%%%%%%%%%%%%%%%%%%%%%%%%%%%%%%%%

The Riccati-type equation for the Schur function
$\Phi_+$ implies the following absolutely convergent
expansion,
\begin{align}
\Phi_+(z,k)&=\sum_{j=1}^\infty \phi_{+,j}(k) z^j, \quad
z\in\D, \; k\in\Z, \\ \phi_{+,1}(k)&=-\ol{\alpha_{k+1}},
\no \\ \phi_{+,2}(k)&=-\rho_{k+1}^2 \, \ol{\alpha_{k+2}},
\lb{2.150} \\ \phi_{+,j}(k)&=\alpha_{k+1}\sum_{\ell=1}^j
\phi_{+,j-\ell}(k+1) \phi_{+,\ell}(k) + \phi_{+,j-1}(k+1),
\; j\geq 3. \no
\end{align}
The corresponding Riccati-type equation for the
Caratheodory function $1/\Phi_-(z,k)$ implies the
absolutely convergent expansion
\begin{align}
1/\Phi_-(z,k)&=\sum_{j=0}^\infty [1/\phi_{-,j}(k)] z^j,
\quad z\in\D, \; k\in\Z, \\ 1/\phi_{-,0}(k)&=\alpha_{k},
\no \\ 1/\phi_{-,1}(k)&=\rho_{k}^2 \, \alpha_{k-1},
\lb{2.152} \\
1/\phi_{-,j}(k)&=-\ol{\alpha_{k}}\sum_{\ell=0}^{j-1}
[1/\phi_{-,j-1-\ell}(k-1)][1/\phi_{-,\ell}(k)] +
[1/\phi_{-,j-1}(k-1)], \; j\geq 2. \no
\end{align}

Next, we introduce the following notation for the half-open
arc on the unit circle,
\begin{equation}
\Arc\big(\big(e^{i\theta_1},e^{i\theta_2}\big]\big)
=\big\{e^{i\theta}\in\dD\,|\, \theta_1<\theta\leq
\theta_2\big\}, \quad \theta_1 \in [0,2\pi), \;
\theta_1<\theta_2\leq \theta_1+2\pi.
\end{equation}
In the same manner we also introduce open and closed arcs
on $\dD$,
$\Arc\big(\big(e^{i\theta_1},e^{i\theta_2}\big)\big)$ and
$\Arc\big(\big[e^{i\theta_1},e^{i\theta_2}\big]\big)$,
respectively. Moreover, we identify the unit circle $\dD$
with the arcs of the form
$\Arc\big(\big(e^{i\theta_1},e^{i\theta_1+2\pi}\big]\big)$,
$\te_1\in[0,2\pi)$.

The following result is the unitary operator analog of a
version of Stone's formula relating resolvents of
self-adjoint operators with spectral projections in the
weak sense (cf., e.g., \cite[p.\ 1203]{DS88}).

%%%%%%%%%%%%%%%%%%%%%%%%%%%%%%%%%%%%%%%%%%%%%%%%%%%%%%%%%%%%%%%%%%%%
\begin{lemma} \lb{l2.23}
Let $U$ be a unitary operator in a complex separable
Hilbert space $\cH$ $($with scalar product denoted by
$(\cdot,\cdot)_\cH$, linear in the second factor$)$, $f,g
\in\cH$, and denote by $\{E_U(\zeta)\}_{\zeta\in\dD}$ the
family of self-adjoint right-continuous spectral
projections associated with $U$, that is,
$(f,Ug)_\cH=\int_{\dD} d(f,E_U(\zeta)g)_\cH\, \zeta$.
Moreover, let $\theta_1 \in [0,2\pi)$,
$\theta_1<\theta_2\leq \theta_1+2\pi$, $F\in C(\dD)$, and
denote by $C(U,z)$ the operator
\begin{equation}
C(U,z)=(U+zI_\cH)(U-zI_\cH)^{-1}
=I_\cH+2z(U-zI_\cH)^{-1}, \quad z \in \bbC\backslash\sigma(U)
\lb{2.163A}
\end{equation}
with $I_\cH$ the identity operator in $\cH$. Then,
\begin{align}
&\big(f,F(U)E_{U}\big(\Arc\big(\big(e^{i\theta_1},
e^{i\theta_2}\big]\big)\big)g\big)_{\cH} \no \\ & \quad =
\lim_{\delta\downarrow 0}\lim_{r\uparrow 1}
\int_{\theta_1+\delta}^{\theta_2+\delta} \f{d\theta}{4\pi}
\, F\big(e^{i\theta}\big)
\big[\big(f,C\big(U,re^{i\theta}\big)g)_{\cH}
            - \big(f,C\big(U,r^{-1}e^{i\theta}\big)g\big)_{\cH}\big].
\lb{2.164}
\end{align}
Similar formulas hold for
$\Arc\big(\big(e^{i\theta_1},e^{i\theta_2}\big)\big)$ and
$\Arc\big(\big[e^{i\theta_1},e^{i\theta_2}\big]\big)$.
\end{lemma}
%%%%%%%%%%%%%%%%%%%%%%%%%%%%%%%%%%%%%%%%%%%%%%%%%%%%%%%%%%%%%%%%%%%%
\begin{proof}
First one notices that
\begin{equation}
C\big(U,re^{i\theta}\big)^* =- C\big(U,r^{-1}e^{i\theta}\big), \quad
r\in (0,\infty)\backslash\{1\}, \; \theta\in [0,2\pi]. \lb{2.165}
\end{equation}
Next, introducing the characteristic function $\chi_A$ of a
set $A\subseteq\dD$ and  assuming $F\geq 0$, one obtains
that
\begin{align}
& \big(F(U)^{1/2}E_{U}\big(\Arc\big(\big(e^{i\theta_1},
e^{i\theta_2}\big]\big)\big)f, C(U,z)
F(U)^{1/2}E_{U}\big(\Arc\big(\big(e^{i\theta_1},
e^{i\theta_2}\big]\big)\big)f\big)_{\cH} \no \\ & \quad =
\int_\dD d\big(f,E_U\big(e^{i\theta}\big)f\big)_{\cH} \,
F\big(e^{i\theta}\big)
\chi_{(e^{i\theta_1},e^{i\theta_2}]}\big(e^{i\theta}\big)
\f{e^{i\theta}+z}{e^{i\theta}-z} \no \\ & \quad = \int_\dD
d\Big(F(U)^{1/2}
\chi_{(e^{i\theta_1},e^{i\theta_2}]}(U)f,E_U\big(e^{i\theta}\big)
F(U)^{1/2}\chi_{(e^{i\theta_1},e^{i\theta_2}]}(U)f\Big)_{\cH}
\, \f{e^{i\theta}+z}{e^{i\theta}-z}, \no \\[1mm] &
\hspace*{9.8cm}  z\in\dD  \lb{2.166}
\end{align}
is a Caratheodory function and hence \eqref{2.164} for $g=f$ follows from
\eqref{A.4a}. If $F$ is not nonnegative, one decomposes $F$ as
$F=(F_1-F_2)+i(F_3-F_4)$ with $F_j\geq 0$ and
applies \eqref{2.166} to each $F_j$, $j\in\{1,2,3,4\}$. The general case
$g\neq f$ then follows from the special case $g=f$ by polarization.
\end{proof}
%%%%%%%%%%%%%%%%%%%%%%%%%%%%%%%%%%%%%%%%%%%%%%%%%%%%%%%%%%%%%%%%%%%%

Next, in addition to the definition of $\wti p_\pm$ and $\wti q_\pm$ in
\eqref{2.99}--\eqref{2.102} we introduce $\wti u_+$ by
\begin{align}
&\binom{\wti u_+(z,\cdot,k_0)}{v_+(z,\cdot,k_0)} =
\binom{\wti q_+(z,\cdot,k_0)}{s_+(z,\cdot,k_0)} + m_+(z,k_0)
\binom{\wti p_+(z,\cdot,k_0)}{r_+(z,\cdot,k_0)} \in
\ell^2([k_0,\infty)\cap\Z)^2, \no \\
& \hspace*{8cm} z\in\bbC\backslash(\dD\cup\{0\}) \label{2.167}
\end{align}
and the functions $\ti t_-$ and $w_-$ by
\begin{align}
&\binom{\wti t_-(z,\cdot,k_0)}{w_-(z,\cdot,k_0)} =
\binom{\wti q_-(z,\cdot,k_0)}{s_-(z,\cdot,k_0)} + m_-(z,k_0)
\binom{\wti p_-(z,\cdot,k_0)}{r_-(z,\cdot,k_0)} \in
\ell^2((-\infty,k_0]\cap\Z)^2, \no \\
& \hspace*{8cm} z\in\bbC\backslash(\dD\cup\{0\}). \label{2.167a}
\end{align}
One then computes for the resolvent of $U_{\pm,k_0}$ in
terms of its matrix representation in the standard basis of
$\ell^2([k_0,\pm\infty)\cap\bbZ)$,
\begin{align}
(U_{+,k_0}-zI)^{-1}(k,k')
& =\f{1}{2z} \begin{cases} \wti p_+(z,k,k_0)v_+(z,k',k_0), &
\text{$k<k'$ or $k=k'$ odd,} \\
r_+(z,k',k_0) \wti u_+(z,k,k_0), & \text{$k'<k$ or $k=k'$ even,}
\end{cases} \no \\
&  \quad \; \, z\in\bbC\backslash(\dD\cup\{0\}), \; k_0\in\bbZ, \;
k,k'\in [k_0,\infty)\cap\bbZ, \lb{2.168} \\
(U_{-,k_0}-zI)^{-1}(k,k') & =\f{1}{2z} \begin{cases} \wti
t_-(z,k,k_0)r_-(z,k',k_0), & \text{$k<k'$ or $k=k'$ odd,}
\\ w_-(z,k',k_0) \wti p_-(z,k,k_0), & \text{$k'<k$ or $k=k'$ even,}
\end{cases} \no \\
&  \quad z\in\bbC\backslash(\dD\cup\{0\}), \; k_0\in\bbZ, \;
k,k'\in (-\infty,k_0]\cap\bbZ. \lb{2.168A}
\end{align}
The proof of these formulas repeats the proof of the
analogous result, Lemma \ref{lA.18}, for the full-lattice
CMV operator $U$ and hence we omit it here.

We finish this section with an explicit connection between
the family of spectral projections of $U_{\pm,k_0}$ and the
spectral function $\mu_\pm(\cdot,k_0)$, supplementing
relation \eqref{A.51}.

%%%%%%%%%%%%%%%%%%%%%%%%%%%%%%%%%%%%%%%%%%%%%%%%%%%%%%%%%%%%%%%%%%%%
\begin{lemma} \lb{l2.24}
Let $f,g \in\ell^\infty_0([k_0,\pm\infty)\cap\Z)$, $F\in C(\dD)$, and
$\theta_1\in [0,2\pi)$, $\theta_1<\theta_2\leq \theta_1+2\pi$. Then,
\begin{align}
\begin{split}
& \big(f,F(U_{\pm,k_0})E_{U_{\pm,k_0}}
\big(\Arc\big(\big(e^{i\theta_1},e^{i\theta_2}\big]\big)\big)
g\big)_{\ell^2([k_0,\pm\infty)\cap\Z)}  \\
& \quad =\big(\hatt f_\pm(\cdot,k_0),M_F
M_{\chi_{\Arc((e^{i\theta_1},e^{i\theta_2}])}}
\hatt g_\pm(\cdot,k_0)\big)_{L^2(\dD;d\mu_\pm (\cdot,k_0))},   \lb{2.169}
\end{split}
\end{align}
where we introduced the notation
\begin{equation}
\hatt h_\pm(\zeta,k_0)=\sum_{k=k_0}^{\pm\infty} r_\pm(\zeta,k,k_0) h(k),
\quad
\zeta\in\dD, \; h\in \ell^\infty_0([k_0,\pm\infty)\cap\Z), \lb{2.170}
\end{equation}
and $M_G$ denotes the maximally defined operator of multiplication by the
$d\mu_{\pm}(\cdot,k_0)$-measurable function $G$ in the
Hilbert space $\Lt{\pm}$,
\begin{align}
\begin{split}
& (M_G\hatt h)(\zeta)=G(\zeta)\hatt h(\zeta) \, \text{ for
a.e.\ $\zeta\in\dD$},  \\ & \hatt h \in \dom(M_G)=\{\hatt
k\in \Lt{\pm} \,|\, G\hatt k \in \Lt{\pm}\}. \lb{2.171}
\end{split}
\end{align}
\end{lemma}
%%%%%%%%%%%%%%%%%%%%%%%%%%%%%%%%%%%%%%%%%%%%%%%%%%%%%%%%%%%%%%%%%%%%
\begin{proof}
It suffices to consider $U_{+,k_0}$ only. Inserting
\eqref{2.168} into \eqref{2.164} and observing
\eqref{2.167} leads to
\begin{align}
&\big(f,F(U_{+,k_0})E_{U_{+,k_0}}\big(\Arc\big(\big(e^{i\theta_1},
e^{i\theta_2}\big]\big)\big)g\big)_{\ell^2([k_0,\infty)\cap\Z)} \no \\
& \quad = \lim_{\delta\downarrow 0}\lim_{r\uparrow 1}
\int_{\theta_1+\delta}^{\theta_2+\delta} \f{d\theta}{4\pi} \,
F\big(e^{i\theta}\big)
\bigg[\sum_{k=k_0}^{\infty}\sum_{k'=k_0}^{\infty}
\ol{f(k)}g(k')[C\big(U_{+,k_0},re^{i\theta}\big)(k,k')  \no \\
& \hspace*{7.7cm}   - C\big(U_{+,k_0},r^{-1}e^{i\theta}\big)(k,k')\bigg]
\no \\
& \quad = \sum_{k=k_0}^{\infty} \ol{f(k)}\bigg\{ \sum_{\substack{k_0\leq
          k'<k \\ k'=k\,\text{even}}} g(k') \lim_{\delta\downarrow 0}
\lim_{r\uparrow 1}\f{1}{4\pi}\int_{\theta_1+\delta}^{\theta_2+\delta}
d\theta \, F\big(e^{i\theta}\big)\wti p_+\big(e^{i\theta},k,k_0\big) \no
\\ & \hspace*{4.6cm} \times r_+\big(e^{i\theta},k',k_0\big)
\big[m_+\big(re^{i\theta},k_0\big)
-m_+\big(r^{-1}e^{i\theta}, k_0\big)\big] \no  \\
& \hspace*{2.6cm} + \sum_{\substack{k_0\leq
          k<k' \\ k'=k\,\text{odd}}} g(k') \lim_{\delta\downarrow 0}
\lim_{r\uparrow 1}\f{1}{4\pi}\int_{\theta_1+\delta}^{\theta_2+\delta}
d\theta \, F\big(e^{i\theta}\big)\wti p_+\big(e^{i\theta},k,k_0\big) \no
\\
& \hspace*{4.6cm} \times r_+\big(e^{i\theta},k',k_0\big)
\big[m_+\big(re^{i\theta},k_0\big)
-m_+\big(r^{-1}e^{i\theta}, k_0\big)\big]\bigg\}.
\end{align}
Here we freely interchanged the $\theta$-integral with the sums over $k$
and $k'$ (the latter are finite) and also replaced
$\wti p_+\big(r^{\pm 1}e^{i\theta},k,k_0\big)$ and
$r_+\big(r^{\pm 1}e^{i\theta},k,k_0\big)$ by
$\wti p_+\big(e^{i\theta},k,k_0\big)$ and
$r_+\big(e^{i\theta},k,k_0\big)$. The latter is permissible since by
\eqref{A.10C},
\begin{equation}
\big|\big(1-r^{\pm 1}\big)\Re\big(m_+\big(r^{\pm
1}e^{i\theta}\big)\big)\big| \underset{r\to 1}{=} \Oh(1), \quad
\big|\big(1-r^{\pm 1}\big)\Im\big(m_+\big(r^{\pm
1}e^{i\theta}\big)\big)\big| \underset{r\to 1}{=} \oh(1).
\end{equation}
Finally, since $\wti
p_+(\zeta,k,k_0)=\ol{r_+(\zeta,k,k_0)}$, $\zeta\in\dD$ by
\eqref{2.59} and $m_+(re^{i\te},k_0) =
-\ol{m_+(\frac{1}{r}e^{i\te},k_0)}$ by \eqref{A.7}, one
infers
\begin{align}
&\big(f,F(U_{+,k_0})E_{U_{+,k_0}}\big(\Arc\big(\big(e^{i\theta_1},
e^{i\theta_2}\big]\big)\big)g\big)_{\ell^2([k_0,\infty)\cap\Z)}
\no \\ & \quad =
\sum_{k=k_0}^{\infty}\sum_{k'=k_0}^{\infty} \ol{f(k)}g(k')
\lim_{\delta\downarrow 0}\lim_{r\uparrow 1}
\int_{\theta_1+\delta}^{\theta_2+\delta} \f{d\theta}{2\pi}
\, F\big(e^{i\theta}\big) \wti
p_+\big(e^{i\theta},k,k_0\big)
r_+\big(e^{i\theta},k',k_0\big) \no \\ & \hspace*{6.45cm}
\times \Re\big(m_+\big(re^{i\theta},k_0\big)\big) \no \\ &
\quad = \sum_{k=k_0}^{\infty}\sum_{k'=k_0}^{\infty}
\ol{f(k)}g(k') \int_{(\theta_1,\theta_2]}
d\mu_+\big(e^{i\theta},k_0\big) \, F\big(e^{i\theta}\big)
\ol{r_+\big(e^{i\theta},k,k_0\big)}
r_+\big(e^{i\theta},k',k_0\big) \no \\ & \quad =
\int_{(\theta_1,\theta_2]} d\mu_+\big(e^{i\theta},k_0\big)
\, F\big(e^{i\theta}\big) \ol{\hatt
f_+\big(e^{i\theta},k_0\big)} \hatt
g_+\big(e^{i\theta},k_0\big) \no \\ & \quad =\big(\hatt
f_+(\cdot,k_0),M_F
M_{\chi_{\Arc((e^{i\theta_1},e^{i\theta_2}])}} \hatt
g_+(\cdot,k_0)\big)_{L^2(\dD;d\mu_+ (\cdot,k_0))},
\end{align}
interchanging the (finite) sums over $k$ and $k'$ and the
$d\mu_+(\cdot,k_0)$-integral once more.
\end{proof}
%%%%%%%%%%%%%%%%%%%%%%%%%%%%%%%%%%%%%%%%%%%%%%%%%%%%%%%%%%%%%%%%%%%%

Finally, this section would not be complete if we wouldn't
briefly mention the analogs of Weyl disks for finite
interval problems and their behavior in the limit where the
finite interval tends to a half-lattice. Before starting
the analysis, we note the following geometric fact: Let
$p,q,r,s\in\bbC$, $|p|\neq |r|$. Then, the set of points
$m(\theta)\in\bbC$ given by
\begin{equation}
m(\theta)=-\f{q+se^{i\theta}}{p+re^{i\theta}}, \quad \theta\in [0,2\pi),
\lb{2.175}
\end{equation}
describes a circle in $\bbC$ with radius $R>0$ and center $C\in\bbC$
given by
\begin{equation}
R=\f{|qr-ps|}{\big||p|^2-|r|^2\big|}, \quad
C=-\f{s}{r}-\f{\ol p}{r}\f{qr-ps}{|p|^2-|r|^2}.  \lb{2.176}
\end{equation}
To introduce the analog of $\UU^{(s)}_{+,k_0}$ and
$\big(\UU^{(s)}_{+,k_0}\big)^\top$ on a finite interval
$[k_0,k_1]\cap\bbZ$, we choose $\alpha_{k_0}=e^{is_0}$,
$\alpha_{k_1+1}=e^{is_1}$, $s_0, s_1 \in [0,2\pi)$. Then
the operator $\UU_{+,k_0}^{(s_0)}$ splits into a direct sum
of two operators $\UU^{(s_0,s_1)}_{[k_0,k_1]}$ and
$\UU_{+,{k_1+1}}^{(s_1)}$
\begin{equation}
\UU_{+,k_0}^{(s_0)} = \UU^{(s_0,s_1)}_{[k_0,k_1]} \oplus
\UU_{+,{k_1+1}}^{(s_1)}  \lb{2.177}
\end{equation}
acting on $\ell^2([k_0,k_1]\cap\bbZ)$ and
$\ell^2([k_1+1,\infty)\cap\bbZ)$, respectively. Then,
repeating the proof of Lemma \ref{lA.2} one obtains the
following result for the CMV operator
$\UU^{(s_0,s_1)}_{[k_0,k_1]}$:
\begin{align}
& \UU^{(s_0,s_1)}_{[k_0,k_1]}\begin{pmatrix} u(z,\cdot) \\
v(z,\cdot)
\end{pmatrix} =z \begin{pmatrix} u(z,\cdot) \\ v(z,\cdot)
\end{pmatrix}, \quad z\in\bbC\backslash\{0\}  \lb{2.178}
\intertext{is satisfied by $\left(\begin{smallmatrix}
u(z,k)\\v(z,k)\end{smallmatrix}\right)_{k\in[k_0,k_1]\cap\bbZ}$
such that} & \begin{pmatrix} u(z,k) \\ v(z,k) \end{pmatrix} =
T(z,k)
\begin{pmatrix} u(z,k-1) \\ v(z,k-1) \end{pmatrix}, \quad
k\in[k_0+1,k_1]\cap\bbZ, \lb{2.179} \\ &
u(z,k_0)=\begin{cases} ze^{is_0} v(z,k_0), & \text{$k_0$
odd,} \\ e^{-is_0} v(z,k_0), & \text{$k_0$ even,}
\end{cases}  \lb{2.180} \\ & u(z,k_1)=\begin{cases} -
e^{is_1} v(z,k_1), & \text{$k_1$ odd,} \\ -ze^{-is_1}
v(z,k_1), & \text{$k_1$ even.} \end{cases} \lb{2.181}
\end{align}
To simplify matters we now put $s_0=0$ in the following. Moreover, we
first treat the case $k_0$ even and $k_1$ odd. Then
$\left(\begin{smallmatrix} p_+(z,k,k_0) \\
r_+(z,k,k_0)\end{smallmatrix}\right)$ satisfies \eqref{2.179} and
\eqref{2.180} and hence there exists a coefficient $m_{+,s_1}(z,k_1,k_0)$
such that
\begin{equation}
\begin{pmatrix} q_+(z,k,k_0) \\ s_+(z,k,k_0) \end{pmatrix} +
m_{+,s_1}(z,k_0,k_1) \begin{pmatrix} p_+(z,k,k_0) \\ r_+(z,k,k_0)
\end{pmatrix}  \lb{2.182}
\end{equation}
satisfies \eqref{2.181}. One computes
\begin{equation}
m_{+,s_1}(z,k_1,k_0)=-\f{q_+(z,k_1,k_0)+s_+(z,k_1,k_0)
e^{is_1}}{p_+(z,k_1,k_0)+r_+(z,k_1,k_0)e^{is_1}}. \lb{2.183}
\end{equation}
By \eqref{2.175}, this describes a (Weyl--Titchmarsh) circle as $s_1$
varies in $[0,2\pi)$ of radius
\begin{align}
R(z,k_1)&=\f{|q_+(z,k_1,k_0)r_+(z,k_1,k_0)-p_+(z,k_1,k_0)
s_+(z,k_1,k_0))|}{\big||p_+(z,k_1,k_0)|^2-|r_+(z,k_1,k_0)|^2\big|} \no\\
&= \f{2}{\big||p_+(z,k_1,k_0)|^2-|r_+(z,k_1,k_0)|^2\big|} \lb{2.184}
\end{align}
since
\begin{equation}
W\bigg(\begin{pmatrix} p_+(z,k_1,k_0) \\ r_+(z,k_1,k_0)
\end{pmatrix}
\begin{pmatrix} q_+(z,k_1,k_0) \\ s_+(z,k_1,k_0) \end{pmatrix} \bigg)
= 2   \lb{2.185}
\end{equation}
if $k_0$ is even and $k_1$ is odd (cf.\ also \eqref{A.82}).

Thus far our computations are subject to $|p_+(z,k_1,k_0)|\neq
|r_+(z,k_1.k_0)|$. To clarify this point we now state the following
result.

%%%%%%%%%%%%%%%%%%%%%%%%%%%%%%%%%%%%%%%%%%%%%%%%%%%%%%%%%
\begin{lemma}  \lb{l2.25}
Let $z\in\bbC\backslash(\dD\cup\{0\})$ and $k_0, k_1\in\bbZ$, $k_1>k_0$.
Then,
\begin{equation}
\big(1-|z|^{-2}\big) \sum_{k=k_0}^{k_1} |p_+(z,k,k_0)|^2=\begin{cases}
|p_+(z,k_1,k_0)|^2 - |r_+(z,k_1.k_0)|^2, & \text{$k_1$ odd,} \\
|r_+(z,k_1.k_0)|^2 - |z|^{-2}|p_+(z,k_1,k_0)|^2, & \text{$k_1$ even.}
\end{cases}    \lb{2.186}
\end{equation}
\end{lemma}
%%%%%%%%%%%%%%%%%%%%%%%%%%%%%%%%%%%%%%%%%%%%%%%%%%%%%%%%%
\begin{proof}
It suffices to prove the case $k_1$ odd. The computation
\begin{align}
& \ol z \sum_{k=k_0}^{k_1} |p_+(z,k,k_0)|^2 =
\sum_{k=k_0}^{k_1}
\ol{(U_{+,k_0}p_+(z,\cdot,k_0))(k)}p_+(z,k,k_0)  \no \\ & =
\sum_{k=k_0}^{k_1-1} \ol{(V_{+,k_0}W_{+,k_0}
p_+(z,\cdot,k_0))(k)}p_+(z,k,k_0) + \ol z
|p_+(z,k_1,k_0)|^2 \no \\ & = \sum_{k=k_0}^{k_1-1}
\ol{(W_{+,k_0}
p_+(z,\cdot,k_0))(k)}(V_{+,k_0}^*p_+(z,\cdot,k_0))(k) + \ol
z |p_+(z,k_1,k_0)|^2  \no \\ & = \sum_{k=k_0}^{k_1}
\ol{p_+(z,k,k_0)}
(W_{+,k_0}^*V_{+,k_0}^*p_+(z,\cdot,k_0))(k) \no \\ & \quad
-
\ol{(W_{+,k_0}p_+(z,\cdot,k_0))(k_1)}
(V_{+,k_0}^*p_+(z,\cdot,k_0))(k_1) + \ol z
|p_+(z,k_1,k_0)|^2  \no \\ & = \sum_{k=k_0}^{k_1}
\ol{p_+(z,k,k_0)} (U_{+,k_0}^*p_+(z,\cdot,k_0))(k) - \ol z
|r_+(z,k_1,k_0))|^2 + \ol z |p_+(z,k_1,k_0)|^2  \no \\ & =
z^{-1} \sum_{k=k_0}^{k_1} |p_+(z,k,k_0)|^2 - \ol z
|r_+(z,k_1,k_0))|^2 + \ol z |p_+(z,k_1,k_0)|^2  \lb{2.187}
\end{align}
proves \eqref{2.186} for $k_1$ odd.
\end{proof}
%%%%%%%%%%%%%%%%%%%%%%%%%%%%%%%%%%%%%%%%%%%%%%%%%%%%%%%%%

A systematic investigation of all even/odd possibilities for $k_0$ and
$k_1$ then yields the following result.

%%%%%%%%%%%%%%%%%%%%%%%%%%%%%%%%%%%%%%%%%%%%%%%%%%%%%%%%%
\begin{theorem}  \lb{t2.26}
Let $z\in\bbC\backslash(\dD\cup\{0\})$ and $k_0, k_1\in\bbZ$, $k_1>k_0$.
Then,
\begin{equation}
m_{+,s_1}(z,k_1,k_0)=\begin{cases}
-\f{q_+(z,k_1,k_0)+s_+(z,k_1,k_0)e^{is_1}}{p_+(z,k_1,k_0)
+r_+(z,k_1,k_0)e^{is_1}}, & \text{$k_1$ odd,} \\
-\f{z^{-1}q_+(z,k_1,k_0)+s_+(z,k_1,k_0)e^{-is_1}}{z^{-1}p_+(z,k_1,k_0)
+r_+(z,k_1,k_0)e^{-is_1}}, & \text{$k_1$ even}
\end{cases}  \lb{2.188}
\end{equation}
lies on a circle of radius
\begin{equation}
R(z,k_1,k_0)=\bigg[\big|1-|z|^{-2} \big| \sum_{k=k_0}^{k_1}
|p_+(z,k_1,k_0)|^2\bigg]^{-1}
\begin{cases} 2|z|, & \text{$k_0$ odd, $k_1$ odd,} \\
2, & \text{$k_0$ even, $k_1$ odd,} \\
2, & \text{$k_0$ odd, $k_1$ even,} \\
2|z|^{-1}, & \text{$k_0$ even, $k_1$ even}
\end{cases}  \lb{2.189}
\end{equation}
with center
\begin{equation}
C(z,k_1,k_0)=\begin{cases}
-\f{s_+(z,k_1,k_0)}{r_+(z,k_1,k_0)} -
\f{\ol{p_+(z,k_1,k_0)}}{r_+(z,k_1,k_0)} \\[1mm]
\hspace*{.5cm} \times \f{2z}{|p_+(z,k_1,k_0)|^2
-|r_+(z,k_1,k_0)|^2}, & \text{$k_0$ odd, $k_1$ odd,}
\\[2mm]
-\f{s_+(z,k_1,k_0)}{r_+(z,k_1,k_0)} -
\f{\ol{p_+(z,k_1,k_0)}}{r_+(z,k_1,k_0)} \\[1mm]
\hspace*{.5cm} \times \f{2}{|p_+(z,k_1,k_0)|^2
-|r_+(z,k_1,k_0)|^2}, & \text{$k_0$ even, $k_1$ odd,}
\\[2mm]
-\f{s_+(z,k_1,k_0)}{r_+(z,k_1,k_0)} -
\f{\ol{z^{-1}p_+(z,k_1,k_0)}}{r_+(z,k_1,k_0)} \\[1mm]
\hspace*{.5cm} \times \f{-2}{|z|^{-2}|p_+(z,k_1,k_0)|^2
-|r_+(z,k_1,k_0)|^2}, & \text{$k_0$ odd, $k_1$ even,}
\\[2mm]
-\f{s_+(z,k_1,k_0)}{r_+(z,k_1,k_0)} -
\f{\ol{z^{-1}p_+(z,k_1,k_0)}}{r_+(z,k_1,k_0)} \\[1mm]
\hspace*{.5cm} \times
\f{-2z^{-1}}{|z|^{-2}|p_+(z,k_1,k_0)|^2
-|r_+(z,k_1,k_0)|^2}, & \text{$k_0$ even, $k_1$ even.}
\end{cases}  \lb{2.190}
\end{equation}
In particular, the limit point case holds at $+\infty$ since
\begin{equation}
\lim_{k_1\uparrow\infty} R(z,k_1,k_0) = 0.  \lb{2.191}
\end{equation}
\end{theorem}
%%%%%%%%%%%%%%%%%%%%%%%%%%%%%%%%%%%%%%%%%%%%%%%%%%%%%%%%%
\begin{proof}
The case $k_0$ even, $k_1$ odd has been discussed explicitly in
\eqref{2.182}--\eqref{2.186}. The remaining cases follow similarly using
Lemma \ref{l2.25} for $k_1$ even and the Wronski relations \eqref{A.82}.
Relation \eqref{2.191} follows since $p_+(z,\cdot,k_0)\notin
\ell^2([k_0,\infty)\cap\bbZ)$, $z\in \bbC\backslash(\dD\cup\{0\})$. The
latter follows from $(U_{+,k_0}p(z,\cdot,k_0))(k)  =z p_+(z,k,k_0)$,
$z\in\bbC\backslash\{0\}$, in the weak sense (cf.\ Remark \ref{rA.3}) and
the fact that $U_{+,k_0}$ is unitary.
\end{proof}
%%%%%%%%%%%%%%%%%%%%%%%%%%%%%%%%%%%%%%%%%%%%%%%%%%%%%%%%%

%%%%%%%%%%%%%%%%%%%%%%%%%%%%%%%%%%%%%%%%%%%%%%%%%%%%%%%%%
\section{Weyl--Titchmarsh Theory for CMV Operators on $\bbZ$}\lb{s3}
%%%%%%%%%%%%%%%%%%%%%%%%%%%%%%%%%%%%%%%%%%%%%%%%%%%%%%%%%

In this section we describe the Weyl--Titchmarsh theory for the
CMV operator $U$ on $\bbZ$. We note that in a context different
from orthogonal polynomials on the unit circle, Bourget, Howland, and
Joye \cite{BHJ03} introduced a set of doubly infinite family of matrices
with three sets of parameters which for special choices of the parameters
reduces to two-sided CMV matrices on $\bbZ$.

We denote by
\begin{align}
& W\left(\binom{u_1(z,k,k_0)}{v_1(z,k,k_0)},
\binom{u_2(z,k,k_0)}{v_2(z,k,k_0)}\right) =
\det\left(\begin{pmatrix}
u_1(z,k,k_0) & u_2(z,k,k_0) \\ v_1(z,k,k_0) & v_2(z,k,k_0)
\end{pmatrix}\right), \lb{A.81} \\
& \hspace*{10cm} \quad k\in\bbZ, \no
\end{align}
the Wronskian of two solutions
$\Big(\begin{smallmatrix}u_1(z,\cdot,k_0)\\
v_1(z,\cdot,k_0)\end{smallmatrix}\Big)$ and
$\Big(\begin{smallmatrix}u_2(z,\cdot,k_0)\\
v_2(z,\cdot,k_0)\end{smallmatrix}\Big)$ of
\eqref{A.28}
for $z\in\bbC\backslash\{0\}$. Then, since
\begin{equation}
\det(T(z,k))=-1, \quad k\in\Z,  \lb{A.81A}
\end{equation}
it follows from Definition \ref{dA.2} that
\begin{align}
W&\left(\binom{p_+(z,k,k_0)}{r_+(z,k,k_0)},
\binom{q_+(z,k,k_0)}{s_+(z,k,k_0)}\right) = (-1)^k
\begin{cases}
2z, & k_0 \text{ odd}, \\ 2,   & k_0 \text{ even},
\end{cases} \lb{A.82}
\\
W&\left(\binom{p_-(z,k,k_0)}{r_-(z,k,k_0)},
\binom{q_-(z,k,k_0)}{s_-(z,k,k_0)}\right) = (-1)^{k+1}
\begin{cases}
2, & k_0 \text{ odd}, \\ 2z,   & k_0 \text{ even},
\end{cases} \lb{A.83} \\
& \hspace*{6.51cm} z\in\bbC\backslash\{0\}, \; k\in\bbZ.
\no
\end{align}
Next, in order to compute the resolvent of $U$, we
introduce in addition to $\wti p_\pm$ and $\wti q_\pm$ in
\eqref{2.99}--\eqref{2.102} the functions $\wti u_\pm$ by
\begin{align}
&\binom{\wti u_\pm(z,\cdot,k_0)}{v_\pm(z,\cdot,k_0)} =
\binom{\wti q_+ (z,\cdot,k_0)}{s_+ (z,\cdot,k_0)} +
M_\pm(z,k_0) \binom{\wti
p_+(z,\cdot,k_0)}{r_+(z,\cdot,k_0)} \in
\ell^2([k_0,\pm\infty)\cap\Z)^2, \no \\ & \hspace*{8.5cm}
z\in\bbC\backslash(\dD\cup\{0\}). \label{3.5A}
\end{align}

%%%%%%%%%%%%%%%%%%%%%%%%%%%%%%%%%%%%%%%%%%%%%%%%%%%%%%%%%%%%%%%
\begin{lemma} \label{lA.18}
Let $z\in\bbC\backslash(\dD\cup\{0\})$ and fix $k_0,
k_1\in\bbZ$. Then the resolvent $(U-zI)^{-1}$ of the
unitary CMV operator $U$  on $\ell^2(\bbZ)$ is given in
terms of its matrix representation in the standard basis of
$\ltz$ by
\begin{align}
& (U-zI)^{-1}(k,k') = \frac{(-1)^{k_1+1}}{z
W\left(\begin{pmatrix}u_+(z,k_1,k_0)\\ v_+(z,k_1,k_0)\end{pmatrix},
\begin{pmatrix}u_-(z,k_1,k_0)\\ v_-(z,k_1,k_0)\end{pmatrix}
\right)} \no \\
& \quad \times \begin{cases}
u_-(z,k,k_0)v_+(z,k',k_0), & k < k' \text{ or } k = k'
\text{ odd},
\\
v_-(z,k',k_0) u_+(z,k,k_0), & k' < k \text{ or } k = k'
\text{ even},
\end{cases} \quad k,k' \in\Z, \label{A.84} \\
& \hspace*{2.65cm} = \frac{-1}{2z[M_+(z,k_0)-M_-(z,k_0)]} \no \\
& \quad \times \begin{cases}
\wti u_-(z,k,k_0)v_+(z,k',k_0), & k < k' \text{ or } k = k'
\text{ odd}, \\
v_-(z,k',k_0) \wti u_+(z,k,k_0), & k' < k \text{ or } k = k'
\text{ even},
\end{cases} \quad k,k' \in\Z, \label{A.84a}
\end{align}
where
\begin{align}
& W\left(\binom{u_+(z,k_1,k_0)}{v_+(z,k_1,k_0)},
\binom{u_-(z,k_1,k_0)}{v_-(z,k_1,k_0)}\right) =
\det\left(\begin{pmatrix}
u_+(z,k_1,k_0) & u_-(z,k_1,k_0) \\ v_+(z,k_1,k_0) &
v_-(z,k_1,k_0)
\end{pmatrix}\right) \no \\
& \quad = (-1)^{k_1}[M_+(z,k_0)-M_-(z,k_0)]
\begin{cases}
2z, & k_0  \text{ odd}, \\ 2,   & k_0 \text{ even},
\end{cases} \lb{A.85}
\end{align}
and
\begin{equation}
W\left(\binom{\wti u_+(z,k_1,k_0)}{v_+(z,k_1,k_0)},
\binom{\wti u_-(z,k_1,k_0)}{v_-(z,k_1,k_0)}\right) = 2(-1)^{k_1}
[M_+(z,k_0)-M_-(z,k_0)].  \lb{A.85a}
\end{equation}
Moreover, since $0\in\bbC\backslash\sigma(U)$, \eqref{A.84}
and \eqref{A.84a} analytically extend to $z=0$.
\end{lemma}
%%%%%%%%%%%%%%%%%%%%%%%%%%%%%%%%%%%%%%%%%%%%%%%%%%%%%%%%%%%%%%%
\begin{proof}
Denote
\begin{align}
w(z,k,k',k_0) &=
\begin{cases}
u_-(z,k,k_0)v_+(z,k',k_0), & k < k' \text{ or } k = k' \text{ odd},
\\
u_+(z,k,k_0)v_-(z,k',k_0), & k' < k \text{ or } k = k' \text{
even},
\end{cases}
\\ \no
& \hspace*{6.26cm} k,k',k_0\in\Z.
\end{align}
We will prove that
\begin{align}
&(U-zI)w(z,\cdot,k',k_0) = (-1)^{k'+1} z
\det\left(\begin{pmatrix}
u_+(z,k',k_0) & u_-(z,k',k_0) \\ v_+(z,k',k_0) &
v_-(z,k',k_0)
\end{pmatrix}\right) \de_{k'}, \label{A.86} \\
& \hspace*{9.6cm} \quad k',k_0\in\Z,  \no
\end{align}
and hence, using \eqref{A.81A}, one obtains
\begin{align}
&(U-zI)w(z,\cdot,k',k_0) = (-1)^{k_1+1} z
\det\left(\begin{pmatrix}
u_+(z,k_1,k_0) & u_-(z,k_1,k_0) \\ v_+(z,k_1,k_0) &
v_-(z,k_1,k_0)
\end{pmatrix}\right) \de_{k'},  \no \\
& \hspace*{9cm} k',k_0,k_1\in\Z.
\end{align}

First, let $k_0\in\Z$ and assume $k'$ to be odd. Then,
\begin{align}
\big((U-zI)w(z,\cdot,k',k_0)\big)(\ell) =
\big((VW-zI)w&(z,\cdot,k',k_0)\big)(\ell) = 0, \quad \ell
\in \Z\backslash \{k',k'+1\}
\end{align}
and
\begin{align}
&\left(\begin{matrix}((U-zI)w(z,\cdot,k',k_0))(k')\\
((U-zI)w(z,\cdot,k',k_0))(k'+1)\end{matrix}\right)
=
\left(\begin{matrix}((VW-zI)w(z,\cdot,k',k_0))(k')\\
((VW-zI)w(z,\cdot,k',k_0))(k'+1)\end{matrix}\right) \no \\
&\quad = \te_{k'+1} z \left(\begin{matrix}
(v_+(z,k',k_0)v_-(z,\cdot,k_0))(k')\\
\big(v_-(z,k',k_0)v_+(z,\cdot,k_0))(k'+1)\end{matrix}\right)
- z \binom{w(z,k',k',k_0)}{w(z,k'+1,k',k_0)} \no \\ & \quad =
z v_-(z,k',k_0) \binom{u_+(z,k',k_0)}{u_+(z,k'+1,k_0)} - z
\binom{v_+(z,k',k_0)u_-(z,k',k_0)}{v_-(z,k',k_0)u_+(z,k'+1,k_0)}
\no \\ & \quad = z
\left(\begin{matrix}\det\left(\begin{pmatrix} u_+(z,k',k_0) &
u_-(z,k',k_0) \\ v_+(z,k',k_0) & v_-(z,k',k_0)
\end{pmatrix}\right) \\ 0\end{matrix}\right).
\end{align}

Next, assume $k'$ to be even. Then,
\begin{align}
\big((U-zI)w(z,\cdot,k',k_0)\big)(\ell) =
\big((VW-zI)w&(z,\cdot,k',k_0)\big)(\ell) = 0, \quad \ell
\in \Z\backslash \{k'-1,k'\}
\end{align}
and
\begin{align}
&\binom {((U-zI)w(z,\cdot,k',k_0))(k'-1)}
{((U-zI)w(z,\cdot,k',k_0))(k')}
=
\binom {((VW-zI)w(z,\cdot,k',k_0))(k'-1)}
{((VW-zI)w(z,\cdot,k',k_0))(k')} \no \\ & \quad = \te_{k'}
z \binom {(v_+(z,k',k_0)v_-(z,\cdot,k_0))(k'-1)}
{(v_-(z,k',k_0)v_+(z,\cdot,k_0))(k')} - z
\binom{w(z,k'-1,k',k_0)}{w(z,k',k',k_0)} \no \\ & \quad = z
v_+(z,k',k_0) \binom{u_-(z,k'-1,k_0)}{u_-(z,k',k_0)} - z
\binom{v_+(z,k',k_0)u_-(z,k'-1,k_0)}{v_-(z,k',k_0)u_+(z,k',k_0)}
\no \\ & \quad = z
\left(\begin{matrix}0\\-\det\left(\begin{pmatrix}
u_+(z,k',k_0) & u_-(z,k',k_0) \\ v_+(z,k',k_0) &
v_-(z,k',k_0)
\end{pmatrix}\right)\end{matrix}\right).
\end{align}

Thus, one obtains \eqref{A.86}.
\end{proof}
%%%%%%%%%%%%%%%%%%%%%%%%%%%%%%%%%%%%%%%%%%%%%%%%%%%%%%%%%%%%%%%

Next, we denote by $d\Omega(\cdot,k)$, $k\in\Z$, the $2
\times 2$ matrix-valued measure,
\begin{align}
d\Omega(\ze,k) &= d
\begin{pmatrix}
\Omega_{0,0}(\ze,k) & \Omega_{0,1}(\ze,k)
\\
\Omega_{1,0}(\ze,k) & \Omega_{1,1}(\ze,k)
\end{pmatrix}  \no \\
&= d
\begin{pmatrix}
(\de_{k-1},E_{U}(\ze)\de_{k-1})_{\ell^2(\Z)}
&(\de_{k-1},E_{U}(\ze)\de_{k})_{\ell^2(\Z)} \\
(\de_{k},E_{U}(\ze)\de_{k-1})_{\ell^2(\Z)} &
(\de_{k},E_{U}(\ze)\de_{k})_{\ell^2(\Z)}
\end{pmatrix}, \quad \ze \in\dD, \label{A.87}
\end{align}
where $dE_{U}(\cdot)$ denotes the operator-valued spectral
measure of the unitary CMV operator $U$ on $\ell^2(\bbZ)$,
\begin{equation}
U=\oint_{\dD} dE_{U}(\zeta)\,\zeta.
\end{equation}
We note that by \eqref{A.87} $d\Omega_{0,0}(\cdot,k)$ and
$d\Omega_{1,1}(\cdot,k)$ are nonnegative
measures on $\dD$ and $d\Omega_{0,1}(\cdot,k)$ and
$d\Omega_{1,0}(\cdot,k)$ are complex-valued measures on
$\dD$.

We also introduce the $2 \times 2$ matrix-valued
function $\cM(\cdot,k)$, $k\in\Z$, by
\begin{align}
& \cM(z,k) =
\begin{pmatrix}
M_{0,0}(z,k) & M_{0,1}(z,k) \\
M_{1,0}(z,k) & M_{1,1}(z,k)
\end{pmatrix} \no \\
& \quad = \begin{pmatrix}
(\de_{k-1},(U+zI)(U-zI)^{-1}\de_{k-1})_{\ell^2(\Z)}
&(\de_{k-1},(U+zI)(U-zI)^{-1}\de_{k})_{\ell^2(\Z)}
\\
(\de_{k},(U+zI)(U-zI)^{-1}\de_{k-1})_{\ell^2(\Z)} &
(\de_{k},(U+zI)(U-zI)^{-1}\de_{k})_{\ell^2(\Z)}
\end{pmatrix} \no
\\ & \quad =
\oint_\dD d\Omega(\ze,k)\, \frac{\ze+z}{\ze-z}, \quad
z\in\bbC\backslash\dD. \lb{A.88}
\end{align}
We note that,
\begin{align}
M_{0,0}(\cdot,k+1) = M_{1,1}(\cdot,k), \quad k\in\bbZ
\end{align}
and
\begin{align}
M_{1,1}(z,k) &= (\de_{k}, (U+zI)(U-zI)^{-1}\de_{k})_{\ell^2(\Z)}
\lb{A.89}
\\ & = \oint_\dD d\Omega_{1,1}(\zeta,k) \, \frac{\zeta+z}{\zeta-z},
\quad z\in\bbC\backslash\dD,\; k\in\Z, \lb{A.90}
\end{align}
where
\begin{equation}
d\Omega_{1,1}(\zeta,k)=d(\de_{k},E_U(\zeta)\de_{k})_{\ell^2(\Z)}, \quad
\zeta\in\dD. \lb{A.91}
\end{equation}
Thus, $M_{0,0}|_{\D}$ and $M_{1,1}|_{\D}$ are Caratheodory
functions. Moreover, by \eqref{A.89} one infers that
\begin{equation}
M_{1,1}(0,k)=1, \quad k\in\Z. \lb{A.93}
\end{equation}

%%%%%%%%%%%%%%%%%%%%%%%%%%%%%%%%%%%%%%%%%%%%%%%%%%%%%%%%%%%%%%%
\begin{lemma} \label{lA.20}
Let $z\in\bbC\backslash\dD$. Then the functions $M_{1,1}(\cdot,k)$ and
$M_\pm(\cdot,k)$,
$k\in\bbZ$, satisfy the following relations
\begin{align}
M_{0,0}(z,k)&=1+\f{[\ol{a_k}-\ol{b_k}M_+(z,k)][a_k
+b_kM_-(z,k)]}{\rho_k^2[M_+(z,k)-M_-(z,k)]}, \lb{A.97a} \\
M_{1,1}(z,k) &= \frac{1-M_+(z,k)M_-(z,k)}{M_+(z,k)-M_-(z,k)},
\label{A.92} \\
M_{0,1}(z,k)&=\f{-1}{\rho_k[M_+(z,k)-M_-(z,k)]}
\begin{cases}
{[1-M_+(z,k)][\ol{a_k}-\ol{b_k}M_-(z,k)]}, & \text{$k$ odd}, \\
{[1+M_+(z,k)][a_k+b_kM_-(z,k)]}, & \text{$k$ even},
\end{cases} \lb{A.97b} \\
M_{1,0}(z,k)&=\f{-1}{\rho_k[M_+(z,k)-M_-(z,k)]}
\begin{cases}
{[1+M_+(z,k)][a_k+b_kM_-(z,k)]}, & \text{$k$ odd}, \\
{[1-M_+(z,k)][\ol{a_k}-\ol{b_k}M_-(z,k)]}, & \text{$k$ even}.
\end{cases} \lb{A.97c}
\end{align}
\end{lemma}
%%%%%%%%%%%%%%%%%%%%%%%%%%%%%%%%%%%%%%%%%%%%%%%%%%%%%%%%%%%%%%%
\begin{proof}
Using \eqref{A.17}, \eqref{A.18}, \eqref{A.28}, and
\eqref{A.48} one finds
\begin{align}
\binom{p_+(z,k_0-1,k_0)}{r_+(z,k_0-1,k_0)} &=\begin{cases}
\frac{1}{\rho_{k_0}}\begin{pmatrix} z\ol{b_{k_0}}\\b_{k_0}\end{pmatrix}, &
\text{$k_0$ odd,} \\[4mm]
\frac{1}{\rho_{k_0}}\begin{pmatrix} b_{k_0}\\ \ol{b_{k_0}}\end{pmatrix}, &
\text{$k_0$ even,} \end{cases} \label{3.35}  \\
\binom{q_+(z,k_0-1,k_0)}{s_+(z,k_0-1,k_0)} &= \begin{cases}
\frac{1}{\rho_{k_0}}\begin{pmatrix} -z\ol{a_{k_0}} \\ a_{k_0}
\end{pmatrix}, & \text{$k_0$ odd,} \\[4mm]
\frac{1}{\rho_{k_0}}\begin{pmatrix} a_{k_0} \\
-\ol{a_{k_0}}\end{pmatrix}, & \text{$k_0$ even.} \end{cases}
\label{3.36}
\end{align}
It follows from \eqref{A.88} that
\begin{align}
M_{\ell,\ell'}(z,k_0) &= \de_{\ell,\ell'} +
2z\big(\de_{k_0+\ell-1},(U-zI)^{-1}\de_{k_0+\ell'-1}\big)_{\ell^2(\Z)}
\no \\
&= \de_{\ell,\ell'} + 2 z
(U-zI)^{-1}(k_0+\ell-1,k_0+\ell'-1), \quad \ell,\ell'=0,1.
\end{align}
Thus, by Lemma \ref{lA.18} and equalities
\eqref{A.48}, \eqref{A.68}, \eqref{3.35}, and \eqref{3.36},
one finds
\begin{align}
(U-zI)^{-1}(k_0,k_0) &=
\frac{[1-M_+(z,k_0)][1+M_-(z,k_0)]}{2z[M_+(z,k_0)-M_-(z,k_0)]},
\\
(U-zI)^{-1}(k_0-1,k_0-1) &=
\f{[\ol{a_{k_0}}-\ol{b_{k_0}}M_+(z,k_0)][a_{k_0}
+b_{k_0}M_-(z,k_0)]}{2z\rho_{k_0}^2[M_+(z,k_0)-M_-(z,k_0)]},
\\
(U-zI)^{-1}(k_0-1,k_0) &= - \f{
\begin{cases}
{[1-M_+(z,k_0)][\ol{a_{k_0}}-\ol{b_{k_0}}M_-(z,k_0)]}, &
\text{$k_0$ odd,}
\\
{[1+M_+(z,k_0)][a_{k_0}+b_{k_0}M_-(z,k_0)]}, & \text{$k_0$
even,}
\end{cases}}{2z\rho_{k_0}[M_+(z,k_0)-M_-(z,k_0)]},
\\
(U-zI)^{-1}(k_0,k_0-1) &= - \f{
\begin{cases}
{[1+M_+(z,k_0)][a_{k_0}+b_{k_0}M_-(z,k_0)]}, & \text{$k_0$
odd,}
\\
{[1-M_+(z,k_0)][\ol{a_{k_0}}-\ol{b_{k_0}}M_-(z,k_0)]}, &
\text{$k_0$ even,}
\end{cases}}{2z\rho_{k_0}[M_+(z,k_0)-M_-(z,k_0)]},
\end{align}
and hence \eqref{A.97a}--\eqref{A.97c}.
\end{proof}
%%%%%%%%%%%%%%%%%%%%%%%%%%%%%%%%%%%%%%%%%%%%%%%%%%%%%%%%%%%%%%%

Finally, introducing the functions $\Phi_{1,1}(\cdot,k)$,
$k\in\bbZ$, by
\begin{equation} \label{A.94}
\Phi_{1,1}(z,k) = \frac{M_{1,1}(z,k)-1}{M_{1,1}(z,k)+1}, \quad
z\in\C\backslash\dD,
\end{equation}
then,
\begin{equation}
M_{1,1}(z,k) = \frac{1+\Phi_{1,1}(z,k)}{1-\Phi_{1,1}(z,k)}, \quad
z\in\C\backslash\dD. \lb{A.95}
\end{equation}
Both, $M_{1,1}(z,k)$ and $\Phi_{1,1}(z,k)$, $z\in\bbC\backslash\dD$,
$k\in\Z$, have nontangential limits to $\dD$ $\mu_0$-a.e.

%%%%%%%%%%%%%%%%%%%%%%%%%%%%%%%%%%%%%%%%%%%%%%%%%%%%%%%%%%%%%%%
\begin{lemma}
The function $\Phi_{1,1}|_{\D}$ is a Schur function and $\Phi_{1,1}$
is related to $\Phi_\pm$ by
\begin{equation}
\Phi_{1,1}(z,k) = \frac{\Phi_+(z,k)}{\Phi_-(z,k)}, \quad
z\in\C\backslash\dD, \; k\in\bbZ. \lb{A.96}
\end{equation}
\end{lemma}
%%%%%%%%%%%%%%%%%%%%%%%%%%%%%%%%%%%%%%%%%%%%%%%%%%%%%%%%%%%%%%%
\begin{proof}
The assertion follows from \eqref{A.78}, \eqref{A.94} and
Lemma \ref{lA.20}.
\end{proof}
%%%%%%%%%%%%%%%%%%%%%%%%%%%%%%%%%%%%%%%%%%%%%%%%%%%%%%%%%%%%%%%

%%%%%%%%%%%%%%%%%%%%%%%%%%%%%%%%%%%%%%%%%%%%%%%%%%%%%%%%%%%%%%%
\begin{lemma} \lb{l3.4}
Let $\zeta\in\dD$ and $k_0\in\Z$. Then the following sets of
two-dimensional Laurent polynomials $\{P(\ze,k,k_0)\}_{k\in\Z}$ and
$\{R(\ze,k,k_0)\}_{k\in\Z}$,
\begin{align}
P(\ze,k,k_0) &= \binom{P_0(\ze,k,k_0)}{P_1(\ze,k,k_0)} =
\begin{cases}
\frac{1}{2\ze}
\begin{pmatrix} -\rho_{k_0} & \rho_{k_0}
\\ \ol{b_{k_0}} & \ol{a_{k_0}} \end{pmatrix}
\begin{pmatrix}q_+(\ze,k,k_0) \\ p_+(\ze,k,k_0)\end{pmatrix},
& k_0 \text{ odd},
\\[4mm]
\frac{1}{2}
\begin{pmatrix} \rho_{k_0} & \rho_{k_0}
\\ -b_{k_0} & a_{k_0} \end{pmatrix}
\begin{pmatrix}q_+(\ze,k,k_0) \\ p_+(\ze,k,k_0)\end{pmatrix},
& k_0 \text{ even},
\end{cases}
\\
R(\ze,k,k_0) &= \binom{R_0(\ze,k,k_0)}{R_1(\ze,k,k_0)} =
\begin{cases}
\frac{1}{2} \begin{pmatrix} \rho_{k_0} & \rho_{k_0}
\\ -b_{k_0} & a_{k_0} \end{pmatrix}
\begin{pmatrix} s_+(\ze,k,k_0) \\
r_+(\ze,k,k_0)\end{pmatrix}, & k_0 \text{ odd},
\\[4mm]
\frac{1}{2} \begin{pmatrix} -\rho_{k_0} & \rho_{k_0}
\\ \ol{b_{k_0}} & \ol{a_{k_0}} \end{pmatrix}
\begin{pmatrix} s_+(\ze,k,k_0) \\
r_+(\ze,k,k_0)\end{pmatrix}, & k_0 \text{ even}
\end{cases}
\end{align}
form complete orthonormal systems in
$L^2(\dD;d\Om(\cdot,k_0)^\top)$ and
$L^2(\dD;d\Om(\cdot,k_0))$, respectively.
\end{lemma}
%%%%%%%%%%%%%%%%%%%%%%%%%%%%%%%%%%%%%%%%%%%%%%%%%%%%%%%%%%%%%%%
\begin{proof}
Consider the following relation,
\begin{align} \label{A.189}
U^\top\de_k = \sum_{j\in\bbZ} U^\top(j,k)\de_j =
\sum_{j\in\bbZ} U(k,j)\de_j, \quad k\in\Z.
\end{align}
By Lemma \ref{lA.1} any solution $u$ of
\begin{align} \label{A.191}
z u(z,k,k_0) = \sum_{j\in\bbZ} U(k,j) u(z,j,k_0), \quad
k\in\bbZ,
\end{align}
is a linear combination of $p_+(z,\cdot,k_0)$ and
$q_+(z,\cdot,k_0)$, and hence, \eqref{A.191} has a unique
solution $\{u(z,k,k_0)\}_{k\in\bbZ}$ with prescribed values
at $k_0-1$ and $k_0$,
\begin{align} \label{A.192}
u&(z,\cdot,k_0) = P_0(z,\cdot,k_0)u(z,k_0-1,k_0) +
P_1(z,\cdot,k_0)u(z,k_0,k_0).
\end{align}
Due to the algebraic nature of the proof of Lemma \ref{lA.1}
and the algebraic similarity of equations \eqref{A.189} and
\eqref{A.191}, one concludes from \eqref{A.192} that
\begin{align} \label{A.193}
\de_k = P_0(U^\top,k,k_0)\de_{k_0-1} +
P_1(U^\top,k,k_0)\de_{k_0}, \quad k\in\Z.
\end{align}
Using the spectral representation for the operator $U^\top$
one then obtains
\begin{align}
P_\ell(U^\top,k,k_0) = \oint_\dD dE_{U^\top}(\ze) \,
P_\ell(\ze,k,k_0), \quad \ell=0,1
\end{align}
and by \eqref{A.193},
\begin{align}
(\de_k,\de_{k'})_{\ell^2(\Z)} &= \sum_{\ell,\ell'=0}^1
\Big(P_\ell(U^\top,k,k_0)\de_{k_0+\ell-1},
P_{\ell'}(U^\top,k',k_0)\de_{k_0+\ell'-1}\Big)_{\ell^2(\Z)}
\no
\\ &= \oint_\dD P(\ze,k,k_0)^*\,d\Omega(\ze,k_0)^\top
P(\ze,k',k_0).
\end{align}

Similarly, one obtains the orthonormality relation for the
two-dimensional Laurent polynomials
$\{R(\ze,k,k_0)\}_{k\in\Z}$ in $L^2(\dD;d\Om(\cdot,k_0))$.

To prove completeness of $\{P(\ze,k,k_0)\}_{k\in\Z}$ we
first note the following fact,
\begin{align}
\spn\{P(z,k,k_0)\}_{k\in\bbZ} &= \spn \left\{
\begin{pmatrix} z^k \\ z^{k-1}\end{pmatrix},
\begin{pmatrix} z^{k-1} \\ z^k\end{pmatrix},
\begin{pmatrix} 1 \\ 0\end{pmatrix},
\begin{pmatrix} 0 \\ 1\end{pmatrix} \right\}_{k\in\Z} \no \\
&= \spn\left\{\begin{pmatrix} z^k \\ 0 \end{pmatrix},
\begin{pmatrix} 0 \\ z^k\end{pmatrix} \right\}_{k\in\bbZ},
\quad k_0\in\Z.
\end{align}
This follows by investigating the leading coefficients of
$p_+(z,k,k_0)$ and $q_+(z,k,k_0)$. Thus, it suffices to
prove that $\Big\{\Big(\begin{smallmatrix}\ze^k \\
0\end{smallmatrix}\Big), \Big(\begin{smallmatrix}0 \\ \ze^k
\end{smallmatrix}\Big)\Big\}_{k\in\Z}$ form a basis in
$L^2(\dD;d\Omega(\cdot,k_0)^\top)$ for all
$k_0\in\Z$.

Let $k_0\in\Z$ and suppose that $F=\Big(\begin{smallmatrix}
f_0 \\ f_1 \end{smallmatrix}\Big)\in
L^2(\dD;d\Omega(\cdot,k_0)^\top)$ is orthogonal to
$\Big\{\Big(\begin{smallmatrix} \ze^k \\ 0
\end{smallmatrix}\Big),\Big(\begin{smallmatrix} 0 \\
\ze^k \end{smallmatrix}\Big)\Big\}_{k\in\Z}$ in
$L^2(\dD;d\Omega(\cdot,k_0)^\top)$, that is,
\begin{equation}
0=\oint_{\dD} \ol{\begin{pmatrix} \zeta^k & 0
\end{pmatrix}}\, d\Omega (\zeta,k_0)^\top F(\zeta)  =\oint_{\dD}
\ol{\zeta^k} \, [f_0(\zeta) d\Omega_{0,0}(\zeta,k_0)
+f_1(\zeta)d\Omega_{1,0}(\zeta,k_0)]
\end{equation}
and
\begin{equation}
0=\oint_{\dD} \ol{\begin{pmatrix} 0 & \zeta^k
\end{pmatrix}}\, d\Omega (\zeta,k_0)^\top F(\zeta)  =\oint_{\dD}
\ol{\zeta^k} \, [f_0(\zeta) d\Omega_{0,1}(\zeta,k_0)
+f_1(\zeta)d\Omega_{1,1}(\zeta,k_0)]
\end{equation}
for all $k\in\Z$. Hence (cf., e.g., \cite[p.\ 24]{Du83}),
\begin{align}
f_0 d\Omega_{0,0}+f_1 d\Omega_{1,0} &= 0, \lb{3.50}
\\
f_0 d\Omega_{0,1}+f_1 d\Omega_{1,1} &= 0. \lb{3.52}
\end{align}
Multiplying \eqref{3.50} by $\ol{f_0}$ and \eqref{3.52} by
$\ol{f_1}$ then yields
\begin{equation}
\abs{f_0}^2 d\Omega_{0,0}+\ol{f_0}f_1 d\Omega_{1,0}
+\ol{f_1}f_0 d\Omega_{0,1}+\abs{f_1}^2 d\Omega_{1,1} = 0
\end{equation}
and hence
\begin{align}
\|F\|_{L^2(\dD;d\Omega(\cdot,k_0)^\top)}^2=\oint_{\dD} F(\zeta)^*\,
d\Omega(\zeta,k_0)^\top F(\zeta) = 0.
\end{align}

Similarly, one proves completeness of
$\{R(\ze,k,k_0)\}_{k\in\Z}$ in $L^2(\dD;d\Om(\cdot,k_0))$.
\end{proof}
%%%%%%%%%%%%%%%%%%%%%%%%%%%%%%%%%%%%%%%%%%%%%%%%%%%%%%%%%%%%%%%

Denoting by $I_2$ the identity operator in $\bbC^2$, we state the
following result.

%%%%%%%%%%%%%%%%%%%%%%%%%%%%%%%%%%%%%%%%%%%%%%%%%%%%%%%%%%%%%%%
\begin{corollary}
Let $k_0\in\bbZ$. Then the operators $U$ and $U^\top$ are unitarily
equivalent to the operator of multiplication by $I_2 id$ $($where
$id(\zeta)=\zeta$, $\ze\in\dD$$)$ on
$L^2(\dD;d\Om(\cdot,k_0))$ and $L^2(\dD;d\Om(\cdot,k_0)^\top)$,
respectively. Thus,
\begin{align}
\si(U) = \supp \, (d\Omega(\cdot,k_0)) =
\supp \, (d\Omega^{\rm tr}(\cdot,k_0)) =
\supp \, (d\Omega(\cdot,k_0)^\top) = \si(U^\top),
\end{align}
where
\begin{equation}
d\Omega^{\rm tr}(\cdot,k_0) = d\Omega_{0,0}(\cdot,k_0) +
d\Omega_{1,1}(\cdot,k_0)
\end{equation}
denotes the trace measure of $d\Omega(\cdot,k_0)$.
\end{corollary}
%%%%%%%%%%%%%%%%%%%%%%%%%%%%%%%%%%%%%%%%%%%%%%%%%%%%%%%%%%%%%%%
\begin{proof}
Consider the linear map $\dot \cU$ from $\ell^\infty_0(\bbZ)$ into the set
of two-dimensional Laurent polynomials on $\dD$ defined by,
\begin{equation}
(\dot \cU f)(\ze) = \sum_{k\in\bbZ} R(\ze,k,k_0)f(k), \quad f\in
\ell^\infty_0(\bbZ).
\end{equation}
A simple calculation for $F(\ze) = (\dot \cU f)(\ze)$, $f\in
\ell^\infty_0(\bbZ)$, shows that
\begin{align}
\sum_{k\in\bbZ} \abs{f(k)}^2 = \oint_\dD
F(\ze)^*d\Omega(\ze,k_0)F(\ze).
\end{align}
Since $\ell^\infty_0(\bbZ)$ is dense in $\ell^2(\bbZ)$,
$\dot \cU$ extends to a bounded linear operator $\cU \colon \ltz
\to L^2(\dD;d\Omega(\cdot,k_0))$. By Lemma \ref{l3.4},
$\cU$ is onto and one verifies that
\begin{align}
(\cU^{-1}F)(k) &= \oint_\dD R(\ze,k,k_0)^* d\Omega(\ze,k_0)
F(\ze).
\end{align}
In particular, $\cU$ is unitary. Moreover, we claim that
$\cU$ maps the operator $U$ on $\ell^2(\bbZ)$ to the
operator of multiplication by $id(\ze)=\ze$, $\ze\in\dD$,
denoted by $M(id)$, on $L^2(\dD; d\Omega(\cdot,k_0))$,
\begin{align}
\cU U \cU^{-1} = M(id),
\end{align}
where
\begin{align}
(M(id) F)(\ze) = \ze F(\ze), \quad F \in
L^2(\dD;d\Omega(\cdot,k_0)).
\end{align}
Indeed,
\begin{align}
&(\cU U \cU^{-1} F(\cdot))(\ze) = (\cU U f(\cdot))(\ze)  \no \\
&\quad = \sum_{k\in\Z} (U f(\cdot))(k) R(\ze,k,k_0) =
\sum_{k\in\Z} (U^\top R(\ze,\cdot,k_0))(k) f(k)  \no  \\
&\quad = \sum_{k\in\Z} \ze R(\ze,k,k_0) f(k) = \ze F(\ze) \no \\
& \quad = (M(id) F(\cdot))(\ze), \quad F \in L^2(\dD;d\Omega(\cdot,k_0)).
\end{align}
The result for the operator $U^\top$ is proved analogously.
\end{proof}
%%%%%%%%%%%%%%%%%%%%%%%%%%%%%%%%%%%%%%%%%%%%%%%%%%%%%%%%%%%%%%%

Finally, we note an alternative approach to (a variant of) the $2\times
2$ matrix-valued spectral function $\Omega(\cdot,k_0)$ associated with
$U$.

First we introduce $\wti \cM(z,k)$,
$z\in\bbC\backslash\dD$, $k\in\bbZ$, defined by
\begin{align}
\wti \cM(z,k)&= \begin{pmatrix}
\wti M_{0,0}(z,k) & \wti M_{0,1}(z,k) \\
\wti M_{1,0}(z,k) & \wti M_{1,1}(z,k) \end{pmatrix} \no \\
&=\begin{cases}
\f{1}{4}\begin{pmatrix} \rho_k & \rho_k \\
-b_k & a_k \end{pmatrix}^* \cM(z,k) \begin{pmatrix} \rho_k & \rho_k
\\[1mm]
-b_k & a_k \end{pmatrix}, & \text{$k$ odd}, \\
\f{1}{4}\begin{pmatrix} -\rho_k & \rho_k \\
\ol{b_k} & \ol{a_k} \end{pmatrix}^* \cM(z,k) \begin{pmatrix} -\rho_k &
\rho_k \\
\ol{b_k} & \ol{a_k} \end{pmatrix}, & \text{$k$ even},
\end{cases}  \no \\
&=\begin{pmatrix}
\f{1}{M_+(z,k)-M_-(z,k)}+\f{i}{2}\Im(\alpha_k) &
\f{1}{2}\f{M_+(z,k)+M_-(z,k)}{M_+(z,k)-M_-(z,k)}+\f{1}{2}\Re(\alpha_k)
\\[2mm]
-\f{1}{2}\f{M_+(z,k)+M_-(z,k)}{M_+(z,k)-M_-(z,k)}-\f{1}{2}\Re(\alpha_k) &
-\f{M_+(z,k)M_-(z,k)}{M_+(z,k)-M_-(z,k)}-\f{i}{2}\Im(\alpha_k)
\end{pmatrix} \no \\
& \hspace*{7cm} z\in\bbC\backslash\dD,  \; k\in\bbZ.
\end{align}
Clearly, $\cM(\cdot,k)$, and hence, $\wti \cM(\cdot,k)$,
$k\in\bbZ$, are $2\times 2$ matrix-valued Caratheodory
functions. Since by \eqref{A.88} $\cM(0,k)=I$, $k\in\Z$,
one computes
\begin{equation}
\wti \cM(0,k)=\f{1}{4}\begin{pmatrix}
\rho_k^2 + \abs{b_k}^2 & -2i\Im(\al_k)
\\ 2i\Im(\al_k) & \rho_k^2 + \abs{a_k}^2
\end{pmatrix}= [\wti \cM(0,k)]^*, \quad k\in\bbZ.
\end{equation}
Hence, the Herglotz representation of $\wti\cM(\cdot,k)$ is given by
\begin{align}
\wti\cM(z,k)=\int_{\dD} d\wti\Omega(\zeta,k) \,
\f{\zeta+z}{\zeta-z}, \quad  z\in\bbC\backslash\dD, \;
k\in\bbZ,
\end{align}
where the measure $d\wti\Omega(\cdot,k)$ can be reconstructed from the
boundary values of $\Re(\wti\cM(\cdot,k))$ via
\begin{align}
&\wti\Omega\big(\big(e^{i\theta_1},e^{i\theta_2}\big],k\big)=
\lim_{\delta\downarrow 0}\lim_{r\uparrow 1} \f{1}{2\pi}
\int_{\theta_1+\delta}^{\theta_2+\delta} d\theta \no \\ &
\quad \times \begin{pmatrix}
\Re\left(\f{1}{M_+(re^{i\theta},k)-M_-(re^{i\theta},k)}\right)
&
\f{i}{2}\Im\left(\f{M_+(re^{i\theta},k)+M_-(re^{i\theta},k)}
{M_+(re^{i\theta},k) -M_-(re^{i\theta},k)}\right) \\[3mm]
-\f{i}{2}\Im\left(\f{M_+(re^{i\theta},k)+M_-(re^{i\theta},k)}
{M_+(re^{i\theta},k) -M_-(re^{i\theta},k)}\right) &
-\Re\left(\f{M_+(re^{i\theta},k)M_-(re^{i\theta},k)}{M_+(re^{i\theta},k)
-M_-(re^{i\theta},k)}\right)
\end{pmatrix}, \\[1mm]
& \hspace*{4.9cm} \theta_1 \in [0,2\pi), \;
\theta_1<\theta_2<\theta_1+2\pi, \; k\in\bbZ. \no
\end{align}
Finally, the analog of Lemma \ref{l2.24} in the full-lattice context
reads as follows.

%%%%%%%%%%%%%%%%%%%%%%%%%%%%%%%%%%%%%%%%%%%%%%%%%%%%%%%%%%%%%%%%%%%%
\begin{lemma} \lb{l3.6}
Let $f,g \in\ell^\infty_0(\Z)$, $F\in C(\dD)$, and
$\theta_1\in [0,2\pi)$, $\theta_1<\theta_2\leq \theta_1+2\pi$. Then,
\begin{align}
\begin{split}
& \big(f,F(U)E_{U}
\big(\Arc\big(\big(e^{i\theta_1},e^{i\theta_2}\big]\big)\big)
g\big)_{\ell^2(\Z)}  \\ & \quad =\big(\hatt
f(\cdot,k_0),M_F
M_{\chi_{\Arc((e^{i\theta_1},e^{i\theta_2}])}} \hatt
g(\cdot,k_0)\big)_{L^2(\dD;d\wti\Omega_\pm (\cdot,k_0))},
\lb{3.65}
\end{split}
\end{align}
where we introduced the notation
\begin{equation}
\hatt h(\zeta,k_0)=\sum_{k\in\bbZ}
\binom{s_+(\zeta,k,k_0)}{r_+(\zeta,k,k_0)} h(k), \quad
\zeta\in\dD, \; h\in \ell^\infty_0(\Z), \lb{3.66}
\end{equation}
and $M_G$ denotes the maximally defined operator of multiplication by the
$d\wti\Omega(\cdot,k_0)$-measurable function $G$ in the
Hilbert space $L^2(\dD;d\wti\Omega(\cdot,k_0))$,
\begin{align}
\begin{split}
& (M_G\hatt h)(\zeta)=G(\zeta)\hatt h(\zeta) \, \text{ for
a.e.\ $\zeta\in\dD$},  \\ & \hatt h \in \dom(M_G)=\{\hatt
k\in L^2(\dD;d\wti\Omega(\cdot,k_0)) \,|\, G\hatt k \in
L^2(\dD;d\wti\Omega(\cdot,k_0))\}  \lb{3.67}
\end{split}
\end{align}
\end{lemma}
%%%%%%%%%%%%%%%%%%%%%%%%%%%%%%%%%%%%%%%%%%%%%%%%%%%%%%%%%%%%%%%%%%%%

Using Lemma \ref{l2.23}, \eqref{2.59}, \eqref{2.60},
\eqref{2.167}, and \eqref{A.84a} one can follow the proof of
Lemma \ref{l2.24} step by step and so we omit the details
(cf.\ also \cite{GZ05a}).

\medskip

Finally, Weyl--Titchmarsh circles associated with finite intervals
$[k_-,k_+]\cap\bbZ$ and the ensuing limits $k_\pm\to\pm\infty$ can be
discussed in analogy to the half-lattice case at the end of Section
\ref{s2}. Without entering into details, we mention that $U$ is of course
in the limit point case at $\pm\infty$.

%%%%%%%%%%%%%%%%%%%%%%%%%%%%%%%%%%%%%%%%%%%%%%%%%%%%%%%%%
%%%%%%%%%%% appendices %%%%%%%%%%%%%%%%%%%%%%%%%%%%%%%%%%
\appendix
%%%%%%%%%%% appendix A
\section{Basic Facts on Caratheodory and Schur Functions}
\lb{A}
\renewcommand{\theequation}{A.\arabic{equation}}
\renewcommand{\thetheorem}{A.\arabic{theorem}}
\setcounter{theorem}{0}
\setcounter{equation}{0}
%%%%%%%%%%%%%%%%%%%%%%%%%%%%%%%%%%%%%%%%%%%%%%%%%%%%%%%%%
%%%%%%%%%%%%%%%%%%%%%%%%%%%%%%%%%%%%%%%%%%%%%%%%%%%%%%%%%

In this appendix we summarize a few basic properties of Caratheodory
and Schur functions used throughout this manuscript.

We denote by $\D$ and $\dD$ the open unit disk
and the counterclockwise oriented unit circle in the complex plane $\C$,
\begin{equation}
\D = \{ z\in\C \st \abs{z} < 1 \}, \quad \dD = \{ \ze\in\C
\st \abs{\ze} = 1 \},
\end{equation}
and by
\begin{equation}
\Cl = \{z\in\C \st \Re(z) < 0\}, \quad \Cr = \{z\in\C \st
\Re(z) > 0\}
\end{equation}
the open left and right complex half-planes, respectively.

%%%%%%%%%%%%%%%%%%%%%%%%%%%%%%%%%%%%%%%%%%%%%%%%%%%%%%%%%%%%%%%
\begin{definition} \lb{dA.9}
Let $f_\pm$, $\varphi_+$, and $1/\varphi_-$ be analytic in
$\D$. \\ $(i)$ $f_+$ is called a {\it Caratheodory
function}
if $f_+\colon \D\to\Cr$ and $f_-$ is called an {\it
anti-Caratheodory function} if $-f_-$ is a Caratheodory
function. \\ $(ii)$ $\varphi_+$ is called a {\it Schur
function} if $\varphi_+\colon\D\to\D$. $\varphi_-$ is
called
an {\it anti-Schur function} if $1/\varphi_-$ is a Schur
function.
\end{definition}
%%%%%%%%%%%%%%%%%%%%%%%%%%%%%%%%%%%%%%%%%%%%%%%%%%%%%%%%%%%%%%%

%%%%%%%%%%%%%%%%%%%%%%%%%%%%%%%%%%%%%%%%%%%%%%%%%%%%%%%%%%%%%%%
\begin{theorem}[\cite{Ak65}, Sect.\ 3.1; \cite{AG81},
Sect.\ 69; \cite{Si04}, Sect.\ 1.3] \label{tA.2} ${}$ \\
Let $f$ be a Caratheodory function. Then $f$ admits the
Herglotz representation
\begin{align}
& f(z)=ic+ \oint_{\dD} d\mu(\zeta) \, \f{\zeta+z}{\zeta-z},
\quad z\in\D, \lb{A.3}
\\
& c=\Im(f(0)), \quad \oint_{\dD}
d\mu(\zeta) = \Re(f(0)) < \infty, \lb{A.4}
\end{align}
where $d\mu$ denotes a nonnegative measure on $\dD$. The
measure $d\mu$ can be reconstructed from $f$ by the formula
\begin{equation}
\mu\big(\Arc\big(\big(e^{i\theta_1},e^{i\theta_2}\big]\big)\big)
=\lim_{\delta\downarrow 0}
\lim_{r\uparrow 1} \f{1}{2\pi}
\oint_{\theta_1+\delta}^{\theta_2+\delta} d\theta \,
\Re\big(f\big(re^{i\theta}\big)\big), \lb{A.4a}
\end{equation}
where
\begin{equation}
\Arc\big(\big(e^{i\theta_1},e^{i\theta_2}\big]\big)
=\big\{e^{i\theta}\in\dD\,|\, \theta_1<\theta\leq
\theta_2\big\}, \quad \theta_1 \in [0,2\pi), \;
\theta_1<\theta_2\leq \theta_1+2\pi. \label{A.5}
\end{equation}
Conversely, the right-hand side of \eqref{A.3} with
$c\in\bbR$ and $d\mu$ a finite $($nonnegative$)$ measure on
$\dD$ defines a Caratheodory function.
\end{theorem}
%%%%%%%%%%%%%%%%%%%%%%%%%%%%%%%%%%%%%%%%%%%%%%%%%%%%%%%%%%%%%%%

We note that additive nonnegative constants on the
right-hand
side of \eqref{A.3} can be absorbed into the measure $d\mu$
since
\begin{equation}
\oint_\dD d\mu_0(\zeta) \, \f{\zeta+z}{\zeta-z}=1, \quad z\in\D,
\lb{A.5a}
\end{equation}
where
\begin{equation}
d\mu_0(\zeta)=\f{d\theta}{2\pi}, \quad
\zeta=e^{i\theta}, \;
\theta\in [0,2\pi] \lb{A.5b}
\end{equation}
denotes the normalized Lebesgue measure on the unit circle
$\dD$.

A useful fact on Caratheodory functions $f$ is a certain
monotonicity property they exhibit on open connected arcs
of the unit circle away from the support of the measure
$d\mu$ in the Herglotz representation \eqref{A.3}. More
precisely, suppose
$\Arc\big(\big(e^{i\theta_1},e^{i\theta_2}\big)\big)\subset
(\dD\backslash\supp(d\mu))$, $\theta_1<\theta_2$, then $f$
has an analytic continuation through
$\Arc\big(\big(e^{i\theta_1},e^{i\theta_2}\big)\big)$ and
it is purely imaginary on
$\Arc\big(\big(e^{i\theta_1},e^{i\theta_2}\big)\big)$.
Moreover,
\begin{equation}
\f{d}{d\theta}f\big(e^{i\theta}\big)=-\f{i}{2}
\int_{[0,2\pi]\backslash(\theta_1,\theta_2)}
d\mu\big(e^{it}\big) \f{1}{\sin^2((t-\theta)/2)}, \quad
\theta\in(\theta_1,\theta_2).
\end{equation}
In particular,
\begin{equation}
-i\f{d}{d\theta}f\big(e^{i\theta}\big)<0, \quad
\theta\in(\theta_1,\theta_2).
\end{equation}

We recall that any Caratheodory function $f$ has finite radial limits to
the unit circle $\mu_0$-almost everywhere, that is,
\begin{equation}
f(\zeta)=\lim_{r\uparrow 1} f(r\zeta) \, \text{ exists and is finite for
$\mu_0$-a.e.\ $\zeta\in\dD$.}
\end{equation}

The absolutely continuous part $d\mu_{\rm ac}$ of the measure
$d\mu$ in the Herglotz representation \eqref{A.3} of the Caratheodory
function $f$ is given by
\begin{equation}
d\mu_{\rm ac}(\zeta) = \lim_{r\uparrow 1}
\Re(f(r\zeta)) \, d\mu_0(\zeta), \quad \zeta \in\dD.
\lb{A.10}
\end{equation}
The set
\begin{equation}
S_{\mu_{\rm ac}}=\{\zeta\in\dD\,|\, \lim_{r\uparrow 1}\Re(f(r\zeta))
=\Re(f(\zeta))>0 \text{ exists finitely}\}  \lb{A.10A}
\end{equation}
is an essential support of $d\mu_{\rm ac}$ and its essential closure,
$\ol{S_{\mu_{\rm ac}}}^e$, coincides with the topological support,
$\supp(d\mu_{\rm ac})$ (the smallest closed support), of
$d\mu_{\rm ac}$,
\begin{equation}
\ol{S_{\mu_{\rm ac}}}^e=\supp \, (d\mu_{\rm ac}).  \lb{A.10B}
\end{equation}
Moreover, the set
\begin{align}
S_{\mu_{\rm s}} = \{\zeta\in\dD\,|\, \lim_{r\uparrow
1}\Re(f(r\zeta)) = \infty \} \label{A.10d}
\end{align}
is an essential support of the singular part $d\mu_{\rm s}$
of the measure $d\mu$, and
\begin{equation}
\lim_{r\uparrow 1} (1-r)f(r\zeta)=\lim_{r\uparrow 1}
(1-r)\Re(f(r\zeta))\geq 0 \, \text{ exists for all
$\zeta\in\dD$.} \lb{A.10C}
\end{equation}
In particular, $\zeta_0\in\dD$ is a pure point of $d\mu$ if
and only if
\begin{equation}
\mu(\{\zeta_0\})=\lim_{r \uparrow 1} \left(\f{1-r}{2}\right)
f(r\zeta_0)>0. \lb{A.10a}
\end{equation}

Given a Caratheodory (resp., anti-Caratheodory) function $f_+$
(resp. $f_-$) defined in $\D$ as in \eqref{A.3}, one extends $f_\pm$ to
all of $\bbC\backslash\dD$ by
\begin{equation}
f_\pm(z)=ic_\pm \pm \oint_{\dD} d\mu_\pm (\zeta) \,
\f{\zeta+z}{\zeta-z}, \quad z\in\bbC\backslash\dD, \;
c_\pm\in\bbR. \lb{A.6}
\end{equation}
In particular,
\begin{equation}
f_\pm(z) = -\ol{f_\pm(1/\ol{z})}, \quad
z\in\C\backslash\ol{\D}. \lb{A.7}
\end{equation}
Of course, this continuation of $f_\pm|_{\D}$ to
$\bbC\backslash\ol\D$, in general, is not an analytic
continuation of $f_\pm|_\D$. With $f_\pm$ defined on
$\bbC\backslash\dD$ by \eqref{A.6} one infers the mapping
properties
\begin{equation}
f_+\colon \D \to \Cr, \quad f_+\colon\bbC\backslash
\ol\D\to\Cl, \quad f_-\colon \D \to \Cl, \quad f_-\colon\bbC\backslash
\ol\D\to\Cr.
\end{equation}

Next, given the functions $f_\pm$ defined in $\bbC\backslash\dD$
as in \eqref{A.6}, we introduce the functions $\varphi_\pm$
by
\begin{equation}
\varphi_\pm(z)=\f{f_\pm(z)-1}{f_\pm(z)+1}, \quad
z\in\bbC\backslash\dD.  \lb{A.18a}
\end{equation}
Then $\varphi_\pm$ have the mapping properties
\begin{align}
\begin{split}
& \varphi_+\colon\D\to\D, \quad
1/\varphi_+\colon\bbC\backslash\ol\D \to\D \quad
(\varphi_+\colon\bbC\backslash\ol\D\to
(\bbC\backslash\ol\D)\cup\{\infty\}), \\ &
\varphi_-\colon\bbC\backslash\ol\D\to\D, \quad
1/\varphi_-\colon\D\to\D \quad (\varphi_-\colon\D\to
(\bbC\backslash\ol\D)\cup\{\infty\}),
\end{split}
\end{align}
in particular, $\varphi_+|_{\D}$ (resp., $\varphi_-|_{\D}$)
is a Schur (resp., anti-Schur) function. Moreover,
\begin{equation}
f_\pm(z)=\f{1+\varphi_\pm (z)}{1-\varphi_\pm (z)}, \quad
z\in\bbC\backslash\dD.
\end{equation}

We also recall the following useful result (see \cite[Lemma
10.11.17]{Si04} and \cite{GZ05} for a proof). To fix some notation
we denote  by $f_+$ and $f_-$ a Caratheodory and anti-Caratheodory
function, respectively, and by $\varphi_+$ and
$\varphi_-$ the corresponding Schur and anti-Schur functions as defined
in \eqref{A.18a}. We also introduce the following notation for open arcs
on the unit circle $\dD$,
\begin{equation}
\Arc\big(\big(e^{i\theta_1},e^{i\theta_2}\big)\big)
=\big\{e^{i\theta}\in\dD\,|\,
\theta_1<\theta< \theta_2\big\}, \quad \theta_1 \in
[0,2\pi], \; \theta_1< \theta_2\leq \theta_1+2\pi. \lb{A.23A}
\end{equation}
An open arc $A\subseteq \dD$ then either coincides with
$\Arc\big(\big(e^{i\theta_1},e^{i\theta_2}\big)\big)$ for some
$\theta_1 \in [0,2\pi]$, $\theta_1< \theta_2\leq \theta_1+2\pi$, or
else, $A=\dD$.

%%%%%%%%%%%%%%%%%%%%%%%%%%%%%%%%%%%%%%%%%%%%%%%%%%%%%%%%
\begin{lemma} \lb{lA.3A}
Let $A\subseteq\dD$ be an open arc and assume that  $f_+$ $($resp.,
$f_-$$)$ is a Caratheodory $($resp., anti-Caratheodory$)$ function
satisfying the reflectionless condition
\begin{equation}
\lim_{r\uparrow 1}\big[f_+(r\zeta)+\ol{f_-(r\zeta)} \big]=0 \,
\text{ $\mu_0$-a.e.\ on $A$}.  \lb{A.24A}
\end{equation}
Then, \\
(i) $f_+(\zeta)=-\ol{f_-(\zeta)}$ for all $\zeta\in A$. \\
(ii) For $z\in\D$, $-\ol{f_-(1/\ol z)}$ is the analytic
continuation of $f_+(z)$ through the arc $A$. \\
(iii) $d\mu_\pm$ are purely absolutely continuous on $A$ and
\begin{equation}
\f{d\mu_{\pm}}{d\mu_0}(\zeta)=\Re(f_+(\zeta)) =-\Re(f_-(\zeta)), \quad
\zeta\in A. \lb{A.25A}
\end{equation}
\end{lemma}
%%%%%%%%%%%%%%%%%%%%%%%%%%%%%%%%%%%%%%%%%%%%%%%%%%%%%%%%

In analogy to the exponential representation of Nevanlinna--Herglotz
functions (i.e., functions analytic in the open complex upper
half-plane $\bbC_+$ with a strictly positive imaginary part
on $\bbC_+$, cf.\ \cite{AD56}, \cite{AD64}, \cite{GT00}, \cite{KK74}) one
obtains the following result.

%%%%%%%%%%%%%%%%%%%%%%%%%%%%%%%%%%%%%%%%%%%%%%%%%%%%%%%%%%%%%%%
\begin{theorem} \lb{tA.3}
Let $f$ be a Caratheodory function. Then $-i\ln(if)$ is a
Caratheodory function and $f$ has the exponential Herglotz
representation,
\begin{align}
& -i\ln(if(z))=id+ \oint_{\dD} d\mu_0(\zeta) \,\Upsilon
(\zeta) \, \f{\zeta+z}{\zeta-z}, \quad z\in\D,
\lb{A.11} \\ &
\;\, d=-\Re(\ln(f(0))), \quad 0 \leq \Upsilon
(\zeta)\leq \pi
\, \text{ for $\mu_0$-a.e.\ $\zeta\in\dD$}. \lb{A.13}
\end{align}
$\Upsilon$ can be reconstructed from $f$ by
\begin{align}
\begin{split}
\Upsilon (\zeta) &= \lim_{r\uparrow
1}\Re[-i\ln(if(r\zeta))]
                \\ & =(\pi/2)+\lim_{r\uparrow 1}\Im[\ln(f(r\zeta))]\,
\text{ for $\mu_0$-a.e.\ $\zeta\in\dD$.}
\end{split}
\end{align}
\end{theorem}
%%%%%%%%%%%%%%%%%%%%%%%%%%%%%%%%%%%%%%%%%%%%%%%%%%%%%%%%%%%%%%%

Next we briefly turn to matrix-valued Caratheodory
functions. We denote as usual $\Re(A)=(A+A^*)/2$,
$\Im(A)=(A-A^*)/(2i)$, etc., for square matrices $A$.

%%%%%%%%%%%%%%%%%%%%%%%%%%%%%%%%%%%%%%%%%%%%%%%%%%%%%%%%%%%%%%%
\begin{definition} \lb{dA.4}
Let $m\in\bbN$ and $\cF$ be an $m\times m$ matrix-valued function
analytic in $\D$. $\cF$ is called a {\it Caratheodory matrix} if
$\Re(\cF(z))\geq 0$ for all $z\in\D$.
\end{definition}
%%%%%%%%%%%%%%%%%%%%%%%%%%%%%%%%%%%%%%%%%%%%%%%%%%%%%%%%%%%%%%%

%%%%%%%%%%%%%%%%%%%%%%%%%%%%%%%%%%%%%%%%%%%%%%%%%%%%%%%%%%%%%%%
\begin{theorem} \label{tA.5}
Let $\cF$ be an $m\times m$ Caratheodory matrix,
$m\in\bbN$. Then $\cF$ admits the Herglotz representation
\begin{align}
& \cF(z)=iC+ \oint_{\dD} d\Omega(\zeta) \,
\f{\zeta+z}{\zeta-z}, \quad z\in\D, \lb{A.240}
\\
& C=\Im(\cF(0)), \quad \oint_{\dD} d\Omega(\zeta) =
\Re(\cF(0)), \lb{A.250}
\end{align}
where $d\Omega$ denotes a nonnegative $m \times m$
matrix-valued measure on $\dD$. The measure $d\Omega$ can
be reconstructed from $\cF$ by the formula
\begin{align}
&\Omega\big(\Arc\big(\big(e^{i\theta_1},e^{i\theta_2}\big]\big)\big)
=\lim_{\delta\downarrow 0} \lim_{r\uparrow 1} \f{1}{2\pi}
\oint_{\theta_1+\delta}^{\theta_2+\delta} d\theta \,
\Re\big(\cF\big(re^{i\theta}\big)\big), \\
& \hspace*{3.9cm} \theta_1 \in
[0,2\pi], \; \theta_1<\theta_2\leq \theta_1+2\pi. \no
\end{align}
Conversely, the right-hand side of equation \eqref{A.240}
with $C = C^*$ and $d\Omega$ a finite nonnegative $m \times
m$ matrix-valued measure on $\dD$ defines a Caratheodory
matrix.
\end{theorem}
%%%%%%%%%%%%%%%%%%%%%%%%%%%%%%%%%%%%%%%%%%%%%%%%%%%%%%%%%%%%%%%

%%%%%%%%%%%%%%%%%%%%%%%%%%%%%%%%%%%%%%%%%%%%%%%%%%%%%%%%%%%%%%%%%%%%%%%
{\bf Acknowledgments.}
We are indebted to Barry Simon for providing us with a copy of his
two-volume treatise \cite{Si04} prior to its publication. We also gratefully
acknowledge the extraordinarily generous and detailed comments by an
anonymous referee, which lead to numerous improvements in the presentation
of the material in this paper.
%%%%%%%%%%%%%%%%%%%%%%%%%%%%%%%%%%%%%%%%%%%%%%%%%%%%%%%%%%%%%%%%%%%%%%%

%%%%%%%%%%%%%%%%%%%%%%%%%%%%%%%%%%%%%%%%%%%%%%%%%%%%%%%%

\end{document}